\documentclass[12pt,a4paper]{amsart}
\usepackage{hyperref,amssymb,mathrsfs,amscd}

\let\g=\mathfrak

\newcommand{\OddId}{{{\scriptstyle \Pi}  }}
\def\O{\emptyset}
\newcommand{\0}{{\boldsymbol{0}}}
\newcommand{\1}{{\boldsymbol{1}}}
\newcommand{\ii}{{\boldsymbol{i}}}
\newcommand{\bb}{{\boldsymbol{b}}}
\newcommand{\bv}{{\boldsymbol{v}}}
\newcommand{\bd}{{\boldsymbol{d}}}
\newcommand{\bk}{{\boldsymbol{k}}}
\newcommand{\bm}{{\boldsymbol{m}}}

\newcommand{\Hom}{\operatorname{Hom}}
\newcommand{\Ad}{\operatorname{Ad}}

\newcommand{\Spf}{\operatorname{Spf}}

\newcommand{\Ber}{\operatorname{Ber}}
\newcommand{\tr}{\operatorname{tr}}
\newcommand{\str}{\operatorname{str}}
\newcommand{\rk}{\operatorname{rk}}

\def\AC{{\mathcal A}}

\def\CC{{\mathcal C}}
\def\DC{{\mathcal D}}
\def\EC{{\mathcal E}}

\def\IC{{\mathcal I}}
\def\JC{{\mathcal J}}

\def\LC{{\mathcal L}}

\def\OC{{\mathcal O}}
\def\PC{{\mathcal P}}

\def\SC{{\mathcal S}}

\def\UC{{\mathcal U}}
\def\VC{{\mathcal V}}
\def\WC{{\mathcal W}}

\let\a=\alpha
\let\be=\beta
\let\ga=\gamma
\let\C=\Gamma
\let\t=\tau

\let\l=\lambda
\let\L=\Lambda
\let\w=\omega
\let\W=\Omega
\let\r=\rho
\let\s=\sigma
\let\d=\delta

\let\de=
\partial
\let\th=\theta
\let\z=\zeta

\newtheorem{theo}{Theorem}[section]
\newtheorem{defi}{Definition}[section]

\newtheorem{prop}{Proposition}[section]
\newtheorem{lemme}{Lemma}[section]

\let\what=\widehat
\def\som{\mathop{\sum}\limits}

\def\tens{\mathop{\otimes}\limits}
\def\sdir{\mathop{\oplus}\limits}
\def\mult{\mathop{\times}\limits}
\def\sup{\mathop{Sup}\limits}

\let\oo=\infty
\def\der{\mathcal Der}

\def\vol{\mathcal V\text{\it{ol}}}
\def\ber{\mathcal B\text{\it{er}}}

\def\hfl#1#2{\smash{\mathop{\hbox to
12mm{\rightarrowfill}}\limits^{\scriptstyle#1}_{\scriptstyle#2}}}

\def\ihfl#1#2{\smash{\mathop{\hbox to
12mm{\hookrightarrowfill}}\limits^{\scriptstyle#1}_{\scriptstyle#2}}}

\def\lihfl#1#2{\smash{\mathop{\hbox to
12mm{\hookleftarrowfill}}\limits^{\scriptstyle#1}_{\scriptstyle#2}}}

\def\lhfl#1#2{\smash{\mathop{\hbox to
12mm{\leqslantftarrowfill}}\limits^{\scriptstyle#1}_{\scriptstyle#2}}}

\title[Equivariant Cohomology and Localization Formula]{Equivariant Cohomology
and Localization Formula in Supergeometry}
\author{P. Lavaud}
\date{\today}

\keywords{14M30, 17B70,  58A50, 58C50}

\begin{document}

\begin{abstract} 
Let $G$ be a compact Lie group. Let $M$ be a smooth $G$-manifold and $\VC\to M$ be an oriented $G$-equivariant vector bundle. One defines the spaces of equivariant forms with generalized coefficients on  $\VC$  and $M$. An equivariant Thom form $\th$ on $\VC$ is a compactly supported closed equivariant form such that its integral along the fibres is the constant function $1$ on $M$. Such a Thom form was constructed by Mathai and Quillen \cite{MQ86}. Its restriction to  $M$ gives a representative of the equivariant Euler class of $\VC$.

In the supergeometric situation we give proper definitions of all the objects involved. But, in this case a Thom form doesn't always exist. In this article, when the action of $G$ on $\VC$ is sufficiently non-trivial, we construct such a Thom form with generalized coefficients. We use it to construct an equivariant Euler form of $\VC$ and to generalize Berline-Vergne's localization formula  
(\cite{BV83a}) to the supergeometric situation. 

\end{abstract}

\maketitle

The aim of this article is to generalize Berline-Vergne's localization formula  
(\cite{BV83a}) to the supergeometric situation. 

Let $G=(G_{\0},\g g)$ be a Lie supergroup (with underlying Lie group $G_\0$ and
Lie superalgebra $\g g$) acting on the right on a globally oriented
supermanifold $M=(M_{\0},\OC_{M})$ (see section \ref{Ori} for precise
definitions). Let $\a$ be a closed equivariant integrable form (the superanalog
of closed equivariant forms with compact support). Let $X$ be an element of $\g
g_\0$ such that  $\overline{\exp(\mathbb R X)}$ is a compact subgroup of $G_\0$.
Let $j:M(X)\hookrightarrow M$ be the subsupermanifold of zeroes of $X$.  Under
some conditions, we  construct an  equivariant Euler form for $T_N M(X)$ (the
normal bundle of $M(X)$ in M). We denote it by $\EC_{\g g}$. Then we obtain, for
$Z$ in a neighborhood of $X$ in $\g g(X)$ (the centralizer of $X$ in $\g g$) the
following equality:
\begin{equation}
\int_{M}\a(Z)=(2\pi)^{\frac{n+m}2}\int_{M(X)}\frac{j^*\a(Z)}{ \EC_{\g g}(Z)}
\end{equation}

\bigskip

Let us describe the most remarkable condition denoted by $(**)$ (cf. Theorems
\ref{Thom} and \ref{Formula}) needed to construct an Euler form. 

We assume that there is a $G(X)$-invariant Euclidean structure (cf. section
\ref{SEucl}) $Q$ on $T_NM(X)$, and  that $T_NM(X)$ has an equivariant
superconnection $\mathbb A_{\g g}$ with equivariant moment $\mu_{\mathbb A}$
(cf. section \ref{Econn}).  We denote by $\C_{T_NM(X)}$ the sheaf of sections of
$T_NM(X)$. The condition $(**)$ is that there is a covering of $M_{\0}$ by open
subsets $W$ such that:

\begin{multline} U_{+}^{T_NM(X)_{\1}}(W)=\Big\{Z\in\g g_\0(X)\,\Big/\,\\
\big(\forall v\in\C_{T_NM(X)}(W)_\1/ 
\forall x\in W,  v_{\mathbb R}(x)\not=0\big),\ 
 Q(v,\mu_{\mathbb A}(Z)v)_{\mathbb R}>0\Big\}.
\end{multline} contains a non empty open subset.

($Q(v,\mu_{\mathbb A}(Z)v)$ is a function on $M$ and $Q(v,\mu_{\mathbb
A}(Z)v)_{\mathbb R}$ is its restriction to its underlying manifold $M_\0$.)

\bigskip

A. S. Schwarz and O. Zaboronky found in \cite{SZ97} a localization formula over
a    supermanifold  in presence of symmetries given by the actions of compact
groups. Here, we consider symmetries given by the actions of more general
groups. In particular, we can have an action of a Lie supergroup which is not a
Lie group.

\bigskip

One motivation for such a formula is to produce  invariant generalized functions
on a Lie superalgebra and to test the validity of a ``Kirillov-like'' character
formula for Lie supergroups. The above condition seems very natural when applied
to a coadjoint orbit of a supergroup $G$. When $\g g$ is nilpotent, this
condition
 allows to associate  an equivalence class of unitary representations of $G$
with the orbit as pointed out to me by M. Duflo. 

\bigskip

In order to obtain such a formula we need to construct a representative of the
equivariant Euler class of an equivariant supervector bundle $\VC\to M$. It is
obtained by restriction to $M_{\0}$ of an equivariant Thom form on $\VC$ (cf.
section
\ref{Eul.sec}).

\bigskip

Let us say some words about the non-equivariant case (cf.
\cite{BS84,VZ88,Vor91,Lav98}). Let $\pi:\VC\rightarrow M$ be a supervector
bundle. We assume that $\VC$ and $M$ are globally oriented. Let
$\what\W_{\int}(\VC)$ (resp. $\what\W_{\int}(M)$) be the space integrable
pseudodifferential forms  on $\VC$ (resp. $M$) (cf. sections \ref{Pseudo} and
\ref{Int}).   As in the purely even situation, there is a ``Thom isomorphism''
between the cohomology of $\what\W_{\int}(\VC)$ and the cohomology of
$\what\W_{\int}(M)$. But unlike the even situation, it is not given by an
integration along the fibres. Indeed, there is no closed form on $\VC$ which is
integrable along the fibres and such that $\pi_*\a=1_M$ (this is the constant
function equal to $1$ on $M$ and $\pi_{*}$ is the integration along the fibres
of $\VC\to M$ cf. section \ref{ImDir}).

Since we work with smooth supermanifolds (cf. \cite{Bat79}), the supervector
bundle $\VC$ is the direct sum of its even part $\VC^\0$ and its odd part
$\VC^\1$. Therefore we have a bundle $\VC\rightarrow \VC^\0\rightarrow M.$ We
denote by $j^\0:\VC^\0\hookrightarrow \VC$ the embedding of $\VC^\0$ in $\VC$ by
mean of the zero-section and  by $\pi^\0:\VC^\0\rightarrow M$ the projection of
$\VC^\0$ onto $M$. Let $\a\in\what\W_{\int}(\VC)$. We take its restriction
$j^{\0*}(\a)$ to $\VC^\0$ and then we take its integral
$\pi^\0_*\big(j^{\0*}(\a)\big)$ along the fibres of $\VC^\0\rightarrow M$. The
application $\a\mapsto
\pi^\0_*\big(j^{\0*}(\a)\big)$ induces the ``Thom isomorphism'' between the
cohomologies of $\what\W_{\int}(\VC)$ and $\what\W_{\int}(M)$. The reverse
isomorphism sends a form $\be\in\what\W_{\int}(M)$ on $\th(\pi^*\be)$ where
$\th$ is a closed form on $\VC$ which is integrable along the fibres (cf.
\ref{ImDir}) and such that $\pi^\0_*\big(j^{\0*}(\th)\big)=1_M$.

\bigskip

Now, let $G$ be a supergroup. Let $\pi:\VC\rightarrow M$ be a $G$-equivariant
supervector bundle. Under condition $(**)$ and other technical ones (cf. Theorem
\ref{Thom}), we construct a closed equivariant form with generalized
coefficients $\th$ which is integrable along the fibres and such that
$\pi_*\th=1$ (cf. section \ref{Thom.sec}). We call it an equivariant Thom form
of $\VC$. Let $\what\W_{G,\int}^{-\oo}(\VC)$ be the space of integrable
equivariant forms with generalized coefficients on $\VC$. Because of the
generalized coefficients, $\what\W_{G,\int}^{-\oo}(\VC)$ is not an
$\what\W_{G,\int}^{-\oo}(M)$-module and the equivariant Thom form does not
provide a Thom isomorphism between the cohomologies of
$\what\W_{G,\int}^{-\oo}(\VC)$  and $\what\W_{G,\int}^{-\oo}(M)$. Nevertheless,
it provides an equivariant Euler form $\EC_{\g g}$ (cf. section \ref{Eul.sec})
and a relation between cohomology classes which is useful to obtain the
localization formula (cf. section \ref{RelCoho}). 

\bigskip

All notations are fixed in sections $1$ to $3$. These sections contain
essentially well known material (cf. for example \cite{Man88} for the most of
it). They are rather detailed and as a result the article is essentially self
contained. 

\bigskip

I wish to thank professor Michel Duflo who introduced me to supergeometry and 
helped me during this research.

\tableofcontents

\section{Prerequisites} We choose  a square root $\ii\in\mathbb C$ of $-1$.

\subsection{Supervector spaces} In this article, unless otherwise specified, all
supervector spaces and superalgebras will be real.  If $V$ is a supervector
space, we denote by $V_\0$ its even part and by $V_\1$ its odd part. If $v$ is a
non zero homogeneous element of $V$, we denote by $p(v)\in \mathbb Z/2\mathbb Z$
its parity.   We put $\dim(V)=(\dim(V_\0),\dim(V_\1))$.

Let $(m,n)\in \mathbb N\times \mathbb N$. We denote by $\mathbb R^{(m,n)}$ the
supervector space of dimension $(m,n)$ such that $V_\0=\mathbb R^m$ and
$V_\1=\mathbb R^n$.

\subsection{Near superalgebras} We say that a commutative superalgebra $\PC$ is
\emph{near} (they are the \emph{alg\`{e}bres proches} of Weil
\cite{Wei53} ) if it is finite dimensional, local, and  with $\mathbb R$ as
residual field. For $\alpha \in\PC$, we denote by $\bb(\alpha )$ the canonical
projection of $\alpha $ in $\mathbb R$ ($\bb(\alpha )$ is the \emph{body} of
$\alpha $, and $\alpha -\bb(\alpha )$ ---a nilpotent element of $\PC$--- the
\emph{soul} of $\alpha $, according to the terminology of \cite{DeW84}). Let
$\a\in\PC_\0$ be an even element. Let $\phi\in\CC^{\oo}(\mathbb{R},W)$ be a
smooth function defined in a neighborhood of $\bb(\alpha )$ in $\mathbb R$, with
values in some Fr\'echet supervector space $W$. We freely use the notation
(where $\phi^{(k)}$ is the $k$-th derivative of $\phi$):
\begin{equation}
\label{FonctLisse} \phi(\alpha )=\sum_{k=0}^\infty \frac{(\alpha -\bb(\alpha
))^k }{k!}\phi^{(k)}(\bb(\alpha ))\in W\otimes \PC.
\end{equation}
Since $(\alpha-\bb(\alpha ))$ is nilpotent, the sum is finite. In particular, if
$\alpha\in \PC_\0$ is invertible, its absolute value $|\alpha |\in \PC_\0$ is
defined by the formula:
\begin{equation}
|\alpha |={\frac{|\bb(\alpha )|}{ \bb(\alpha )}}\alpha ,
\end{equation}
and if $\bb(\a)>0$, its square root is defined by the finite sum 
\begin{displaymath}
\sqrt\a=\sqrt{\bb(\a)}\Big(1+\frac12(\frac{\a}{\bb(\a)}-1)-\frac
1{2^22!}\big(\frac{\a}{\bb(\a)}-1)^2+\frac 3 {2^3 3!}(\frac{\a}{\bb(\a)}-1)^3
-\frac {3.5} {2^4 4!}(\frac{\a}{\bb(\a)}-1)^4+\dots\Big),
\end{displaymath}
 where $\sqrt \l$ is the unique positive square root of $\l>0$. A notation like
$\sqrt{|\a|}$ (for $\a\in\PC_{\0}$ invertible) is not ambiguous, because if
$f\in\CC^\oo(\mathbb{R}^+)$ and $g\in\CC^\oo(\mathbb{R},\mathbb{R}^+)$ , for
$\a\in\PC_{\0}$ $f\circ g(\a)= f(g(\a))$.

\subsection{Modules over a superalgebra}

\subsubsection{Definition} We fix  a commutative superalgebra  $\AC$. We
define
an $\AC$-module $V$ as an  $\AC$-bimodule (in particular $V$ is a supervector
space) such that the left and right $\AC$-module structures satisfy  the  rule 
of signs: for $a\in \AC$ and $v\in V$ non zero and homogeneous:
\begin{displaymath} av=(-1)^{p(a)p(v)}va.
\end{displaymath}

If $V$ and $W$ are $\AC$-modules we denote  by $\Hom_{\AC}(V,W)$ the
superalgebra of morphisms of right $\AC$-modules. For $\phi\in \Hom_{\AC}(V,W)$
$v\in V$ and $a\in\AC$, $\phi(va)=\phi(v)a$. 

We put $\g{gl}(V)=\Hom_{\AC}(V,V)$.

We denote by $V^*$ its dual module $\Hom_{\AC}(V,\AC)$. It is naturally a left 
$\AC$-module. We give it a structure of $\AC$-(bi)module by the rule of signs.

Let $u\in\g{gl}(V)$ be non zero and homogeneous. We denote by $u^*$ the
endomorphism of $V^*$ such that for any non zero homogeneous $\phi\in V^*$ and
any $v\in V$, $u^*(\phi)(v)=(-1)^{p(u)p(\phi)}\phi(u(v))$. Let
$u,v\in\g{gl}(V)$ be non zero and homogeneous, we
have $(u\circ v)^*=(-1)^{p(u)p(v)}v^*\circ u^*$. We
extend
linearly the map $u\mapsto u^*$ to an even isomorphism supervector spaces of 
$\g{gl}(V)\to\g{gl}(V^*)$.
\subsubsection{Tensor product} For two $\AC$-modules $V$ and $W$ their tensor
product $V\tens_{\AC} W$ is defined by the relation:
\begin{displaymath}
\forall a\in\AC,\, \forall v\in V,\,\forall w\in W,\, va\tens w=v\tens aw.
\end{displaymath} We identify $V\tens_{\AC} W$ and $W\tens_{\AC} V$ by the usual
rule of signs (for non zero homogeneous $v\in V$ and $w\in W$ we identify $v\tens
w$ and 
$(-1)^{p(v)p(w)}w\tens v$).

When $\AC=\mathbb R$ and $V,W$ are Fr\'echet supervector spaces we denote by 
$V\what\tens W$ the completion of $V\tens W$.

 We denote by $S(V)$ the symmetric $\AC$-algebra of $V$. Let us recall its
definition. We put $V^{\tens_{\AC}0}=\AC$ and $V^{\tens_{\AC}n+1}
=V\tens V^{\tens_{\AC}n}$. We denote by $T(V)$ the tensor $\AC$-superalgebra of
$V$: 
\begin{displaymath} T(V)=\sdir_{n\in\mathbb N}V^{\tens_{\AC}n}.
\end{displaymath}

We denote by $\JC$ the ideal of $T(V)$ generated by 
\begin{displaymath}
\{v\tens w-(-1)^{p(v)p(w)}w\tens v\,/\, v,w\in V\text{ non zero and homogeneous}\}.
\end{displaymath} Then:
\begin{displaymath} S(V)=T(V)/\JC.
\end{displaymath} We denote by $S^n(V)$ the canonical projection of
$V^{\tens_{\AC}n}$ in $S(V)$.

When $\AC=\mathbb R$ we recall that $S(V)$ is equal to $S(V_\0)\otimes
\Lambda(V_\1)$, where $S(V_\0)$ and $\Lambda(V_\1)$ are the classical symmetric
and exterior algebras of the corresponding ungraded vector spaces. We  
use the notation $\Lambda(U)$ only for ungraded vector spaces $U$. So, if $V$ is
a supervector space, $\Lambda(V)$ is the exterior algebra of the underlying
vector space.

\subsubsection{Change of parity} For an $\AC$-module $V$, we denote by $\Pi V$
the $\AC$ supermodule
with reverse parity. It is defined by $(\Pi V)_{\0}=V_{\1}$ and $(\Pi
V)_{\1}=V_{\0}$. Let $\OddId:V\to \Pi V$ the ``odd identity''  from $V$ to $\Pi
V$. Let $a\in\AC,$ $v,w\in V$ be non zero homogeneous elements. The $\AC$-module
structure of $\Pi V$ is given by the relations:

\begin{equation}
\begin{split}
\OddId(va)&=\OddId(v)a;\\
\OddId(av)&=(-1)^{p(a)}a\OddId(v).
\end{split}
\end{equation}
\medskip

For $\phi\in V^*$ non zero and homogeneous and $v\in V$ we put:
\begin{equation}
(\OddId \phi)(\OddId v)=(-1)^{p(\phi)}\phi(v).
\end{equation}
This can be extended $\AC$-linearly to an odd isomorphism from $V^*$ to $(\Pi
V)^*$. This gives an identification between $(\Pi V)^*$ and $\Pi(V^*)$ (and
thus we denote it by $\Pi V^*$).

\medskip

 Similarly we put for any non zero homogeneous $\phi\in\g{gl}(V)$ and any $v\in
V$:
\begin{equation}
\label{eq:glPi}
\phi(\OddId v)=(-1)^{p(\phi)}\OddId\phi(v).
\end{equation}
It induces an even  isomorphism $\g{gl}(V)\to\g{gl}(\Pi V)$. We identify this
way this two algebras.

\subsubsection{Finite rank free modules} Let $V$ be a free $\AC$-module. We say
that $V$ is of rank $(k,l)\in \mathbb N^2$ if there exists a basis
$(e_{1},\dots,e_{k},f_{1},\dots,f_{l})$ of $V$ such that $p(e_{i})=0$ and
$p(f_{j})=1$. We also denote by $(e_{i},f_{j})$ such a standard basis. We write
$\rk(V)=(k,l)$. We put $V^{\mathbb R}=\mathbb R e_{1}\oplus
\dots\oplus\mathbb R e_{k}\oplus\mathbb Rf_{1}\oplus\dots\oplus\mathbb R f_{l}$.
Then, we have $V\simeq V^{\mathbb R}\tens_{\mathbb R}\AC$. We stress that the
definition of $V^{\mathbb R}$ depends on  the choice of the basis
$(e_{1},\dots,e_{k},f_{1},\dots,f_{l})$ of $V$. 

In this case $V^*\simeq \AC\tens_{\mathbb R}(V^{\mathbb R})^*$ and: $S(V)\simeq
S(V^{\mathbb R})\tens_{\mathbb R}\AC.$

If $V$ is free of rank $(k,l)$, $\Pi V$ is free of rank $(l,k)$.

\medskip

We recall that an $\AC$-module $V$ is said to be {\em locally free} is for any
maximal ideal $\g m$ of 
$\AC$, $V\tens_{\AC}\AC_{\g m}$ is free ($\AC_{\g m}$ is the localization of
$\AC$ at $\g m$).

\subsubsection{Bilinear forms}\label{sec:bilin} Let $V$ be an $\AC$-module. Let
$B$
be a {\em bilinear form} on $V$. It means that for $v,w\in V$ and $a,b\in \AC$:
\begin{equation}
B(va,wb)=(-1)^{p(a)p(w)}B(v,w)ab.
\end{equation}

We say that $B$ is symmetric (resp. antisymmetric) if for $v,w$ non zero
homogeneous:
\begin{equation}
B(v,w)=(-1)^{p(v)p(w)}B(w,v) \quad (\text{resp. }B(v,w)=-(-1)^{p(v)p(w)}B(w,v)).
\end{equation}
In this case, we  denote by $\g{osp}(V,B)$ (resp. $\g{spo}(V,B)$) the Lie
subsuperalgebra of $\g{gl}(V)$ such that for $\bk\in\mathbb Z/2\mathbb Z$:
\begin{multline}
\g{osp}(V,B)_\bk=\,(\text{resp. } \g{spo}(V,B)_{\bk}=)\\
\big\{X\in\g{gl}(V)_\bk\,/\,\forall  v,w\in V\setminus\{0\},\text{  homogeneous },
B(Xv,w)+(-1)^{\bk p(v)}B(v,Xw)\big\}.
\end{multline}
 When $B$ is clear from the context we will denote it by $\g{osp}(V)$ (resp.
$\g{spo}(V)$).

 \medskip

We say that $B$ is even if for any non zero homogeneous $v,w\in V$,
$p(B(v,w))=p(v)+p(w)$. If $B$ is even antisymmetric and non degenerate, we say
that $(V,B)$ (or $V$ if $B$ is clear from the context) is a {\em symplectic
supervector space} (cf. section
\ref{SympVect} for more details on symplectic supervector spaces).

\medskip

Let $V$ be an $\AC$-module and $B$ an even bilinear form on $V$. We define a
bilinear form $\OddId B$ on $\Pi V$ by 
\begin{equation}
\label{eq:PiB}
\OddId B(\OddId v,\OddId w)=(-1)^{p(v)}B(v,w).
\end{equation}
 If $B$ is symmetric (resp. antisymetric) $\OddId B$ is antisymmetric (resp.
symmetric).

\medskip

When $B$ is non degenerate it defines an isomorphism $B^*:V\to V^*$ defined for
$v,w\in V$ by:
\begin{equation}
B^*(v)(w)=B(v,w).
\end{equation}

\subsection{Symplectic oriented supervector spaces} \label{SympVect} Since this
particular case  is very important in this article we pay closer attention  to
it.  In this section $\AC=\mathbb R$ and $V$ is finite dimensional  but all the
definitions  makes sense if $\AC$ is a superalgebra and $V$ is a free $\AC$-module  of finite rank.

Let $V=V_\0\oplus V_\1$ be a symplectic supervector space. It means that $V_\0$
and $V_\1$
are orthogonal, that the restriction of $B$ to $V_\0$ is a non degenerate 
antisymmetric bilinear form, and that the restriction of $B$ to $V_\1$ is a non
degenerate     symmetric bilinear form. 

\medskip Such a space is a direct sum of $(2,0)$-dimensional symplectic
supervector spaces (i.e. $2$-dimensional symplectic  vector spaces), and of
$(0,1)$-dimensional   symplectic supervector spaces  (i.e. $1$-dimensional
quadratic  vector spaces). We first review these building blocks.

\subsubsection{Symplectic $2$-dimensional vector spaces}\label{2symp-a} Let
$V=V_\0$ a purely even $2$-dimensional   symplectic space. A
\emph{symplectic basis} of $V$ is a basis $(e_1,e_2)$ such that $B(e_1,e_2)=1$,
$B(e_1,e_1)=0$,  $B(e_2,e_2)=0$.
 It defines a dual \emph{symplectic coordinate system} $(x^1,x^2)$, and an
orientation of $V^*$.

\subsubsection{Symplectic $1$-dimensional odd vector spaces} Let $V=V_\1$ a
purely odd $1$-dimensio\-nal   symplectic supervector space (i.e. a
$1$-dimensional quadratic space). A symplectic basis of $V$ is a basis $(f  )$
such that $B(f ,f )= 1$. However, such a basis does not always exists, and we 
allow $f$ to be in $V_\1
\cup \ii V_\1 \subset V_\1\otimes \mathbb C$. Let $(\xi)\in V_\1^*
\cup \ii V_\1^* $ be the dual basis. 

We will call the choice of $(\xi)$ (the other possible choice is $(-\xi)$) an
\emph{orientation} of $V_\1$. If $B$ is positive definite, then $(\xi)$ is a
basis of $V_\1^*$, and so defines an orientation in the usual sense.  If $B$ is
negative definite, then $(-\ii\xi)$ is a basis of $V_\1^*$, and so defines an
orientation in the usual sense.

\subsubsection{General case} Let us go back to the general case.

The dimension $m$ of $V_\0$ is even. We choose a  symplectic basis
$(e_1,\dots,e_m)$ of $V_\0$, that is  $V_\0$ is the direct sum of $m/2$
symplectic vector spaces generated by the pairs $(e_1,e_2) ,
 (e_3,e_4) , \dots$, and  $B(e_1,e_2)=1, B(e_3,e_4)=1, \dots$. The dual basis
$(x^i)$ of $V_\0^*$ is called a symplectic coordinate system. 

\medskip On $V_\1$, we define a symplectic basis $(f_1,\dots f_n)$ as an
orthonormal basis of $V_\1\otimes \mathbb C$ such that $f_i\in V_\1$ or $f_i\in
\ii V_\1$ for all $i$. Let $(\xi^1,\dots\xi^n)$ be the dual basis. The pair of
functions $\pm \xi^1\dots\xi^n $ does not depend on the symplectic basis
$(f_1,\dots,f_n)$. A choice of one of the two elements of $\pm \xi^1\dots\xi^n $
is called an \emph{orientation} of $V_\1$. If $V_\1$ is oriented, an oriented
symplectic coordinate system on $V_\1$  is a basis for which the orientation is
$\xi^1\dots\xi^n $. If $B|_{V_{\1}}$ has signature $(p,q)$, then
$(-\ii)^{q}\xi^1\dots\xi^n\in\L^n(V_{\1}^*)$ and so defines an orientation of
$V_{\1}$ in the usual sense.

\medskip Let us   remark that in the specially interesting case where $B$ is
positive definite,   a symplectic basis is a basis of $V_\1$, and not only of
$V_\1\otimes \mathbb C$.

\medskip We define an \emph{oriented symplectic supervector space} as a
symplectic supervector space $(V,B)$ provided with an orientation of $V_\1$. 

A symplectic oriented basis $(e_{i},f_{j})=(e_1,\dots,e_m, f_1,\dots,$ $f_n)$ is
a basis of $V\otimes \mathbb C$ such that $(e_1,\dots,e_m)$ is a symplectic
basis of $V_\0$, and $(f_1,\dots,f_n)$ an oriented symplectic basis of
$V_\1\otimes
\mathbb C$. We use the corresponding dual system of coordinates
$(x,\xi)$. Then $V$ becomes a symbol   bearing a supermanifold 
structure, a symplectic structure, an orientation....

\subsubsection{Moment application}\label{sec:moment}

A particular  linear even isomorphism $\mu$ of $\g{spo}(V)$ to $ S^2(V^*)$,
called the moment mapping,  will play an important role. Thus, $\mu$ is an
element of $S(\g{spo}(V)^*\tens V^*)$. It is defined by the formula

\begin{equation}
\label{eq:moment}
\mu(X)(v)=-\frac{1}{2} B(v,Xv),
\end{equation}

where, for any commutative superalgebra $\PC$ (not necessarily near), 
$X\in\big(\g{spo}(V)\tens\PC\big)_{\0}$
and
$v\in\big(V\tens\PC)_{\0}$ are $\PC$-valued points of $\g{spo}(V)$  and $V$
(cf. next section for a precise definition of $\PC$-valued points),
and $B(v,Xv)\in
\PC$ is defined by the natural extension of scalars.
Considering a basis $G_k$ of $\g{spo}(V)$, the dual basis $Z^k$, the generic
point $X=G_k Z^k$, a basis $g_i$ of $ V $, the dual basis $z^i$, and the generic
point $v=g_i z^i$, we obtain: $$
\mu =-\frac{1}{2} B(g_i,G_k g_j)z^j  Z^k z^i. $$

\medskip

Let us explain the choice of the constant $-\frac12$ in definition of $\mu$ and
why we call $\mu$ the moment mapping. 

The symplectic form on $V$ gives to the associated supermanifold a structure of
symplectic supermanifold (cf. next section for a precise definition of a
supermanifold).
 We define a Poisson bracket on $S(V^*)$ by the following. Let $f\in V^*$ we
denote by $v_{f}$ the element of $V$ such that for any $w\in V$:
\begin{equation}
\label{eq:IsoSymp} 
B(v_{f},w)=f(w).
\end{equation}
In other words $B^*(v_{f})=f$. This gives an isomorphism from $V^*$ onto $V$.
For $f,g\in V^*$, we put:
\begin{equation}
\{f,g\}=B(v_{f},v_{g})\in\mathbb R\subset S(V^*)
\end{equation}
and we extend it to a Poisson bracket on $S(V^*)$.

With the above definitions we have:
\begin{equation}
\label{eq:Poisson}
\{\mu(X),\mu(Y)\}=\mu([X,Y]).
\end{equation}
Thus $\mu$ is a morphism of Poisson algebras.

\subsection{Supermanifolds}\label{superman} By a \emph{supermanifold} we mean a
smooth real supermanifold as in \cite{Kos77}, \cite{Ber87} (cf. below for a
precise definition).

\subsubsection{Affine supermanifolds} Let $V$ be a finite dimensional
supervector space. We denote also by $V$ the associated supermanifold. We say
that $V$ is an \emph{affine supermanifold}. We recall some relevant  definitions
in this particular and fundamental case.

Let $\UC \subset V_\0$ be an open set. We put $$
\CC_V^\infty(\UC)=\CC^\infty(\UC)\otimes \Lambda(V_\1^*),$$
 where $\CC^\infty(\UC)$ is the usual algebra of smooth real valued functions
defined in $\UC$, and $\Lambda(V_\1^*)$ is the exterior algebra of $V_\1^*$. We
say that $\CC_V^\infty(\UC)$ is the superalgebra of \emph{smooth functions on
$V$ defined in $\UC$}. The supermanifold $V$ is by definition the topological
space $V_\0$ equipped with the sheaf of superalgebras $\CC_V^\infty$. We denote
by $V(\UC)$  the supermanifold such that $V(\UC)_{\0}=\UC$ and sheaf of
functions $\CC_{V}^{\oo}$ restricted to the open subsets of $\UC$.

Notice that if $\UC$ is not empty, there is a canonical inclusion $S(V^*)
\subset\CC^\infty(\UC)$. The corresponding elements are called \emph{polynomial
functions}. One can also define rational functions. Similarly, if $W$ is a
Fr\'{e}chet supervector space (for instance $W=\mathbb C$), we denote by
 $\CC_V^\infty(\UC,W)=\CC^\oo(\UC,W)\tens\L(V_1)^*$ the space of $W$-valued
smooth functions.
 
\medskip Let $\PC$ be  a near superalgebra. We put $$
V_\PC=(V\otimes \PC)_\0. $$
It is called the \emph{set of points of $V$ with values in $\PC$}. Extending the
body $\bb :\PC \to \mathbb R$ to a map $V\otimes \PC
\to V$, and restricting to $V_\PC$, we obtain a map, also called the body and
denoted by $\bb$, $$
\bb :V_\PC \to V_\0. $$
Let  $\UC \subset V_\0$ be an open set. We denote by $V_\PC(\UC)
\subset V_\PC$ the inverse image of $\UC$ in $V_\PC$. We have
$V_{\PC}(\UC)=V(\UC)_{\PC}$. It is known that $V_\PC(\UC)$ is canonically
identified to the set of even algebra homomorphisms   $\CC_V^\infty(\UC) \to
\PC$. Let $v\in V_\PC(\UC)$. We will denote the corresponding character by $\phi
\mapsto \phi(v)$, and say that $\phi(v) \in \PC$ is the value of $\phi\in
\CC_V^\infty(\UC) $ at the point $v$.

\medskip 

For example, let $ v=v_i p^i \in V_\PC$ (we use Einstein's summation
rule, and considering tensorisation by $\PC$ as an extension of scalars, we
write $v_i p^i$ instead of $v_i \otimes p^i$) where the $(v_i,p^i)\in V\times
\PC$ are a finite number of pair of homogeneous elements with the same parity.
Then $$
\bb(v)=v_i \bb(p^i) , $$
 which is in $V_\0$ since $\bb(p^i)=0$ if $p^i$ is odd. Let $\phi\in V^*$. We
denote by the same letter the corresponding element in $\CC_V^\infty(V_{0})$.
Then $\phi(v)=\phi(v_i)p^i$, and this formula in fact completely determines the
bijection between $V_\PC(\UC)$ and the set of even homomorphisms of algebras 
$\CC_V^\infty(\UC) \to \PC$ (cf.
\cite{Wei53}).

\medskip 

For   $\phi \in \CC_V^\infty(\UC,W)$, we denote by $\phi_\PC$ the
corresponding function $v \mapsto \phi(v)$ defined in $V_\PC(\UC)$. Then
$\phi_\PC \in \CC^\infty(V_\PC(\UC),W\tens\PC)$. The importance of this
construction is that for $\PC$ large enough (for instance if $\PC$ is an
exterior algebra $\Lambda \mathbb R^N$ with $N\geq
\dim(V_\1)$), the map $\phi \mapsto \phi_\PC$ is injective, which allows more or
less to treat   $\phi$ as an ordinary function.

\medskip

 We emphasize  the special case $\PC=\mathbb R$. Then $V_\mathbb
R=V_\0$, $V_{\mathbb R}(\UC)=\UC$, and $\phi_\mathbb R$ is the projection of
$\phi\in \CC_V^\infty(\UC)$ to $\CC^\infty(\UC)$ which naturally extends the
projection of $\Lambda V_\1^*$ to $\mathbb R$.

\medskip To help the reader, we give two typical examples.

\medskip
 Let  $V=\mathbb R^{(1,0)}$. Then $V_\0=V=\mathbb R$, and $V_\PC= \PC_\0$. Let
$\UC \subset\mathbb R$ be an open subset,
$\phi\in\CC_V^\infty(\UC,W)=\CC^\infty(\UC,W)$, and $\alpha\in V_\PC=
\PC_\0$ such that $\bb(\alpha)\in \UC$. Then $\phi(\alpha)\in \PC$ is defined by
formula (\ref{FonctLisse}).

\medskip
 Let   $V=\mathbb R^{(0,1)}$.  Then $V_\0=\{0\}$, and $V_\PC=
\PC_\1$.  Any $\phi\in\CC_V^\infty( \{0\},W) $ can be written as $\phi=c+\xi d$,
where $c$ and $d$ are elements of $W$, and $\xi$ the standard coordinate (the
identity function) on $\mathbb R$. Then, for $\alpha \in  \PC_\1$,
$\phi(\alpha)\in \PC$ is defined by formula $$
\phi(\alpha)=c+\alpha d .$$

\medskip
We return now to the general case.

We define a Fr\'echet topology on $\CC_{V}^{\oo}(\UC,W)$ by saying that
$(\phi_{n})_{n\in\mathbb N}\in\CC_{V}^{\oo}(\UC,W)^{\mathbb N}$ converges to
$\phi\in\CC_{V}^{\oo}(\UC,W)$ if for any normed near superalgebra $\PC$,
$((\phi_{n})_{\PC})_{n\in\mathbb N}\in\CC^{\oo}(V_{\PC}(\UC),W\tens\PC)^{\mathbb
N}$ converges uniformly to $\phi_{\PC}\in \CC^{\oo}(V_{\PC}(\UC),W\tens\PC)$on
any compact subset $K\subset V_{\PC}(\UC)$.

\medskip Let $\AC$ be a commutative superalgebra. We still use the notation
$V_\AC=(V\otimes \AC)_\0$.  Polynomial functions $S(V^*)$ can be evaluated on
$V_\AC$, but,  in general, smooth functions can be evaluated  on $V_\AC$ only if
$\AC$ is a near algebra.

\medskip The particular case $\AC=S(V^*)$ is important, because $V_\AC$ contains
a particular point, the \emph{generic point}, corresponding to the identity in
the identification of $Hom(V,V)$ with $V\otimes V^*$. Let us call $v$ the
generic point. Then we have $f(v)=f$ for any polynomial function $f \in S(V^*)$.

\medskip

\medskip

\subsubsection{Coordinates}\label{coordinates}

 Let $V$ be a finite dimensional supervector space. By  a  basis $(g_i)_{i\in
I}$ of $V$, we mean an indexed basis consisting of homogeneous elements. The  
\emph{dual basis} $(z^i)_{i\in I}$ of $V^*$ is defined by the usual relation
$z^j(g_i)=\delta^j_i$ (the Dirac symbol). We will also say that the basis
$(g_i)_{i\in I}$ is  the
\emph{predual basis} of the basis $(z^i)_{i\in I}$ (the dual basis, in the
canonical identification of $V$ to the dual  of $V^*$ is 
$((-1)^{p(g_i)}g_i)_{i\in I}$). A basis $(z^i)_{i\in I}$ of $V^*$ will be also
called \emph{a system of coordinates on $V$}. The corresponding vector fields
(i.e. derivations of the algebra of smooth functions on $V$) are denoted by
$\frac{\de}{\de z^i}$. They are characterized by the rule $$
\frac{\de}{\de z^j}(z^i)=\delta^i_j. $$
Notice that the generic point $v$ of $V$ is then given by the formula
\begin{eqnarray}
\label{eq:generic} v=g_i z^i \in V_{S(V^*)}.
\end{eqnarray}

\medskip We will mainly use standard coordinates. Let $(m,n)=\dim(V)$. Then they
are basis of $V^*$ of the form $(x^1,\dots,x^m,\xi^1,\dots,\xi^n)$, where
$(x^1,\dots,x^m)$ is a basis of $V_\0^*$, and $(\xi^1,\dots,\xi^n)$ a basis of
$V_\1^*$. Such a basis will be sometimes denoted by the symbol $(x^{i},\xi^{j})$
or $(x,\xi)$. For the corresponding predual basis
$(e_1,\dots,e_m,f_1,\dots,f_n)$ of $V$, then $(e_1,\dots,e_m)$ is a basis of
$V_\0$ and $(f_1,\dots,f_n)$ is a basis of $V_\1$. These notations will be used
in particular for the canonical basis of $\mathbb R^{(m,n)}$. Let $I= (i_1,
\dots,i_n)\in\{0,1\}^n$.
Then we denote by $\xi^I$ the monomial $(\xi^1)^{i_1}\dots( \xi^n)^{i_n}$ of
$S(V^*)$. Let $\UC \subset V_\0$ be an open set.
Let $W$ be a Fr\'{e}chet supervector space. Then any $ \phi\in
\CC_V^\infty(\UC,W)$ is of the form
\begin{eqnarray}
\label{eq:funct}
\phi=\sum_I\xi^I \phi_I(x^1,\dots,x^m),
\end{eqnarray}
with $\phi_I$ is an ordinary $W$-valued smooth function defined in the
appropriate open subset of $\mathbb R^m$. Notice that $\phi_{\mathbb R}=
\phi_{(0,\dots,0)}(x^1,\dots,x^m)$   does not depend on the choice of the odd
coordinates
$\xi^i$. We emphasize the fact that we write $\phi_I$ to the right of  $\xi^I$
(recall that $\xi^I \phi_I= \pm \phi_I \xi^I $, according to the sign rule).

\subsubsection{General supermanifolds}

We define {\em supermanifold} $M$ of dimension $(m,n)$ as a pair 
$(M_{\0},\OC_{M})$ where $M_{\0}$ is a smooth manifold of dimension $m$ and
$\OC_{M}$ is a sheaf of commutative superalgebras on $M_{\0}$ such that locally
on a sufficiently small open subset $\UC$ of $M_{\0}$:
\begin{equation}
\OC_{M}(\UC)\simeq\CC^{\oo}(\UC)\tens\Lambda((\mathbb{R}^n)^*)
\end{equation}
where $\CC^{\oo}(\UC)$ denotes the algebra of smooth functions on $\UC$.
Sometimes, when $\UC=M_{\0}$ we will use the abuse of notations
$\OC(M)=\OC_{M}(M_{\0})$. If moreover $\UC$ is a coordinates set of $M_{\0}$, we
say that it is a coordinates set of $M$. In this case, let $(x^1,\dots,x^{m})$
be coordinates on $\UC$ and $(\xi^1,\dots, \xi^n)$ be the standard basis of
$(\mathbb{R}^n)^*$, we say that $(x^1,\dots,x^m,$ $\xi^1,\dots,\xi^n)$ is a
coordinates system on $\UC$. We will also use the notation $(x^{i},\xi^{j})$ or
$(x,\xi)$ to denote such a coordinates system. 

We stress that if $\UC$ is a coordinates set of $M$, the supermanifold
$M(\UC)=(\UC,\OC_{M}|_{\UC})$  isomorphic to the affine supermanifold
$\mathbb{R}^{(m,n)}$.

\medskip

Let $\PC$ be a near superalgebra. As for affine supermanifolds we can define the
set
 $M_{\PC}$ of points of $M$ with values in $\PC$. It is the set of even algebra
homomorphisms 
$\OC_{M}(M_{\0})\to \PC$.

 Let $v\in M_{\PC}$. We still denote by $\phi\mapsto\phi(v)$ the corresponding
character.

For $\PC=\mathbb{R}$, we have $M_{\mathbb{R}}=M_{\0}$. 

The canonical projection $\bb:\PC\to\mathbb{R}$ induces a projection
$\bb:M_{\PC}\to M_{\mathbb{R}}=M_{\0}$. Let $v\in M_{\PC}$ and
$\phi\in\OC_{M}(M_{\0})$, then by definition $\phi(\bb(v))=\bb(\phi(v))$. This
definition coincides in case  of an affine supermanifold $V$ with the body map
$\bb:V_{\PC}\to V_{\0}$. Let $\UC$ be an open subset of $M_{\0}$. As in the
affine case we denote by $M_{\PC}(\UC)$ the inverse image of $\UC$ in $M_{\PC}$.
It is canonically isomorphic to the set of even algebra homomorphisms
$\OC_{M}(\UC)\to \PC$.

 The set $M_{\PC}(\UC)$ is canonically provided with a structure of a smooth
manifold
  (cf. \cite{Wei53}). As in the affine case,  for $\phi\in\OC_{M}(\UC)$,  we
denote by 
  $\phi_{\PC}:v\in M_{\PC}(\UC)\mapsto\phi(v)$ the
   corresponding function. It is a smooth function on $M_{\PC}(\UC)$.

We stress that on a coordinates set $\UC$ of $M$ the map
$\phi\mapsto\phi_{\mathbb R}$ is the projection $\OC_{M}(\UC)\to\CC^{\oo}(\UC)$
which extends the canonical projection from $\L((\mathbb R^m)^*)$ to $\mathbb
R$.

\medskip

Let $\UC$ be an open subset of $M_{\0}$. As in the affine case, we provide
$\OC_{M}(\UC)$ with a Fr\'echet space topology. 

\subsubsection{Morphisms of supermanifolds}

Let $M,N$ be two supermanifolds. A morphism   $\pi:M\to N$ of supermanifolds is
a morphism $(\pi_{\0},\pi^*)$ topological spaces with sheafs of commutative
superalgebras. Then $\pi_{\0}:M_{\0}\to N_{\0}$ is a morphism of smooth
manifolds and for an open subset $\UC\in N_{\0}$,
$\pi^*:\OC_{N_{\0}}(\UC)\to\OC_{M_{\0}}(\pi_{\0}^{-1}(\UC))$ is an even morphism
of superalgebras.

\medskip

Let 
$\PC$ be a near superalgebra. Then $\pi$ induces a morphism 
$\pi_{\PC}:M_{\PC}\to N_{\PC}$ of smooth manifolds. In particular we have
 $\pi_{\mathbb R}=\pi_{\0}$.

Let $\phi\in\OC(N)$. We will also use the notation $v\mapsto\phi(\pi(v))$ to
denote the function $v\mapsto(\pi^*\phi)(v)$. Such a notation is not confusing
because for $v\in M_{\PC}$, $(\pi^*\phi)(v)=\phi(\pi_{\PC}(v))$.

\subsection{Supergroups}
\subsubsection{Defintion} We recall (cf. \cite{Kos77}) that a {\em supergroup} 
$G$ can be defined a a couple $(G_{\0},\g g)$ where $G_{\0}$ is a Lie group, $\g
g$ is a Lie superalgebra such that $\g g_{\0}$ is the Lie algebra of $G_{\0}$
and there is a representation $\Ad:G_{\0}\to GL(\g g_{\0})\times GL(\g g_{\1})$
such that its
differential is the restriction to $\g g_{\0}$ of the adjoint representation of 
$\g g$.

Such a supergroup defines a supermanifold $G=G_{\0}\times\g g_{\1}$ (where $\g
g_{1}$ is the purely odd affine supermanifold associated to the odd supervector
space $\g g_{\1}$).

Let $\PC$ be any near superalgbra then $G_{\PC}$ is canonically a Lie group with
Lie algebra $\g g_{\PC}$.

\medskip {\em Example 1: $GL(V)$.} Let $V$ be a supervector space. We denote by
$GL(V)$ the supergroup with underlying Lie group $GL(V)_{\0}=GL(V_{\0})\times
GL(V_{\1})$ (where $V_{\0}$ and  $V_{\1}$ are considered as ungraded vector
spaces) and  Lie superalgebra $\g{gl}(V)$. When $V=\mathbb R^{(m,n)}$ we denote
it by $GL(m,n)$.

In this case $GL(V)_{\PC}=GL(V_{\PC})$ is the group of invertibles elements of
$\g{gl}(V)_{\PC}=\g{gl}(V_{\PC})$.

\medskip

{\em Example 2: $SpO$ and $SpSO$.} Let $V$ a symplectic supervector space We
denote by $SpO(V)$ the subsupergroup of $GL(V)$ such that
$SpO(V)_{\0}=Sp(V_{\0})\mult O(V_{\1})$ and Lie superalgebra $\g{spo}(V)$. We
denote by $SpSO(V)$ its connected component: $SpSO(V)_{\0}=Sp(V_{\0})\mult
SO(V_{\1})$ and the same Lie superalgebra.

When $V$ is identified to $\mathbb R^{(m,n)}$ by the choice of a symplectic
basis   we denote them by $SpO(m,n)$ and $SpSO(m,n)$ respectively.

\medskip
\subsubsection{Representation} Let $G=(G_{\0},\g g)$ be a Lie supergroup. Let
$V$ be a Fr\'echet supervector space. A representation $\r$ of $G$ in $V$ is a
differentiable representation $\r_{\0}:G_{\0}\to GL(V_{\0})\times GL(V_{\1})$
of $G_{\0}$ in $V$ and a representation also denoted by $\r:\g g\to\g{gl}(V)$ of
$\g g$
in $V$ such that $\r|_{\g g_{\0}}$ is the differential $d\r_{\0}$ of $\r_{\0}$.

Let $\PC$ be any near superalgbra. The representation $\r$ induces a
representation $\r_{\PC}$ of $G_{\PC}$ in $V_{\PC}$. It is defined by a morphism
$\r_{\PC}:G_{\PC}\to GL(V_{\PC})=GL(V)_{\PC}$.

\medskip

\subsubsection{Action on a supermanifold} Let $G=(G_{\0},\g g)$ be a supergroup
and $M$ be a supermanifold. We say that $M$ is a {\em $G$-supermanifold} or that
$G$ {\em  acts (on the right) on} $M$ if $G_{\0}$ acts on the right  of
$(M_{\0},\OC_{M}(M_{\0}))$  and there is a morphism of Lie superalgbras
$X\mapsto X_{M}$ of $\g g$ into the derivations of $\OC_{M}(M_{\0})$. We denote
by $x\in M_{\0}\mapsto xg$  the action of $g\in G_{\0}$ on $M_{\0}$ and
$\phi\in\OC_{M}(M_{\0})\mapsto g\phi$ its action on $\OC_{M}(M_{\0})$. 
For $g\in G_{\0},\, x\in M_{\0}$ and $\phi\in\OC_{M}(M_{\0})$, we have
$(g\phi)(x)=\phi(xg)$. Moreover,  we assume the following compatibility between
the actions of $G_{\0}$ and $\g g$. For any $X\in\g g_{\0}$ and
$\phi\in\OC_{M}(M_{\0})$,
\begin{equation}\label{X_{M}}
X_{M}\phi=\frac d{dt}\exp(tX)\phi\big|_{t=0}
\end{equation}
(where $\exp$ is the expononential map from $\g g_{\0}$ into $G_{\0}$).

Let $\PC$ be a near superalgebra. Such an action induces a Lie group right
action of $G_{\PC}$ on $M_{\PC}$.

\medskip

 Let $\UC$ be a $G_{\0}$-invariant subset of $M_{\0}$. Let
$\phi\in\OC_{M}(\UC)$. We say that $\phi$ is $G$-invariant if $\phi$ is
$G_{\0}$-invariant and  for any $X\in\g g$, $X_{M}\phi=0$. We denote by
$\OC_{M}(\UC)^G$ the set of $G$-invariant functions on $M$ defined on $\UC$.

Equivalently this means that for any near superalgebra $\PC$, $\phi_{\PC}$ is a
$G_{\PC}$-invariant function.

\subsection{Rapidly decreasing functions}

 We say (cf. for example \cite[Chapter 7]{Hor83}) that
$\phi\in\CC^{\oo}(\mathbb{R}^m)$ is rapidly decreasing if for any
$(\a_1,\dots,\a_m)\in\mathbb{N}^m$ and any $(\be_1,\dots,\be_m)\in\mathbb{N}^m$,
\begin{equation}
\label{eq:bounded}
\sup\Big|(x^1)^{\be_1}\dots(x^m)^{\be_m}
\frac{\de^{\a_1}}{\de (x^1)^{\a_1}}\dots \frac{\de^{\a_m}}{\de
(x^m)^{\a_m}}\phi(x^1,\dots,x^m)\Big|<+\oo.
\end{equation}
where $(x^1,\dots,x^m)$ are the canonical coordinates on $\mathbb{R}^m$.

\medskip

Let $V$ be a supervector space. Let $\phi\in\CC_V^{\oo}(V_{\0})$ be a smooth
function on $V$. Let $(x^1,\dots,x^m,\xi^1,\dots,\xi^n)$ be a basis of $V^*$. We
put $\phi=\som_{I}\xi^I\phi_I(x^1,\dots,x^m)$ where
$\phi_I\in\CC^\oo(\mathbb{R}^m)$.

We say that $\phi$ is {\em rapidly decreasing} if for any $I$, $\phi_I$
 is a rapidly decreasing function on $\mathbb{R}^m$.

This definition does not depend on the choice of the basis $(x^{i},\xi^j)$ of
$V^*$. We denote by $\mathscr{S}_{V}(V_{\0})$ or $\mathscr{S}(V)$ the set of
rapidly decreasing functions on $V$.
\medskip

Equivalently, $\phi\in\CC^{\oo}_{V}(V_{\0})$ is rapidly decreasing if for any
near superalgebra $\PC$, $\phi_{\PC}$ is rapidly decreasing on $V_{\PC}$.

\medskip 

Let $W$ be a Fr\'echet supervector space. Similarly we define the set
$\mathscr{S}_{V}(V_{\0},W)$ of rapidly decreasing functions with values in $W$.
In this case, condition (\ref{eq:bounded}) must be satisfied when the absolute
value $|.|$ is replaced by any seminorm defining the topology of $W$.

\subsection{Supertrace and Berezinians}\label{sec:Ber} Let $V$ be a supervector
space.

Let $\AC$ be a commutative superalgebra. We write an element of
$\g{gl}(V)_{\AC}$ in the form
\begin{equation}
M=
\begin{pmatrix} A&B\\
C&D\\
\end{pmatrix}\in\g{gl}(V)_{\AC}.
\end{equation}
where $A\in \g{gl}(V_0)\otimes \AC_0$, $D\in \g{gl}(V_1)\otimes
\AC_0$, $B\in Hom(V_\1,V_\0)\otimes \AC_\1$, and  $C\in Hom(V_\0,V_\1)\otimes
\AC_\1$

We recall the definition of the supertrace:
\begin{defi} The supertrace of $M\in\g{gl}(V)_{\AC}$ is defined by
\begin{equation}
\str(M)=\tr(A)-\tr(D),
\end{equation}
where $\tr$ is the ordinary trace.
\end{defi}

\medskip

Berezin introduced the following generalization of the determinant (cf.
\cite{Ber87,BL75,Man88}), called the
\emph{Berezinian}.

\medskip

If $D$ is invertible, we define:
\begin{equation}
\label{eq:bera}
\Ber(M)=\det(A-BD^{-1}C)\det(D)^{-1}.
\end{equation}

\medskip

Assume moreover that  $\AC$ is a near superalgebra. If $D$ is invertible, we
define (cf. \cite{Vor91}):
\begin{align}\label{eq:berd}
\Ber_{(1,0)}(M) &=\Big|\det(A-BD^{-1}C)\Big|\det(D)^{-1}.
\end{align}

All these functions are multiplicative.

\medskip

Recall that $\g{gl}(V)_\0$ consists of the matrices $
\begin{pmatrix}A&0\\
0&D
\end{pmatrix}$ with $A\in \g{gl}(V_\0)$ and $D\in \g{gl}(V_\1)$. We consider the
 open set
 $\UC'= \g{gl}(V_\0)\times GL(V_\1)$. Formula (\ref{eq:bera}) defines a rational
function on the open set $\UC'$  of the supermanifold $\g{gl}(V)$.  We still
denote by $\Ber$ and $\Ber_{(1,0)}$   the elements of
$\CC^\infty_{\g{gl}(V)}(\UC')$ whose evaluation in $\g{gl}(V)_\AC$ is given as
above.

Since $GL(V)_{\0}\subset\UC'$ these functions are well defined on $GL(V)$.

\medskip

We keep the preceding notations. We assume moreover that $V$ is a symplectic
supervector space, then if $M\in SpO(V)_{\AC}$:

\begin{equation}
\Ber(M)=\Ber_{(1,0)}(M)=\det(D-CA^{-1}B)\in\{\pm 1\};
\end{equation}
and if if $M\in SpSO(V)_{\AC}$:
\begin{equation}
\Ber(M)=\Ber_{(1,0)}(M)=1.
\end{equation}

\subsection{Berezinian modules}

(cf. \cite{Man88}) Let $\AC$ be a commutative superalgebra. Let $V$ be an
locally free  $\AC$-module of finite type. We put $K(V)=S(\Pi V)\tens S(V^*)$.
Let $(e_{i})$ and $(x^{i})$ be families of vectors in $V$ and $V^*$ respectively
such that for any $v\in V$ we have $v=\som_{i}e_{i}x^{i}(v)$. We put:

\begin{displaymath} d_{V}=\som_{i}\OddId e_{i} x^{i}.
\end{displaymath}

Then, by definition:
 \begin{equation}
\ber(V)=H(K(V),d_{V}),
\end{equation}
  where $H$ denotes the homology of $(K(V),d_{V})$.

\medskip

Let $\w\in K(V)$, we denote by $\ber(\w)$ its canonical image in $\ber(V)$.

In particular, assume that $V$ is free. Let
$(e_{1},\dots,e_{m},f_{1},\dots,f_{n})$ be a standard basis of $V$, then
\begin{equation}
\ber(V)=\AC\, \ber(\OddId e_{1}\dots \OddId e_{m}\xi^{1}\dots\xi^{n}).
\end{equation}

\medskip
Let $V,V'$ be two free  $\AC$-modules of finite type. Let $\phi:V\to V'$ be an isomorphism.
Then we define:
\begin{equation}
\begin{split}
\ber(\phi):\ber(V)&\to\ber(V')\\
\ber(\OddId e_{1}\dots\OddId e_{n}\xi_{1}\dots\xi_{m})&\mapsto
\ber(\OddId \phi(e_{1})\dots\OddId \phi(e_{n})(\phi^{-1})^*(\xi_{1})\dots,(\phi^{-1})^*(\xi_{m})).
\end{split}
\end{equation}

\medskip

In the particular case where $V'=V$, the map $\ber(\phi)$ coincides  with the multiplication by $\ber(\phi)\in\AC$ defined in the preceding section. We put $V^{\mathbb R}=\mathbb R e_{1}\oplus\dots\oplus \mathbb R e_{m}
\oplus \mathbb R f_{1}\oplus \dots \oplus \mathbb R f_{n} $ thus $V=V^{\mathbb
R}\tens \AC$. Let $\phi\in GL(V)\simeq GL(V^{\mathbb R})_{\AC}$. We have:

\begin{equation}
 \ber(\OddId \phi(e_{1})\dots \OddId
\phi(e_{m})(\phi^{-1})^*(\xi^{1})\dots(\phi^{-1})^*(\xi^n))=\Ber(\phi)
  \ber(\OddId e_{1}\dots \OddId e_{m}\xi^{1}\dots\xi^{n}).
\end{equation}

\bigskip

Now, for later use, we stress: 
\begin{equation}
K(\Pi V^*)=S(\Pi (\Pi V^*))\tens S((\Pi V^*)^*) =S(V^*)\tens S(\Pi V)=K(V)
\end{equation}

Now, $(\OddId x^i)$ is a family of $\Pi V^*$ and $(\OddId e_{i})$ is a family of
$(\Pi V^*)^*=\Pi V$. On the other hand,   we have $\OddId (\OddId x^{i})=x^{i}$.
Thus, by rule of signs:
\begin{equation}
d_{\Pi V^*}=\som_{i}x^{i} \OddId e_{i}=d_{V}.
\end{equation}
It follows that $(K(V),d_{V})$ and $(K(\Pi V^*),d_{\Pi V^*})$ are canonically
isomorphic. Therefore:
\begin{equation}\label{eq:BerPi}
\ber(V)=H(K(V),d_{V})=H(K(\Pi V^*),d_{{\Pi V^*}})=\ber(\Pi V^*).
\end{equation}

\bigskip

\section{Some supergeometry}

\subsection{Supervector bundles}
\subsubsection{Definition} Let $M$ be a supermanifold. We recall that a {\em
supervector bundle} $\VC$  on $M$ is  defined by the locally free sheaf $\C_\VC$
of $\OC_M$-modules of its sections. More precisely, $\C_{\VC}$ is a sheaf of
$\OC_{M}$ modules such that locally, $\C_{\VC}(\UC)$ is a free $\OC_{M}(\UC)$
module.

In this article, unless otherwise specified, the supervector bundle will be of
finite rank denoted by $(k,l)$. This means that  locally, $\C_{\VC}(\UC)$ is a
free $\OC_{M}(\UC)$ module of rank $(k,l)$. We put $\rk(\VC)=(k,l)$.

\medskip

 In particular the {\em tangent bundle} $TM$ of $M$ is the supervector bundle
which sheaf of sections is $\C_{TM}=\der_M$, the sheaf of derivations of
$\OC_M$. We stress that since by derivation  we mean  left derivation $\der_{M}$
is naturally a sheaf of left $\OC_{M}$-modules. As usual we provide it with a
structure of right $\OC_{M}$-module by the rule of signs.

\medskip

Let $x\in M_\0$. We put $\g m_{x}=\{\phi\in\OC_{M}(M_{\0})\,/\,\phi(x)=0)\}$. We
recall that $\g m_{x}$ is a maximal ideal of $\OC_{M}(M_{\0})$ and that
$\OC_{M}(M_{\0})/\g m_{x}\simeq \mathbb R$. 

Then we put $\VC_x=\C_\VC(M_\0)/\g
m_{x}=\C_\VC(M_\0)\tens_{\OC_{M}(M_\0)}\OC_{M}(M_{\0})/\g m_{x}$. It is the
fibre of $\VC$ at $x$. If $\VC=TM$ we denote it  by $T_xM$. It is a supervector
space with $\dim(\VC_x)=\rk(\VC)$. For $v\in\C_\VC(M_\0)$ we denote by $v(x)$
its image by the canonical projection  on $\VC_{x}$.

\medskip

{\em Example:} assume that $\VC$ is trivial. In this case there is a supervector
space $V$ such that for any open subset $\UC\subset M_{\0}$, we have
$\C_{M}(\UC)=V\tens\OC_{M}(\UC)$. Then for any $x\in M_{\0}$, $\VC_{x}=V$.
Moreover, let $(e_{i})$ be a basis of $V$ and $v=\som_{i}
e_{i}\phi_{i}\in\C_{\VC}(M_{\0})$ ($\phi_{i}\in\OC_{M}(M_{\0})$), we have
$v(x)=\som_{i}e_{i}\phi_{i}(x)$.

In general $\VC$ is not trivial but locally trivial: there exists a
supervector space $V$ such that for
$\UC$ sufficiently small, $\C_{\VC}(\UC)\simeq V\tens \OC_{M}(\UC)$. We call
such an open subset
$\UC$ a {\em trivialization subset} of $\VC$ and such a supervector space $V$
the {\em generic fibre }of $M$.
\medskip 
If $\UC$ is any open subset of $M$ we denote by
$\VC(\UC)$ the restriction of $\VC$ to $M(\UC)$.

\bigskip 
Let $M$ be a supermanifold. Let $\VC$ be a supervector bundle on $M$.

We denote by $\VC^*$ the dual supervector bundle of $\VC$. It is the supervector
bundle whose sheaf of sections is the dual sheaf of $\C_{\VC}$:
$\C_{\VC^*}(\UC)=\big(\C_{\VC}(\UC)\big)^*=\Hom(\C_{\VC}(\UC),\OC_{M}(\UC))$.

\medskip

We denote by $\Pi \VC$ the supervector bundle whose sheaf of sections is
$\Pi\C_{\VC}$: $\C_{\Pi\VC}(\UC)=\Pi\C_{\VC}(\UC)$.

\medskip
We denote by $S(\VC)$ the supervector bundle whose sheaf of sections is
$S(\C_{\VC})$: $\C_{S(\VC)}(\UC)=S\big(\C_{\VC}(\UC)\big)$.

\medskip

Let $\VC\to M$ and $\WC\to M$ be two supervector bundles. We denote by 
$\VC\tens\WC$ the supervector bundle on $M$ whose sheaf of sections is
$\C_{\VC}\tens_{\OC_{M}}\C_{\WC}$.

\medskip

We denote by $\ber(\VC)$ the {\em Berezinian bundle} of $\VC$. It is the vector
bundle on $M$ whose sheaf of sections is 
\begin{equation}
\UC\mapsto \C_{\ber(\VC)}(\UC)=\ber(\C_{\VC}(\UC)).
\end{equation}
Since $\ber(\C_{\VC})=\ber(\Pi\C_{\VC}^*)$ (cf. (\ref{eq:BerPi})) we canonically
identify $\ber(\VC)$ and $\ber(\Pi\VC^*)$.

\subsubsection{Supermanifold structure} Let $M$ be a supermanifold. Let $\VC$ be
a supervector bundle on $M$. As in the purely even case $\C_{\VC}$ determines a
supermanifold also denoted by $\VC$ with a canonical projection $\pi:\VC\to M$. 

Let $\UC$ be an open subset of $M_{\0}$.  Thus $\pi_{\0}^{-1}(\UC)$ is an open
subset of $\VC_{\0}$. We have a canonical injection
$S(\C_{\VC^*}(\UC))\hookrightarrow \OC_{\VC}(\pi_{\0}^{-1}(\UC))$. We call the
elements of $S(\C_{\VC^*}(\UC))$ the functions on $\VC$ polynomial in the
fibres.

{\em Example:} If $\VC$ is trivial,  $\VC=M\times V$ and
\begin{displaymath} S(\C_{\VC^*}(\UC))=S(V^*)\tens\OC_{M}(\UC)\subset
\CC^{\oo}_{V}(V)\what\tens\OC_M(\UC)\simeq \OC_{\VC}(\pi_{\0}^{-1}(\UC)).
\end{displaymath}

\medskip

Let $\PC$ be a near superalgebra. Then $\pi_{\PC}:\VC_{\PC}\to M_{\PC}$ is a
vector bundle on $M_{\PC}$. Let  $x\in M_{\PC}$. We denote by $(\VC_{\PC})_{x}$
the fibre of $\VC_{\PC}$ at $x$. We have $\bb(x)\in M_{\0}$ and:
\begin{displaymath}
(\VC_{\bb(x)})_{\PC}=(\VC_{\PC})_{x}=\{v\in\VC_{\PC}\,/\,\pi_{\PC}(v)=x\}.
\end{displaymath}

\subsubsection{Bilinear forms} A {\em bilinear form} $B$ on $\VC$ is a bilinear
form on the sheaf $\OC_{M}$-modules $\C_{\VC}$. Let $x\in M_{\0}$ and
$v_{x},w_{x}\in \VC_{x}$. We choose $v,w\in\C_{\VC}(M_{\0})$ such that
$v(x)=v_{x}$ and $w(x)=w_{x}$. We put
\begin{equation}
B_{x}(v_{x},w_{x})=B(v,w)_{\mathbb R}(x).
\end{equation}
This definition does not depends on the choice of $v$ and $w$. This defines a
bilinear form $B_{x}$ on the supervector space $\VC_{x}$.

If $B$ is antisymmetric, we denote by $\g{spo}(\VC)$ the supervector bundle such
that for any open subset $\UC$ of $M_{\0}$:
$\C_{\g{spo}(\VC)}(\UC)=\g{spo}(\C_{\VC}(\UC),B)$. If $B$ is symmetric we define
$\g{osp}(\VC)$ in the same way.
\medskip

The form $B$ determines a form $\OddId B$ on $\Pi\VC$ in the usual way.

\subsubsection{Rapidly decreasing functions along the fibres}

Let $\pi:\VC\to M$ be a supervector bundle on $M$. Let $W$ be a Fr\'echet
supervector space. Let $\phi\in\OC_{\VC}(\VC_{\0},W)$. We say that $\phi$ is
{\em rapidly decreasing along the fibres} if for any trivialization subset
$\UC$, such that $\C_{\VC}(\UC)\simeq V\tens\OC_{M}(\UC)$ ($V$ is the generic
fibre of $\VC$),
$\phi\in\mathscr{S}_{V}(V_{\0})\mathop{\what\tens}\limits_{\OC_{M}(\UC)}\OC_{M%
}(\UC,W)$. We denote by $\SC_{\VC	}(\VC_{\0},W)$ the set of rapily decreasing functions 
along the fibres of $\VC$.

Equivalently this means that for any near superalgbra $\PC$ and any $x\in
M_{\PC}$, $\phi_{\PC}|_{(\VC_{\PC})_{x}}\in
\CC^{\oo}_{(\VC_{\PC})_{x}}((\VC_{\PC})_{x},W)$ is rapidly decreasing.

\subsubsection{$G$-equivariant vector bundle} Let $G=(G_{\0},\g g)$ be a
supergroup. Let $M$ be a $G$-supermanifold. Let $\VC\to M$ be a vector bundle.
We say that $\VC$ is an equivariant vector bundle if there is a representation
of $G$ in the Fr\'echet supervector space $\C_{\VC}(M_{\0})$ satisfying the
following conditions. Let $g\in G_{\0}$, we denote by
$v\in\C_{\VC}(M_{\0})\mapsto gv$ its action on $\C_{\VC}(M_{\0})$ and we denote
by $\LC^{\VC}:\g g\mapsto \g{gl}(\C_{\VC}(M_{\0}))$ the representation of $\g
g$. For any non zero and homogeneous $v,w\in\C_{\VC}(M_{\0})$, any
$\phi\in\OC_{M}(M_{\0})$ any $g\in G_{\0}$ and any non zero and homogeneous $X\in
\g g$ we require that:
\begin{equation}
\begin{split} g(v\phi)&=(gv)(g\phi)\\
\LC^{\VC}(X)(v\phi+w)&=(\LC^{\VC}(X)v)\phi+(-1)^{p(X)p(v)}v(X_{M}\phi).
\end{split}
\end{equation}

\subsection{Euclidean superstructure}\label{SEucl}

Let $\VC\rightarrow M$ be a supervector bundle. An {\em Euclidean  structure} 
on $\VC$ is an an even non degenerate symmetric bilinear form $Q$ on $\C_\VC$
such that for any $x\in M_\0$ the symmetric bilinear form $Q_{x}$ on $\VC_x$ is
positive definite on the even part $(\VC_{x})_{\0}$. We stress that since $Q$ is
non degenerate, $Q_{x}$ is non degenerate and thus if we forget the graduation,
$((\VC_{x})_{\1},Q_{x}|_{(\VC_{x})_{\1}})$ is a symplectic vector space.

\medskip

If we assume only  that for any $x$ in $M_0$ the restriction of the form $Q_{x}$
to $(\VC_x)_{\0}$ is positive definite (we no longer assume that it is non
degenerate, and thus $(\VC_{x})_{\1}$ is no longer symplectic), we say that $Q$
is a {\em weak Euclidean structure}.

\medskip

By definition a (weak) Euclidean structure on  a supermanifold $M$ is 
a (weak) Euclidean structure on $TM$.

We stress that if $M$ has a (weak) Euclidean structure, then $M_\0$ is a
Riemannian manifold.

\subsection{Pseudodifferential forms} \label{Pseudo}

\subsubsection{Definition}

Let $M$ be a supermanifold and $TM$ its tangent bundle. We put:
\begin{equation}
\what M=\Pi TM.
\end{equation} 
  If $dim(M)=(n,m)$, then $dim\big(\what{M}\big)=(n+m,n+m)$. The
pseudodifferential forms on $M$ are the functions on $\what M$ (cf.
\cite{BL77a,BL77b,Vor91}).

 Let $\what\pi:\what M\rightarrow M$ be the canonical projection of $\what M$
onto $M$. Let $\UC$ be an open subset of $M_{\0}$. We have $\what
{M}(\UC)_{\0}=\what\pi_{\0}^{-1}(\UC)$. We put

\begin{equation}
\what\W_M(\UC)=\OC_{\what M}(\what  {M}(\UC)_{\0}).
\end{equation} 
It is by definition the algebra of {\em pseudodifferential forms} on $\UC$. We
put $\what\W(M)=\what\W_{M}(M_{\0})$. Let $W$ be a Fr\'echet supervector space
we put similarly:
\begin{equation}
\what\W_M(\UC,W)=\OC_{\what M}(\what  {M}(\UC)_{\0},W);
\end{equation} 
and $\what\W(M,W)=\what\W_{M}(M_{\0},W)$.

\medskip

In particular, when $W=\C_{\VC}(M_{\0})$ is the $\OC_{M}(M_{\0})$-module of
sections of a supervector bundle (possibly not of finite rank) on $M$,  we put:
\begin{equation}
\what\W(M,\VC)=\what\W(M,\C_{\VC}(M_{\0}))
=\C_{\VC}(M_{\0})\mathop{\what\tens}\limits_{\OC_{M}(M_{\0})}\what\W(M).
\end{equation}

Let us look at this pseudodifferential forms in coordinates.
\medskip

Let $\UC$ be a coordinates set of $M$. Let $(x^1,\dots,x^m,\xi^1,\dots,\xi^n)$
be a local coordinates system on $\UC$ (with respective parities $p(x^i)=0$ and
$p(\xi^j)=1$). 

We denote by $(dx^1,\dots,$ $dx^m,d\xi^1,\dots,d\xi^n)$ the associated
coordinates on the fibres of $\what M$ ($p(dx^{i})=1$ and $p(d\xi^j)=0$). They
are defined as follows. Let $(\frac{\de}{\de x^1},\dots,\frac{\de}{\de
x^m},\frac{\de}{\de \xi^1},\dots,\frac{\de}{\de \xi^n})$ be the local basis of
sections of $TM$ such that $\frac{\de}{\de x^i}x^j=\d_{i}^j$ (the Dirac symbol)
and $\frac{\de}{\de \xi^i}\xi^j=\d_{i}^j$.  Then $(\OddId\big(\frac{\de}{\de
x^1},\dots,\OddId\frac{\de}{\de x^m},\OddId\frac{\de}{\de
\xi^1},\dots,\OddId\frac{\de}{\de \xi^n})$ is a local basis of sections of $\Pi
TM$.
Let $(\big(\OddId\frac{\de}{\de x^1}\big)^*,\dots,\big(\OddId\frac{\de}{\de
x^m}\big)^*,\big(\OddId\frac{\de}{\de
\xi^1}\big)^*,\dots,\big(\OddId\frac{\de}{\de \xi^n}\big)^*)$ be its dual basis.
We have 
\begin{displaymath} dx^{i}=\big(\OddId\frac{\de}{\de x^i}\big)^* \text{ and }
d\xi^j=\big(\OddId\frac{\de}{\de \xi^j}\big)^*;
\end{displaymath} these are elements of  $\C_{(\Pi TM)^*}(\UC)$ and thus
functions on $\what M=\Pi TM$ that are linear along the fibres.

{\em Case where $M$ is an affine manifold $V$.} (In fact the case of a
coordinates set is isomorphic to it) Let $(x^1,\dots,x^m,\xi^1,\dots,\xi^n)$ be
a basis of $V^*$ (coordinates on $V$). We have $TV\simeq V\oplus V$ and $\what
V\simeq V\oplus \Pi V$. Then $(\OddId x^1,\dots,\OddId
x^m,\OddId\xi^1,\dots,\OddId\xi^n)$  is a basis of $(\Pi V)^*$. We put
\begin{displaymath} dx^{i}=\OddId x^{i}\text{ and }d\xi^{j}=\OddId \xi^{j}.
\end{displaymath}

\medskip

A pseudodifferential form $\w$ can be written locally as:

\begin{equation}
\w=\som_{I,J}dx^I\xi^J \w_{I,J}(x,d\xi)
\end{equation}
 where the $\w_{I,J}$ are smooth functions of the variables $x^i$ and $d\xi^j$.

\bigskip

Let $M,N$ be two supermanifolds. Let $\pi:M\to N$ be a morphisms of
supermanifolds. Let $\UC$ be an open subset of $N_{\0}$. Then $\pi$ induces as
usual a morphism $\pi^*:\what\W_{N}(\UC)\to\what\W_{M}(\pi_{\0}^{-1}(\UC))$ (cf.
\cite {Man88}).

\subsubsection{Exterior differential...} Let $\UC$ be a coordinates set of $M$.
On $\what\W_M(\UC)$ the exterior differential is given by the odd vector field
on $\what M$:
\begin{equation} 
d=\som_idx^i{\frac{\de}{\de x^i}}+\som_jd\xi^j{\frac{\de}{\de \xi^j}}.
\end{equation}  
Let $\z=\som_if_i{\frac{\de}{\de x^i}}+\som_jg_j{\frac{\de}{\de
\xi^j}}\in\C_{TM}(\UC)$ ($f_{i},g_{j}\in\OC_{M}(\UC)$) be  an homogeneous vector
field on $M$ (that is a derivation of $\OC_M(\UC)$). The inner product
$\iota(\z)$ and the Lie derivative $\LC(\z)$ of $\z$ are defined by the
following vector fields on $\what M$:
\begin{align}\label{contr}
\iota(\z)&=(-1)^{p(\z)}\Big(\som_if_i{\frac{\de}{\de
dx^i}}+\som_jg_j{\frac{\de}{\de d\xi^j}}\Big),\\
\LC(\z)&=[d,\iota(\z)]=d\iota(\z)+(-1)^{p(\z)}\iota(\z)d.
\end{align} 
These definitions does not depend on the choice of coordinates and thus can be
extended to any open subset $\UC$ of $M_{\0}$.

Thus $\iota$ (resp. $\LC$) is an odd (resp. even) morphism  of
$\OC_{M}(\UC)$-modules from $\C_{TM}(\UC)$ to $\C_{T\what M}(\what M(\UC))$. We
recall the Cartan relations:
\begin{align}
\LC(\z)f=&\z f.\\
\iota(\z)(df)&=(-1)^{p(\z)}\z f,\\
[\LC(\z),d]&=0,\\
[\iota(\gamma),\LC(\z)]&=\iota([\gamma,\z]),\\
[\iota(\gamma),\iota(\z)]&=0,\\
[\LC(\gamma),\LC(\z)]&=\LC([\gamma,\z]).
\end{align} 
In particular $\LC$ is a morphism of sheaf of Lie superalgebras.

\subsection{Orientation}\label{Ori}

\label{Ori.def} (cf. \cite{Vor91}) We say that a supermanifold $M$ is {\em
oriented} if $M_\0$ is oriented. 

We say that $M$ is {\em globally oriented} if $\what M$ is oriented.

\medskip

If $V$ is a supervector space, we say that it is oriented (resp. globally
oriented)  if $V_\0$ (resp. $V$ as a non graded vector space) is oriented.

We stress that in this case the definitions orientation and global orientation
for the supervector space $V$ or for its associated affine supermanifold
coincides.

{\em Example:} Let $V$ be an oriented symplectic supervector space. Then $V$ is
oriented in the above sense and also globally oriented. Thus the terminology is
not too confusing.

\medskip

For a supervector bundle $\VC\rightarrow M$ we say that $\VC$ is an oriented
(resp. globally oriented) supervector bundle if $\VC_\0\rightarrow M_\0$ (resp.
$(\what\VC)_\0\rightarrow(\what M)_\0$) is an oriented vector bundle.

{\em Example:} If $\VC=M\times V$ is a trivial supervector bundle. it is
oriented (resp. globally oriented) if $V$ is oriented (resp. globally oriented).

\medskip

\subsection{Oriented symplectic supervector bundle} \label{sec:SympOrBun} Let
$\VC\to M$ be a supervector bundle. We say that $\VC$ is symplectic if there is
an even non degenerate antisymmetric bilinear form $B$ on $\VC$. If moreover
$\VC$ is a globally oriented supervector bundle, we say that $\VC$ is an
oriented symplectic supervector bundle. This is equivalent to say that its
fibres are oriented symplectic supervector spaces or that $\C_{\VC}$ is a sheaf
of locally free oriented symplectic  $\OC_{M}$-modules.

\medskip

We stress that if $\VC$ as an Euclidean structure $Q$, then $(\Pi\VC,\OddId Q)$
is a symplectic supervector bundle. Moreover, if $\VC$ is an oriented
supervector bundle, $\Pi\VC$ is an oriented symplectic supervector bundle.

\bigskip

We stress that $M$ is globally oriented if $\ber(TM)$ is a trivial vector bundle
on $M$. In this case, up to a multiplication by a function
$\phi\in\OC_{M}(M_{\0})$ such that $\phi_{\mathbb R}>0$ there are exactly two
different basis of $\C_{\ber(TM)}(M_{\0})$. The choice of a global orientation
of $M$ correponds to the choice of such a basis.

 Let $(x^{i},\xi^{j})$ be a local coordinates system of $M$. Let
$(\frac{\de}{\de x^{1}},\dots,\frac{\de}{\de x^{m}},\frac{\de}{\de
\xi^{1}},\dots,\frac{\de} {\de\xi^{n}})$ be the corresponding local basis of
sections of $TM$. We say that $(x^{ i},\xi^{j})$ is globally oriented local
system of coordinates if locally the orientation of $M$ is defined by the basis:
 \begin{displaymath}
\ber\Big(\OddId\frac{\de}{\de x^{1}}\dots\OddId\frac{\de}{\de
x^{m}}\xi^{1}\dots\xi^{n}\Big)
\end{displaymath} of $\C_{\ber(TM)}(M_{\0})$.

\medskip

For a supervector bundle it is globally oriented (resp. oriented) if
$\ber(\VC)\to M$ (resp. $\ber(\VC_{\0})\to M_{\0}$) is trivial. In this case, the
choice
of an orientation corresponds to a choice of a basis of $\ber(\VC)$ (up to
multiplication  by $\phi\in\OC_{M}(M_{\0})$ with $\phi_{\mathbb R}>0$). We say
that a local basis of section $(e_{i},f_{j})$ with dual basis $(x^{i},\xi^{j})$
is globally oriented if the orientation is defined by 
 the basis:
 \begin{displaymath}
\ber\Big(\OddId e_{1}\dots\OddId e_{m}\xi^{1}\dots\xi^{n}\Big)
\end{displaymath} of $\C_{\ber(\VC)}(M_{\0})$.

\medskip
\subsection{Superconnections}

\subsubsection{Definitions}\label{sec:Superconn-Def}
 Let $\VC\rightarrow M$ be a supervector bundle (possibly infinite dimensional).
Let $\C_\VC$ be its sheaf of sections. Let
$\Gamma_\VC\tens\limits_{\OC_M}\what\W_M$ be the sheaf of pseudodifferential
forms on $M$ with values in $\VC$. For an open subset $\UC$ of $M_{\0}$ we put
\begin{displaymath}
\what\W_M(\UC,\VC)=\Gamma_\VC(\UC)\mathop{\what\tens}\limits_{\OC_M(\UC)}\what\W_M(\UC)
\end{displaymath}

According to Quillen (cf.
\cite{MQ86,BGV92}), a {\em superconnection} on $\VC$ is an odd endomorphism of
the sheaf of supervector spaces $\Gamma_\VC\mathop{\what\tens}\limits_{\OC_M}\what\W_M$ such
that, for $\a\in\what\W_M(\UC)$ and $\w\in\what\W_M(\UC,\VC)$ non zero and
homogeneous:
\begin{equation}
\mathbb A(\w\a)=(\mathbb A\w)\a+(-1)^{p(\w)}\w(d\a).
\end{equation}

\medskip

Locally, on a trivialization subset $\UC$ of $\VC$ there is
$\w\in\what\W_{M}(\UC,\g{gl}(\VC))_{\1}$ such that:
\begin{equation}
\mathbb A=d+\w.
\end{equation}

\medskip

The superconnection $\mathbb A$ acts on
$\g{gl}(\Gamma_\VC)\tens\limits_{\OC_M}\what\W_M$ by means of the
supercommutation bracket of endomorphisms of $\C_\VC\tens_{\OC_M}\what\W_M$.

\medskip

The curvature of a superconnection is the operator $F=\mathbb
A^2\in\g{gl}(\VC)\tens\what\W_{M}(\UC)$. Locally, $F=d\w+\w^2$. It  satisfies the
Bianchi identity:
\begin{equation}
\mathbb A (F)=0.
\end{equation}

\medskip Let $B$ be an even bilinear form on $\VC$. We extend it 
$\what\W_{M}$-linearly to a bilinear form on $\C_\VC\tens_{\OC_M}\what\W_M$.

We say that $\mathbb A$ preserves $B$ (or leave $B$ invariant) if for any non
zero
 homogeneous $\w,\w'\in\what\W_{M}(\UC,\VC)$ we have:
\begin{equation}
B(\mathbb A\w,\w')+(-1)^{p(\w)}B(\w,\mathbb A\w')=dB(\w,\w).
\end{equation}

In this case when $B$ is symmetric (resp. antisymmetric) we have: 
\begin{equation}
F=\mathbb A^2\in \g{osp}(\VC)\tens\what\W_{M}(\UC) \text{(resp. }\in
\g{spo}(\VC)\tens\what\W_{M}(\UC) \text{).}
\end{equation}

\subsubsection{``Induced'' superconnections}\label{IndConn}

Let $\VC\to M$ be a supervector bundle and $\mathbb A$ be a superconnection on
$\VC$.
We denote by $\mathbb A^{\OddId}$ the  superconnection on $\Pi\VC$ defined for
$\w\in\what\W(\UC,\VC)$by:
\begin{equation}
\mathbb A^{\OddId}(\OddId \w)=-\OddId (\mathbb A\w).
\end{equation}

\medskip

Recall that
 $\what\W(\UC,\Pi\VC)$ is the $\what\W_{M}(\UC)$-module 
$\C_{\Pi\VC}(\UC)\tens_{\OC_M(\UC)}\what\W_M(\UC)$. Now, $\mathbb A$  determines
also   a  superconnection $\mathbb A^{\OddId*}$ (resp. $\mathbb A^*$) on
$\Pi\VC^*$ (resp. $\VC^*$) by the following
formula. Let $\a\in\what\W(\UC,\Pi\VC^*)$ (resp. $\a\in\what\W(\UC,\VC^*)$) non
zero and homogeneous and
  $\be\in\what\W(\UC,\Pi\VC)$ (res.   $\be\in\what\W(\UC,\VC)$):
\begin{equation}
\label{A*}
\mathbb A^{\OddId*}(\a)(\be)=(-1)^{p(\a)}\a(\mathbb A^{\OddId}\be)\qquad
(\text{resp. } \mathbb A^{*}(\a)(\be)=(-1)^{p(\a)}\a(\mathbb A\be)).
\end{equation}

\medskip

We extend $\mathbb A^{\OddId*}$ (resp. $\mathbb A^*$)  to a derivation of the
algebra
$\what\W(\UC,S(\Pi\VC^*))=S\big(\what\W(\UC,\Pi\VC^*)\big)$ (resp. 
$\what\W(\UC,S(\VC^*))=S\big(\what\W(\UC,\VC^*)\big)$). Thus
$\mathbb A^{\OddId*}$ (resp. $\mathbb A^{*}$) is a superconnection on
$S(\Pi\VC^*)$ (resp. $S(\VC^*)$). It is
defined for $\a,\be\in\what\W(\UC,S(\Pi\VC^*))$ 
(resp. $\a,\be\in\what\W(\UC,S(\VC^*))$) non zero and homogeneous,
\begin{equation}
\label{A*S}\begin{split}
\mathbb A^{\OddId*}(\a\be)&=(\mathbb A^{\OddId *}\a)\be+(-1)^{p(\a)}\a(\mathbb
A^{\OddId*}\be)\\
(\text{resp. } \mathbb A^{\OddId*}(\a\be)&=(\mathbb A^{\OddId *}\a)\be+(-1)^{p(\a)}\a(\mathbb
A^{\OddId*}\be).
).
\end{split}
\end{equation}

\section{Integration}

\subsection{Integration on a supervector space}
 
\subsubsection{Definition} Let $V$ be a supervector space, $W$ be Fr\'echet
supervector space and $\UC$ be an open subset of $V_{\0}$.

We will denote by $\CC_{V,c}^\infty(\UC,W)$
  the subspace of $\CC_V^\infty(\UC,W)$ of function with compact support. 

\medskip

The  \emph{distributions on $V$ defined in $\UC$} are the elements of the
(Schwartz's) dual of $\CC_{V,c}^\infty(\UC )$. If $t$ is   a  distribution, we
will use the notation $$
t(\phi)=\int_{V } t(v) \phi(v) $$
for $\phi\in \CC_{V,c}^\infty(\UC)$. We will also use complex valued
distributions, defined in an obvious way.

A \emph{Berezin integral} (or  \emph{Haar}, or  \emph{Lebesgue}) is by
definition a distribution on $V$ which is   invariant by translations (i.e. 
which vanishes on functions of the form $\de_X\phi$ where $\phi\in
\CC^{\oo}_{V,c}(V)$ and $\de_X$ is the vector field on $V$ with constant
coefficients corresponding to $X\in V$: for $f\in V^*$,
$\de_{X}f=(-1)^{p(X)p(f)}f(X)$). Up to a multiplicative constant, there is
exactly one    Berezin integral, and it is an important matter in this article
to choose a particular one for the symplectic oriented supervector spaces (see
below and \cite{Lav03}).

\medskip A choice of a standard system of coordinates determines a specific
choice $d_{(x,\xi)}$ of a Berezin integral by the formula
\begin{equation}
\label{eq:haar}
\begin{split}
 \int_V d_{(x,\xi)}(v)\phi(v)&=(-1)^{\frac{n(n-1)}2}\int_{\mathbb R^m}
   \big|\,dx^1\ldots dx^m\big| \phi_{(1,\dots,1)}(x^1,\dots,x^n)\\
&=\int_{\mathbb{R}^m}  \big|\,dx^1\ldots dx^m\big|
\Big(\frac{\de}{\de \xi^1}\dots\frac{\de}{\de \xi^n}
\phi\Big)_\mathbb R (x^1,\dots,x^m),
\end{split}
\end{equation}
for $\phi\in \CC_{V,c}^\infty(V_\0)$, where $\big|\,dx^1\ldots dx^m\big|$ is the
Lebesgue measure on $\mathbb{R}^m$. 

Note that this formula can also be applied to any $\phi\in
\CC_{V,c}^\infty(\UC,W)$ with a result in $W$ and also to any rapidly decreasing
 
$\phi\in\CC_{V}^\infty(\UC,W)$.

The choice of sign is such that Fubini's formula holds. More precisely, let
$V,W$ be two supervector spaces of dimensions $(m,n)$ and $(p,q)$. Let
$(x^1,\dots,x^m,\xi^1,\dots,\xi^n)$ be standard coordinates on $V$ and
$(y^1,\dots,y^p,\eta^1,\dots,\eta^q)$ be standard coordinates on $W$. Then
$(x,y,\xi,\eta)=(x^1,\dots,x^m,y^1,\dots,y^p,\xi^1,\dots,\xi^n,\eta^1,\dots%
,\eta^q)$ defines standard coordinates on $V\times W$. Let $\phi(v,w)$ is a
smooth compactly supported function on $V\times W$. Then:
\begin{equation}
\label{eq:Fubini}
\int_{V\times W}d_{(x,y,\xi,\eta)}(v,w)\phi(v,w)=\int_Vd_{(x,\xi)}(v)\Big(\int_W
d_{(y,\eta)}(w)\phi(v,w)\Big).
\end{equation}
We write:
\begin{equation}
\label{eq:Fubini-int} d_{(x,y,\xi,\eta)}(v,w)=d_{(x,\xi)}(v)d_{(y,\eta)}(w).
\end{equation}

In particular, since $V=V_\0\oplus V_\1$, $d_{x,\xi}=d_x d_\xi$ and formula
(\ref{eq:haar}) is a particular case of formula (\ref{eq:Fubini}).

\medskip Let us stress that if $V_\1$ is not $\{0\}$, in the setting of
supermanifolds  there is no natural notion of measure  on $V$ and no natural
notion of positive distribution on $V$. Thus we use these notions  only for
(ungraded) vector spaces, or for the even part $V_\0$ of a supervector space,
which is then regarded as an  ungraded  vector space. Otherwise, we use the
terms \emph{distribution} or \emph{integral}.

\medskip In this article, we will be in fact interested by complex valued
distributions. Then we allow   standard basis $(e_1,\dots,e_m,f_1,\dots,f_n)$ of
$V\otimes \mathbb C$, where $(e_1,\dots,e_m)$ is a basis of $V_\0$ and
$(f_1,\dots,f_n)$ is a basis of $V_\1\otimes \mathbb C$. Then the  dual basis
$(x,\xi)$ provides a coordinate system $(x)$ on $V_\0$ and a dual basis $(\xi)$
of $V_\1^*\otimes \mathbb C$. Any $f \in
\CC_{V,c}^\infty(V_\0,\mathbb C)$ can be written in the form (\ref{eq:funct}),
and the (complex) Berezin integral $d_{(x,\xi)}$ is again well defined by  
formula (\ref{eq:haar}).

\medskip

\subsubsection{Change of variables} Now, we recall the formula for change of
variable in integration (cf. \cite{Ber87} for example).

Let $V$ be a supervector space. Let $(x^{i},\xi^{j})$ be a standard basis of
$V^*$. Let $h:V\to V$ be an isomorphism of supermanifolds. We put
$y^{i}=\phi^*x^{i}$ and $\eta^{j}=\phi^{*}\xi^{j}$.

Let $J(h)$ be the jacobian matrix of $h$ defined by:
\begin{equation}
J(h)=
\begin{pmatrix}
\frac{\de y_{j}}{\de x_{i}}&\frac{\de y_{j}}{\de \xi_{i}}\\
\frac{\de \eta_{j}}{\de x_{i}}&\frac{\de \eta_{j}}{\de \xi_{i}}
\end{pmatrix}\in
\CC^{\oo}_{V}(V_{\0},\g{gl}(V))_{\0}=\g{gl}(V)_{\CC_{V}^{\oo}(V_{\0})}.
\end{equation}
Moreover since $h$ is an isomorphism, $J(h)$ is invertible, and
$\Ber_{(1,0)}(J(h))$ is the function on $V$ such that for any near superalgebra
$\PC$ and any $v\in V_{\PC}$,
$Ber_{(1,0)}(J(h))(v)=Ber_{(1,0)}\big(J(h)(v)\big).$

Now the formula for change of variables is:

\begin{equation}
\label{ChangVar}
\int_{V}d_{(x,\xi)}(v)Ber_{(1,0)}(J(h))(v)\phi
(h(v))=\int_{V}d_{(x,\xi)}(v)\phi(v).
\end{equation}

\subsubsection{Link with $\ber(V)$}

Assume that $V$ is oriented. Let $h:V\to V$ be an isomorphism of supermanifolds.
Then if $h$ preserves orientation, that is if $(x^1,\dots,x^{m})$ and
$(h^{*}x^{1},\dots,h^{*}x^{m})$ define the same orientation, we have:
\begin{equation}
\Ber(J(h))=\Ber_{(1,0)}(J(h)).
\end{equation}
Thus map:
\begin{displaymath} d_{(x,\xi)}\mapsto \ber(\pi e_{1}\dots\pi
e_{m}\xi^{1}\dots\xi^{n})
\end{displaymath} induces an isomorphism from the $(1,0)$ or $(0,1)$-dimensional
(depending on the parity of $n$) vector space of Berezin integrals on $V$ and
$\ber(V)$. This isomorphism is even or odd depending on the parity of $m$. When
$V$ is oriented we will identify $d_{(x,\xi)} $ and $\ber(\pi e_{1}\dots\pi
e_{m}\xi^{1}\dots\xi^{n})$.

\subsection{Integration in symplectic oriented supervector spaces}
\label{SympInt} Let $V=V_\0\oplus V_\1$ be an oriented symplectic supervector
space (cf. subsection \ref{SympVect}).

\medskip Since such a space is a direct sum of $(2,0)$-dimensional symplectic
supervector spaces, and of $(0,1)$-dimensional   symplectic supervector spaces 
(i.e. $1$-dimensional quadratic  vector spaces). We first review these building
blocks.

\subsubsection{Symplectic $2$-dimensional vector spaces}\label{2symp} Let
$V=V_\0$ a purely even $2$-dimensional   symplectic space. Let  $(x^1,x^2)$ be a
symplectic coordinate system.   It defines a a
\emph{Liouville integral} (a particular normalization of the Berezin integral):
$$
\phi\in \CC^\infty_c(V)=\CC^\infty_c(\mathbb R^2)\mapsto
 \int_V d_V(v) \phi(v)
 =\frac{1}{2\pi }\int |dx^1  dx^2|\phi(x^1,x^2).
 $$

\subsubsection{Symplectic $1$-dimensional odd vector spaces} Let $V=V_\1$ a
purely odd $1$-dimensio\-nal  oriented symplectic supervector space.  Let
$(\xi)\in V_\1^*
\cup \ii V_\1^* $ be a symplectic oriented coordinates system. It defines a
\emph{Liouville integral} (which is complex valued if $B$ is negative definite)
$d_{V }$:
 $$
\phi=a+\xi b\in \Lambda(V^*\otimes \mathbb C)  \mapsto \int_V d_{V}(v)
\phi(v)=b.$$

\subsubsection{General case} Let us go back to the general case. Since $V_\0$ is
a classical symplectic space, there is a canonical normalization of Lebesgue
integral on $V_\0$, the Liouville integral,   which we recall now.

The dimension $m$ of $V_\0$ is even. Let $(x^i)$ be  a symplectic coordinate
system. The Liouville integral on $V_\0$ is $$
\frac{1}{(2
\pi)^{m/2}}|dx^1\dots dx^m|.$$
  The Liouville integral does not depend on the choice of the symplectic basis
of $V_0$.

\medskip Let $(\xi^1,\dots\xi^n)$ be  an oriented symplectic coordinate system
on $V_\1$.

We call the corresponding  Berezin integral $d_{\xi}$ the Liouville integral of
the oriented symplectic space $V_\1$.

\medskip Let $(x,\xi)$ be an oriented symplectic coordinates system on $V$. The
associated Berezin integral $\frac{1}{(2\pi)^{\frac n2}}d_{(x,\xi)}$ will be denoted by $d_V$.

\medskip

\subsection{Generalized functions} Let $V$ be a finite dimensional supervector
space and $\UC\subset V_0$ be an open set. Let $(x,\xi)$ be a standard
coordinates system on $V$.  We will say that a distribution $t$ on $V$ defined
in $\UC$ is \emph{smooth} (resp. \emph{smooth compactly supported}) if there is
a function $\psi\in\CC^{\oo}_V(\UC)$ (resp. $\psi\in\CC^{\oo}_{V,c}(\UC)$) such
that $t(v)=d_{x,\xi}(v)\psi(v)$. In particular, for any
$\phi\in\CC^{\oo}_{V,c}(\UC)$ (resp. $\phi\in\CC^{\oo}_V(\UC)$):
\begin{equation}
t(\phi)=\int_Vd_{x,\xi}(v)\psi(v)\phi(v).
\end{equation}
This definition does not depend on the standard coordinates system $(x,\xi)$.

\medskip

By definition, the \emph{generalized functions on $V$ defined on $\UC$} are the
elements of the (Schwartz's) dual of the space of smooth compactly supported
distributions. For a generalized function $\phi$ and a smooth compactly
supported distribution $t$, we write:
\begin{equation}
\phi(t)=(-1)^{p(t)p(\phi)}\int_V t(v)\phi(v).
\end{equation}
(The spaces of distributions and thus of generalized functions are naturally
$\mathbb{Z}/2\mathbb{Z}$-graded.)

We denote by $\CC^{-\oo}_V(\UC)$ the set of generalized functions on $V$ defined
on $\UC$.

Let us remark that, as $\CC^\oo_V(\UC)=\CC^{\oo}(\UC)\tens\Lambda( V_\1^*)$, we
have:
\begin{equation}
\CC^{-\oo}_V(\UC)=\CC^{-\oo}(\UC)\tens\Lambda (V_\1^*).
\end{equation}

\medskip

Let $W$ be a Fr\'echet supervector space. A $W$-valued generalized function  is
a continuous homomorphism (in sense of Schwartz) from the space of smooth
compactly supported distributions on $V$ to $W.$ We denote by
$\CC^{-\oo}_V(\UC,W)$ the set of $W$-valued generalized functions. If $W$ is
finite dimensional, we have $\CC^{-\oo}_V(\UC,W)=\CC^{-\oo}_V(\UC)\tens W$. We
will be  in particular concerned with the cases $W=\mathbb{C}$ and
$W=\what\W(M)$ for some supermanifold $M$.

\medskip

Let $G$ be a Lie supergroup acting on $V$. Then, we have a representation of  $G$
on $\CC^{\oo}_{V,c}(V_{\0})$ and thus a representation of $G$ on the
$\CC^{\oo}_{V,c}(V_{\0})$-module of compactlty supported distributions and
finally  on  $\CC^{-\oo}_{V}(V_{\0})$.

We denote by $\CC^{-\oo}_{V}(V_{\0})^G$ the set of the $G$-invariant generalized
functions on $V$.

\medskip

\subsection{Integration of pseudodifferential forms}\label{Int}

 Let $M$ be a supermanifold. Let $W$ be a Fr\'echet supervector space. Let 
 $\w\in\what\W(M,W)$ be a pseudodifferential form on $M$ with values in $W$.

We say that $\w$ is {\em  integrable}  if it is compactly supported on $M_\0$
and if rapidly decreasing along the fibres  as a function on the supervector
bundle $\what M\rightarrow M$.

We denote by $\what\W_{\int}(M,W)$ (resp. $\what\W_{\int}(M)$ when $W=\mathbb
R$) the set of integrable pseudodifferential forms on $M.$

\medskip

{\em Example:} Assume that $M$ is an affine supermanifold. Let $(x^{i},\xi^{j})$
be a coordinates system on $M$. Let $\w\in\what\W(M)$. We put:
\begin{displaymath}
\w=\som_{I,J}dx^{I}\xi^{J}\w_{I,J}(x^{i},d\xi^{j}).
\end{displaymath}
 Then $\w$ is integrable if it is a compactly supported function
in the $x^{i}$ and rapidly decreasing in the $d\xi^{j}$.

\medskip

Let $M$ be a globally oriented supermanifold of dimension $(m,n)$. Let $\w$ be
an integrable $W$-valued pseudodifferential form which support is included in a
coordinates domain $\UC$ with coordinates $(x^1,\dots,x^m,\xi^1,\dots,\xi^n)$. 
We put (cf. {BL77a,BL77b,Vor91}):
\begin{equation}
\int_M\w=\int_{\what M(\UC)}d_{(x,d\xi,dx,\xi)}\,\w.
\end{equation} 

It does not depend on the globally oriented coordinate system.

We extend this definition to a general integrable form on $M$ by means of a
partition of unity.

\subsection{Direct image of pseudodifferential forms}\label{ImDir} Let $M$ be a
supermanifold of dimension $(m,n)$ and $\pi:\VC\rightarrow M$ be a vector bundle
of rank $(k,l)$. 

\medskip

Let $W$ be a Fr\'echet supervector space. We say that
$\w\in\what\W_\VC(\VC_{\0},W)$ is {\em integrable along the fibres} if $\w$ is
rapidly decreasing along the fibres as a function on the supervector bundle
$\what\VC\rightarrow\what M$. 

We denote by $\what\W_{\pi_*}(\VC,W)$ (resp. $\what\W_{\pi_*}(\VC)$ when
$W=\mathbb R$) the set of pseudodifferential forms on $\VC$ that are integrable
along the fibres.

When $M$ and $\VC$ are globally oriented, and $\w\in\what\W_\VC(\VC,W)$ is
integrable along the fibres, we define the direct image $\pi_*\w$ of $\w$  as
the unique pseudodifferential form on $M$ such that for any 
pseudodifferential form $\a$ on $M$ which is compactly supported in $(\what
M)_{\0}$
we have (cf.
\cite{BGV92} with an other normalization):
\begin{equation}
\int_M(\pi_*\w)\a=\frac 1{(2\pi)^{\frac{k+l}2}}\int_\VC\w(\pi^*\a).
\end{equation}

 Let $\a\in\what\W_\VC(\VC,W)$ be integrable along the fibres and
$\be\in\what\W_M(M)$, then:
\begin{equation}
\pi_*(\a\pi^*(\be))=(\pi_*\a)\be.
\end{equation}

Indeed, since $\pi^*:\what\W_M(M)\to\what\W_\VC(\VC)$ is a
 morphism of superalgebras, if $\a\in\what\W_M(M)$ is compactly supported on
$(\what M)_{\0}$, 
we have:
\begin{equation}
\int_\VC\a(\pi^*\be)(\pi^*\ga)=\int_\VC\a\pi^*(\be\ga)
=(2\pi)^{\frac{k+l}2}\int_M(\pi_*\a)\be\ga.
\end{equation}

The application $\pi_*:\what\W_\VC(\VC,\AC)\longrightarrow\what\W_M(M,\AC)$ is a
morphism of $\what\W_M(M)$-modules with parity $k+l(mod\
\mathbb Z/2\mathbb Z)$. Indeed, the integration on $\VC$ is an operator of
parity $k+l+m+n(mod.\ \mathbb Z/2\mathbb 
Z)$ and the integration on $M$ has parity $m+n(mod.\
\mathbb Z/2\mathbb  Z)$.

We denote by $d_\VC$ (resp. $d_M$) the exterior differential on $\VC$ (resp.
$M$). Let $\w\in\what\W_{\VC}(\VC,\AC)$ be integrable along the fibres, then we
have:
 \begin{equation}
d_M\pi_*\w=(-1)^{k+l}\pi_*(d_\VC\w).
\end{equation}

\subsection{Integration along the fibres}\label{IntFib}

\subsubsection{Definition}

Let $\pi:\VC\rightarrow M$ be a supervector bundle  of rank $(k,l)$.

We call {\em volume form on the fibres $\VC$} the sections of the bundle  
\begin{displaymath}
 \vol (\VC)=\ber(\VC)\mult_{M}\VC\to \VC. 
\end{displaymath} 
We have for an open subset $\UC\subset M_{\0}$: 
\begin{displaymath}
\C_{\vol(\VC)}(\pi_{\0}^{-1}(\UC))=\C_{\ber(\VC)}(\UC)\tens_{\OC_{M}(\UC)}\OC_%
{\VC}(\pi_{\0}^{-1}(\UC)).
\end{displaymath} 
The canonical inclusion $M\hookrightarrow \VC$ by means of the
zero section (it corresponds to the canonical projection
$S(\C_{\VC^*})\to\OC_{M}$) induces a canonical inclusion
 \begin{displaymath}
\ber(\VC)=\ber(\VC)\mult_{M}M\hookrightarrow \vol(\VC).
\end{displaymath}
 We call the sections of $\ber(\VC)$ the {\em Berezinians volume forms along the
fibres of $\VC$.}

Let $\UC\subset M$ be an open subset. We say that a volume form $\DC\in \C_{\vol(\VC)}(\pi_{\0}^{-1}(UC))$ is rapidly decreasing along the fibres if 
 \begin{displaymath}
\DC\in\C_{\ber(\VC)}(\UC)\tens_{\OC_{M}(\UC)}\mathscr{S}_{\VC}(\pi_{\0}^{-1%
}(\UC)).
\end{displaymath}

\bigskip

Now we assume that $\VC\to M$ is an oriented supervector bundle.

\medskip

Let $\UC\subset M$ be an open subset. To any  volume form 
$\DC\in\C_{\vol(\VC)}(\pi_{\0}^{-1}(\UC))$ on the fibres of $\VC$ we
associate canonically  a distribution on $\VC$ defined on $\pi_{\0}^{-1}(\UC)$
with values in $\OC_{M}(\UC)$.

\medskip

Assume that $\UC$ is a trivialization subset of  $\VC$. Let $(e_{i},f_{j})$ be a
standard basis of $\C_{\VC}(\UC)$ (as $\OC_{M}(\UC)$-module). Let
$(x^{i},\xi^{j})$ be its dual basis.  
Let $V=\mathbb R e_{1}\oplus\dots\oplus\mathbb Re_{k}\oplus\mathbb R
f_{1}\oplus\dots\oplus\mathbb R f_{l}$ be the generic fibre of $\VC$.

We recall that we identify the Berezin integral $d_{(x,\xi)}$ on $V$ with a
section of $\ber(V)\subset\ber(\VC)\simeq\ber(V)\tens\OC_{M}(\UC)$. Then there exists
$h_{\UC}\in\OC_{\VC}(\pi_{\0}^{-1}(\UC))$ such that
$\DC|_{\UC}=d_{(x,\xi)}h_{\UC}$.

On the other hand we have
$\OC_{\VC}(\pi_{\0}^{-1}(\UC))=\CC^\oo_{V}(V_{\0},\OC_{M}(\UC))$.

Let $W$ be a Fr\'echet supervector space. Let 
\begin{displaymath}
\phi\in\CC_{V,c}(V_{\0},\OC_{M}(\UC)\tens W)\subset
\OC_{\VC}(\pi_{\0}^{-1}(\UC),W)
\end{displaymath}
 be a function on $\VC$ with values in $W$ that is  rapidly decreasing along the
fibres. We assume moreover that $\phi$ is  supported in $\UC$. Let $\PC$ be a near superalgebra. We have
$\VC_{\PC}(\pi_{\0}^{-1}(\UC))=V_{\PC}\mult M_{\PC}(\UC)$. For $v\in V_{\PC}$,
we have $\phi(v)\in\CC^{\oo}_{M_{\PC}}(M_{\PC}(\UC),W_{\PC})$ and  we write
$\phi(v,m)$ for its value at $m\in M_{\PC}(\UC)$.

\medskip 

We put for any near superalgebra $\PC$ and any $m\in M_{\PC}(\UC)$:
\begin{equation}
\label{eq:IntFib}
\Big(\int_{\VC/M}\DC\phi\Big)(m)=\int_{\VC/M}\DC_{m}(v) \phi(v,m) =
\int_{V}\, d_{(x,\xi)}(v)\,h_{\UC}(v,m) \phi(v,m)  \in\OC_{M}(\UC).
\end{equation} 
This formula defines $\int_{\VC/M}\DC \phi$ as a function in $\OC_{M}(\UC)$. The
formula  for change of variables shows that it does not depends on the choice of
the basis $(e_{i},f_{j})$.

\medskip

We extend it to any $\phi\in\OC_{\VC}(\VC_{\0})$ that is rapidly decreasing
along the fibres by means of a partition of unity on $M_{\0}$ associated with
trivialization subsets.  It is called the {\em integral along the fibres} of
$\phi$ against the volume form  $\DC.$

\medskip

In the same way we define the integral along the fibres of a rapidly decreasing
volume form $\DC$ and we denote it by $\int_{\VC/M}\DC\in \OC(M)$.

\medskip

{\em Example 1:} Let $M$ be a point and $\VC=\what V$ for some supervector
spaces $V$. Assume that $\VC$ is oriented (this means that $V$ is globally
oriented). Let $(x,\xi)$ be standard globally oriented coordinates on $V$. Then
$d_{(x,d\xi,dx,\xi)}$ defines a canonical volume form  on $\what V$. It is
canonical in the sense that it does not depends on the choice of standard
globally oriented coordinates $(x,\xi)$. Then, integration of integrable
pseudodifferential forms on $V$ is integration of functions on $\VC=\what V$
against this canonical volume form.

\medskip

{\em Example 2:} More generally, we can do the same construction when $M$ is not
trivial. We consider $\VC\to M$ a globally oriented supervector bundle and
$\what\VC\to\what M$ be defined as usual (it is an oriented supervector bundle).
We  replace $(x,\xi)$ by a standard globally oriented basis of sections of
$\VC^*\to M$. 
We obtain a canonical volume form  on the fibres of $\what \VC\to \what M$
(since $\VC$ is a globally oriented supervector bundle, $\ber(\what\VC)\to \what
M$ is trivial). We denote it $d_{(x,d\xi,dx,\xi)}$. Let $\pi:\VC\to M$ be the
projection of $\VC$ onto $M$. Then, the direct image $\pi_{*}\w$ of a
pseudodifferential form integrable along the fibres is the integration of $\w$
along the fibres of $\what\VC\to\what M$ against
$\frac1{(2\pi)^{\frac{k+l}2}}d_{(x,d\xi,dx,\xi)}$.

\subsubsection{Symplectic oriented case}

Assume that $\VC$ is a symplectic oriented supervector bundle. 

Let $\UC\subset M_{\0}$ be a trivialization subset and $V$ be the generic fibre.
Then $V$ is an oriented symplectic supervector space. In formula
(\ref{eq:IntFib}) we can assume that $(e_{i},f_{j})$  is a symplectic oriented
basis. In this case we take the Liouville integral $d_{V}=\frac1{(2\pi)^{\frac
n2}}d_{(x,\xi)}$ on $V$ (instead of $d_{(x,\xi)}$) which does not depend on the
choice of the symplectic oriented symplectic basis $(e_{i},f_{j})$. We can take
$h_{\UC}$ to be constant equal to $1$. Thus we obtain a canonical oriented
symplectic volume form on the fibres of $\VC$ called the {\em Liouville volume
form} and denoted by $\DC_{\VC}\in\C_{\ber(\VC)}(M_{\0})$.

\subsubsection{Volume forms and superconnection}

Let $\mathbb A$ be a superconnection on $\VC$.  We recall that $\mathbb A$
determines a superconnection $\mathbb A^*$ on $S(\VC^*)$ by formulas
(\ref{A*}) and (\ref{A*S}). Then, we extend $\mathbb A^*$ by   continuity to a
derivation of $\what\W_{M}(\UC,\OC(\VC))$.

\medskip

There is a superconnection $\mathbb A^{\ber}$ on $\ber(\VC)$ defined locally as
follows. Let $\UC$ be a trivialization set of $\VC$. Let $V$ be a generic fibre
of $\VC$. On $\UC$, we have $\mathbb A=d+\w$ where $d$ is the exterior
differential and $\w$ is an odd pseudodifferential form on $U$ with values in
$\g{gl}(V)$. Then we put $\mathbb A^{\ber}=d+str(\w)$. This does not depends on
the trivialization and thus defines a superconnection on $\ber(\VC)$.

\medskip
We extend it to  a superconnection $\mathbb A^{\vol}$ on $\vol(\VC)$ by the
following.

Let $\phi\in\what\W(M,\OC(\VC))=\OC(\VC)\what\tens_{\OC(M)}\what\W(M)$ and $\DC$
be a volume form on the fibres. 
Then if $\rk(\VC)=(k,l)$:
\begin{equation}
\mathbb A^{\vol}(\DC\phi)=(\mathbb A^{\vol}\DC)\phi
+(-1)^{k}\DC(\mathbb A^*\phi).
\end{equation} 
Since, locally, $\mathbb A^*$ is the sum of $d$ (the differential of
$\what\W(M)$) and of a derivation of $\OC(\VC)$ with 
coefficients in $\what\W(M)$) (in particular, they are constant in
the fibres), we have:
\begin{equation}
\int_{\VC/M}\mathbb A^{\vol}(\DC\phi)=d\int_{\VC/M}\DC\phi.
\end{equation}
 Moreover, if we suppose that $\mathbb A^{\vol}\DC=0$, then we have:
\begin{equation}
(-1)^k\int_{\VC/M}\DC(\mathbb A\phi) =d\int_{\VC/M}\phi\DC.
\end{equation}

\bigskip

Now, assume that $\VC$ is an oriented symplectic supervector space. Let $\mathbb
A$ be a superconnection on $\VC$ which leaves the symplectic structure
invariant. Let $\UC\subset M_{\0}$ be a trivialization subset of $\VC$. On
$\UC$, $\mathbb A=d+\w$ with $\w\in \what\W(\UC,\g{spo}(\VC))$. Let $\DC_{\VC}$
be the Liouville volume form on the fibres of $\VC$. Since  on
$\g{spo}(\C_{\VC}(\UC))$ the supertrace is zero, we have:
\begin{equation}
\mathbb A\DC_{\VC}=str(\w)\,\DC_{\VC}=0.
\end{equation}  
It follows, that for any
$\phi\in\what\W(M,\OC(\VC))=\OC(\VC)\mathop{\what\tens}\limits_{\OC(M)}\what%
\W(M)$ such that $\phi$ and $\mathbb A^*\phi$ are rapidly decreasing along the
fibres:
\begin{equation}
\label{eq:IntFibA}
\int_{\VC/M}\DC_{\VC}(\mathbb A^*\phi)=d\int_{\VC/M}\DC_{\VC}\phi.
\end{equation}

This result extends naturally to the equivariant situation defined in the next
section.

\subsection{Fourier transform}
\subsubsection{Fourier transform of smooth rapidly decreasing distributions}

Let $V$ be a supervector space of dimension $(n,m)$. Let $t$ be a smooth rapidly
decreasing distribution on $V$. We define the {\em Fourier transform of $t$} as
the function $\what t$ on $V^*$ defined  by:
\begin{equation}
\what t(h)=\int_{V}t(v)\exp(-\ii h(v)).
\end{equation}
This means that for any near superalgebra $\PC$ and any $h\in V_{\PC}^*$,
$\what t(h)$ is given by the integral on the right  hand side.

We stress that, this implies that $\what t$ is a rapidly decreasing function on
$V^*$.

\medskip

\subsubsection{Fourier transform of rapidly decreasing functions}

Let $V$ be a supervector space of dimension $(n,m)$. Let $\phi$ be a rapidly
decreasing function on $V$. We define the {\em Fourier transform of $\phi$} as the
distribution $\what \phi$ on $V^*$ defined  by the following. Let $(e_{i},f_{j})
$ be a standard basis of $V$. Let $(x^{i},\xi^{j})$ be its dual basis.  We put:
\begin{equation}
\what \phi(h)=\frac {(-1)^{\frac{n(n-1)}2}\ii^{n}}{(2\pi)^{m}}
d_{(e_{i},f_{j})}(h)\int_{V}d_{(x,\xi)}(v)\phi(v)\exp(\ii h(v)).
\end{equation}
The formula for change of variables (\ref{ChangVar}) implies that this formula
does not depends on the choice of $(e_{i},f_{j})$. 

Fourier inversion formula reads:
\begin{equation}
\what{\what \phi}(v)= \phi(v).
\end{equation}
(Cf. below the case of a
volume forms on the fibres of a supervector bundle for an other 
justification of the normalization)

\medskip

{\em Example: $dim(V)=(0,1)$.} Since the purely even case is well known we
examine the purely odd case. The simplest one is a one dimensional purely odd
vector space $V$. Let $(f)$ be a base of $V$. Let $(\xi)$ be its dual basis. Let
$\phi=a+\xi b\in\CC_{V}^{\oo}(V_{\0})$ be any function on $V$. Then, $\phi$ is
rapidly decreasing.

Let 
$h=-\xi f \in V^*\tens V$ be the generic point of $V^*$. Then 

\begin{equation}
\begin{split}
\what\phi=\what \phi(h)&=\ii d_{(f)}(h)\int_{V}d_{(\xi)}(v)\phi(v)\exp(\ii
h(v))\\
&= \ii d_{(f)}(h) \int_{V}d_{(\xi)}(v)\big((a+\xi b)(1-\ii \xi f)\big)(v)\\
& = \ii d_{(f)}(h)  (b-\ii af) (h)
\end{split}
\end{equation}

Let $v=f\xi\in V\tens V^*$ be the generic point of $V$. Then we have:

\begin{equation}
\begin{split}
\what{\what\phi}=\what {\what\phi}(v)&= \ii \int_{V^*} d_{(f)}(h)  (b-\ii af)
(h)\exp(-\ii v(h))\\
&= \ii \int_{V^*} d_{(f)}(h) \big ((b-\ii af) (1-\ii f\xi)\big)(h)\\
&= (a+b\xi)(v)=\phi(v)
\end{split}
\end{equation}

\medskip

\subsubsection{Fourier transform on the fibres of a supervector bundle}

Let $\VC\to M$ be a supervector bundle of rank $(k,l)$. Let $\DC$ be a rapidly
decreasing volume form along the fibres of $\VC$. We define the Fourier
transform of $\DC$ as the function $\what\DC$ on $\VC^*$ defined by:
\begin{equation}
\what\DC(h)=\int_{\VC/M}\DC(v)\exp(-\ii h(v)).
\end{equation}
We stress that this definition coincides with the precedent one if $M$ is a
point and thus $\VC$ is a supervector space.

\bigskip

The supervector space $\widetilde\VC=\VC^*\oplus\VC$ has a standard symplectic
structure $B$ defined as follows. Let $\UC\subset M_{\0}$ be a trivialization
subset of $\VC$. Let $(g_{i})$ be an homogeneous basis of sections of $\VC$ on
$\UC$. Let $(z^{i})$ be its dual basis. We put
$B(\C_{\VC},\C_{\VC})=B(\C_{\VC^*},\C_{\VC^*})=0$ and $B(z^{i},g_{j})=\d_{i}^j$. This
does not depend on the choice of $(g_{i})$ and thus defines $B$ globally.

Moreover we give to $\widetilde\VC$ the  orientation defined by the basis
$(z^{i},g_{i})$ (it is independent to the choice of the basis $(g_{i})$). We
recall that its dual basis is $((-1)^{p(g_{i})}g_{i},z^{i})$.

\medskip

Let $\DC_{\widetilde\VC}$ be the
canonical Liouville volume form along the fibres of $\widetilde\VC$. 
Let $\DC$ be a Berezin volume form on $\VC$. Since
$\ber(\widetilde\VC)=\ber(\VC^*)\tens_{M}\ber(\VC)$, there is an unique Berezin
volume form $\DC^*$ on $\VC^*$ such that $\DC^*\tens\DC=\DC_{\widetilde\VC}$.

\medskip

Let $\phi\in\mathscr{S_{\VC}}(\VC_{\0})$. We define the  Fourier transform of $\phi$ as the
volume form $\what\phi$ along the fibres of $\VC^*$ defined by:
\begin{equation}
\what\phi(h,m)= \DC_{m}^*(h)\int_{\VC/M}\DC_{m}(v)\phi(v,m)\exp(\ii h(v))
\end{equation}
Once more this definition does not depends on $\DC$ and coincides with the
preceding one when $M$ is a point.

\medskip We check this assertion in the fondamental cases of a $(1,0)$, $(0,1)$
and $(0,n)$ dimensional  supervector spaces.

\medskip

{\em Case of a $(1,0)$ dimensional vector space.} Let $(e)$ be a basis of $\VC$.
Let $(x)$ be its dual basis. Then $(x,e)$ is an oriented symplectic basis ot
$\widetilde\VC=\VC^*\oplus \VC$. Its dual basis is $(e,x)$, and
$\DC_{\widetilde\VC}=\frac 1{2\pi}d_{(e,x)}=\frac 1{2\pi}d_{(e)}d_{(x)}$. Thus
${d_{(x)}}^*=\frac 1{2\pi}d_{(e)}$.

\medskip

{\em Case of a $(0,1)$ dimensional vector space.} Let $(f)$ be a basis of $\VC$.
Let $(\xi)$ be its dual basis. Then $(\frac{f+\xi}{\sqrt 2},\ii
\frac{\xi-f}{\sqrt 2})$ is an oriented symplectic basis of
$\widetilde\VC=\VC^*\oplus \VC$. In fact it is a basis of $\widetilde
\VC\tens\mathbb C$. Its dual basis is $(\frac{\xi-f}{\sqrt 2},\ii
\frac{\xi+f}{\sqrt 2})$ and
 $\DC_{\widetilde\VC}= d_{(\frac{\xi-f}{\sqrt 2},\ii \frac{\xi+f}{\sqrt 2})}=\ii
d_{(f)}d_{(\xi)}$. Thus ${d_{(\xi)}}^*=\ii d_{(f)}$.

\medskip 

{\em Case of a $(0,n)$ dimensional vector space.} We have to be very carrefull
on orientation. Let $(f_{i})$ be a basis of $\VC$ with dual basis $(\xi^{j})$.
The symplectic form $B$ on the odd supervector space $\widetilde \VC$ is a
quadratic form of signature $(n,n)$. The orientation of the real supervector
space $\widetilde \VC$ corresponding to the symplectic basis
$(\frac{f+\xi}{\sqrt 2},\ii \frac{\xi-f}{\sqrt 2},\dots,\frac{f+\xi}{\sqrt
2},\ii \frac{\xi-f}{\sqrt 2})$ is given by:
\begin{equation}
(-\ii)^n(-\ii f_{1}\xi^{1})\dots(-\ii
f_{n}\xi^{n})=(-1)^{n}(-1)^{\frac{n(n-1)}2}f_{1}\dots f_{n}\xi^{1}\dots\xi^{n}.
\end{equation}
On the other hand, the orientation of $\widetilde \VC$ corresponding to
$(\xi^{1},\dots,\xi^{n},f_{1},\dots, f_{n})$ is given by:
\begin{equation}
(-1)^{n}f_{1}\dots f_{n}\xi^{1}\dots\xi^{n}.
\end{equation}
Thus we obtain:
\begin{equation}
{d_{(\xi^{j})}}^*=(-1)^{\frac{n(n-1)}2}\ii^n d_{(f_{j})}.
\end{equation}

\bigskip

Now the Fourier inversion formula reads:

\begin{equation}
\label{InvFour}
\what{\what \phi}= \phi.
\end{equation}

\bigskip

\section{Equivariant cohomology}

\subsection{Definitions} Let $M$ be a supermanifold and $G=(G_{\0},\g g)$ be a
supergroup acting on $M$ on the right. Let $\UC\subset\g g_{\0}$ be a $G_{\0}$
invariant open subset. Let $\OC_{\g g}\big(\UC,\what\W(M,\mathbb C)\big)$ be the
superalgebra of the functions on  $\UC$ with values in the complex
pseudodifferential forms on $M$:
\begin{equation}
\OC_{\g g}\big(\UC,\what\W(M,\mathbb C)\big)=\OC_{\g g\times\what M}(\UC\times
(\what M)_{\0},\mathbb C) =\OC_{\g
g}(\UC)\what\tens\what\W_{M}(M_{\0})\tens\mathbb C.
\end{equation}

We consider the subalgebra of the $G$-invariant elements of $\OC_{\g
g}\big(\UC,\what\W_{M} (M)_{\0},\mathbb C)\big)$. We denote  it by:
\begin{equation}
\what\W_G(\UC,M)=\bigl(\OC_{\g g}\big(\UC,\what\W(M,\mathbb C)\big)\bigr)^G.
\end{equation}

We call the elements of this algebra the {\em  equivariant forms}, on $M$.

We extend by $\mathbb C\tens\OC_{\g g}(\UC)$-linearity  and continuity the
exterior differential on $\what\W(M)$ to
 $\OC_{\g g}\big(\UC,\what\W(M,\mathbb C)\big)$. We recall that $d$ is an odd vector field on $\what M$.
 Weconsider it now as a vector field on $\g g\times \what M$. Let 
$(G_{i})$ be an homogeneous basis of $\g g$. We denote by $(g^{i})$ its dual
basis. We recall definition (\ref{X_{M}}) for $G_{iM}=(G_{i})_{M}$. We put:

\begin{equation}
\iota=\som_i\ \iota(G_{iM})\ g^{i},
\end{equation}
(cf. (\ref{contr})).  Then $\iota$ is a derivation of $\OC_{\g g\times\what
M}(\UC\times (\what M)_{\0},\mathbb C)$, called operator of contraction. It is
an odd vector field on $\g g\times \what M$ which satisfies $\iota^2=0$.

We define the equivariant differential on $\OC_{\g g\times\what M}(\UC\times
(\what M)_{\0},\mathbb C)$ by the odd vector field:
\begin{equation}
d_{\g g}=d-\ii\iota.
\end{equation} 
We have:

\begin{equation} 
d_{\g g}^2=-\ii\LC,
\end{equation} 
 with 

\begin{equation}
\LC=\som_i \LC(G_{iM})\ g^{i}.
\end{equation} 
The differentials $d$ and $\iota$ commute with the action of $G$ and so leave
$\what\W_G(\UC,M)$ stable. Therefore, $d_{\g g}$ induces a derivation of
$\what\W_G(U,M)$. Moreover, on $\what\W_G(U,M)$, we have $d_{\g g}^2=0$. 

We denote by $\what H_G(\UC,M)$ the cohomology of $\big(\what\W_G(\UC,M),d_{\g
g}\big)$ and we call it the equivariant cohomology of $M$ (with coefficients in
$\OC_{\g g}(\UC)^G$).

\medskip

We need variants of this definition. We define:

\medskip
\subsubsection{} The space $\what\W_{G,\int}(U,M)$ of {\em integrable
equivariant forms}:
\begin{equation}
\what\W_{G,\int}(U,M) 
= \OC_{\g g}\Big(\UC,\what\W_{\int}(M,\mathbb C)\Big)^G
\end{equation}
We have $\OC_{\g g}\Big(\UC,\what\W_{\int}(M,\mathbb C)\Big)\subset
\OC_{\g g}\Big(\UC,\what\W(M,\mathbb C)\Big)=\OC_{\g g\times\what M}(\UC\times
(\what M)_{\0},\mathbb C)$. The integrable equivariant forms are $G$-invariant
functions on the supervector bundle $\g g\times \what M\rightarrow \g g\times M$
which are rapidly decreasing in the fibres and compactly supported on $M_{\0}$.

\medskip
\subsubsection{}The space $\what\W_G^{-\oo}(\UC,M)$ of {\em equivariant forms
with generalized coefficients}: More precisely: 
\begin{equation}
\what\W^{-\oo}_{G}(\UC,M)=\CC^{-\oo}_{\g g}(\UC,\what\W(M,\mathbb C))^G.
\end{equation}

This is the space of $G$-invariant continuous homomorphisms from the space of
smooth compactly supported distributions on $\g g(\UC)$ with values in
$\what\W(M,\mathbb C)$. Let $\a\in\what\W^{-\oo}_{G}(\UC,M)$. Let $t$ be any
smooth compactly supported distribution $t$ on $\g g$  with support in $\UC$,
then:
\begin{equation}
\a(t)=(-1)^{p(t)p(\a)}\int_{\g g}t(X)\a(X)\,\in \what\W(M,\mathbb C).
\end{equation}

We recall that $\CC^{-\oo}(\UC)=\CC^{-\oo}_{\g g_{\0}}(\UC)\tens \L(\g g_{\1}^*)$. We put
 for any smooth compactly supported distribution $t$ on  $\UC$:
\begin{equation}
\a(t)=(-1)^{p(t)p(\a)}\int_{\g g_{\0}}t(X)\a(X)\,\in \L(\g
g_{\1}^*)\tens\what\W(M,\mathbb C).
\end{equation}

\medskip

{\it Example:} Let $M$ be a point, then:
\begin{equation}
\W_G^{-\oo}(\UC,M)=\CC^{-\oo}_{\g g}(\UC)^{G}=\Big(\CC^{-\oo}(\UC)\tens\L(\g
g_1^*)\Big)^G,
\end{equation}

\medskip

\subsubsection{}The space $\what\W^{-\oo}_{G,\int}(\UC,M)$ of {\em integrable
equivariant forms with generalized coefficients}:

More precisely: 
\begin{equation}
\what\W^{-\oo}_{G,\int}(\UC,M)=\CC^{-\oo}_{\g g}(\UC,\what\W_{\int}(M,\mathbb
C))^G.
\end{equation} 

This is the space of $G$-invariant continuous homomorphisms from the space of
smooth compactly supported distributions on $\g g(\UC)$ with values in
$\what\W_{\int} (M,\mathbb C)$.  

Let $\a\in\what\W^{-\oo}_{G,\int}(\UC,M)$. Let $t$ be any smooth compactly
supported distribution $t$ on $\g g$  with support in $\UC$, then:
\begin{equation}
\a(t)=(-1)^{p(t)p(\a)}\int_{\g g}t(X)\a(X)\,\in \what\W_{\int}(M,\mathbb C).
\end{equation}

As before, for any smooth compactly supported distribution $t$ on  $\UC$:
\begin{equation}
\a(t)=(-1)^{p(t)p(\a)}\int_{\g g_{\0}}t(X)\a(X)\,\in \L(\g
g_{\1}^*)\tens\what\W_{\int}(M,\mathbb C).
\end{equation}

\medskip

Let $\a\in\what\W^{-\oo}_{G,\int}(\UC,M)$. Then $\int_{M}\a\in\CC^{-\oo}_{\g
g}(\UC,\mathbb C)$ is  defined for any smooth compactly supported distribution
$t$ on  $\UC$ by the formula:

\begin{equation}
\Big(\int_{M}\a\Big)(t)=\int_{M}\big(\a(t)\big)=(-1)^{p(\a)p(t)}\int_{M}\int_{\g
g}t(X)\a(X).
\end{equation}

\medskip
\subsubsection{}Let $\pi:\VC\rightarrow M$ be a $G$ equivariant supervector
bundle. We will consider also the space $\what\W_{G,\pi_*}^{-\oo}(U,\VC)$ of
{\em equivariant forms on $\VC$ with generalized coefficients which are 
integrable along the fibres}:

More precisely:
\begin{equation}
\what\W^{-\oo}_{G,\pi_*}(\UC,\VC)=\CC^{-\oo}_{\g
g}(\UC,\what\W_{\pi_{*}}(\VC,\mathbb C))^G
\end{equation}

If $\a\in\what\W^{-\oo}_{G,\pi_*}(\UC,\VC)$,  for any  
compactly supported smooth distribution $t$ on $\g g_\0$ which support is
included in $\UC$:
\begin{equation}
(\pi_*\a)(t)=\pi_*\big(\a(t)\big)\in\what\W(M).
\end{equation}
 This defines $\pi_*\a$ as an element of $\what\W_{G}^{-\oo}(\UC,M)$.

\bigskip

The differential $d_{\g g}$ is defined on these spaces as well as the
corresponding cohomology spaces. If $M$ has a global orientation, we have (cf.
\cite{Lav98}): 
\begin{equation}
\int_Md_{\g g}\a=0,\text{ for }\a\in\what\W^{-\oo}_{G,\int}(\UC,M).
\end{equation}

\subsection{Some proprieties}

Let $G$ be a supergroup and  $M,N$ be two $G$-supermanifolds. Let
$\pi:N\rightarrow M$ be an equivariant morphism of $G$-supermanifolds. Let
$\a\in\what\W_{G}(M)$. We have $\pi^*\a\in\what\W_{G}(N)$ and:
\begin{equation}
\pi^*d_{\g g}\a=d_{\g g}\pi^*\a.
\end{equation}
Therefore $\pi^*$ induces an application in equivariant cohomology.

In particular, if $M$ is a point $H_G(\UC,M)=\OC_{\g g}(\UC)^G$ and
$\what H_G(\UC,N)$ is an $\OC_{\g g}(\UC)^G$-algebra (cf. for example
\cite{BGV92} in the classical situation).
\medskip

Moreover, if $\pi$ defines an equivariant supervector bundle and if $M$ and $N$
have a $G$-invariant global orientation,  if the fibres of $N$
are of dimension $(k,l)$, the following equation for  $\a\in\what\W_{G,\pi_{*}}(\VC)$:
\begin{equation}\label{eq:IntFibEq}
d_{\g g}\pi_*\a=(-1)^{(k+l)}\pi_*d_{\g g}\a.
\end{equation}

\medskip 

\subsection{Equivariant superconnection}\label{Econn}(cf
\cite{BGV92}) Let $G=(G_{\0},\g g)$ be a supergroup, $M$ a $G$-supermanifold and
$\VC\rightarrow M$ an equivariant supervector bundle. There is a natural
representation of $G$ in the space of pseudodifferential forms on $M$ and also
on  the space $\what\W(M,\VC)$ of pseudodifferential forms with value in $\VC$
(that is sections of $\VC\mult_{M}\what M\to \what M$). 
 We denote by $\LC^\VC$ the associated representation of $\g g$ in
$\what\W(M,\VC)$. This is a morphism of Lie superalgebras, from $\g g$  into
first order differential 
operators on $\what\W(M,\VC)$.

Let us consider equivariant forms on $M$ with values in $\VC$. For an open
$G_{\0}$-invariant subset $\UC\subset\g g_{\0}$ we put:
\begin{equation}
\what\W_G(\UC,M,\VC)=\OC_{\g g}\Big(\UC,\what\W(M,\VC\tens\mathbb C)\Big)^G.
\end{equation}

\medskip

Let $\mathbb A$ be a $G$-invariant superconnection on $\VC$. We extend it
$\OC_{\g g}(\UC)$-linearly to an endomorphism of $\OC_{\g
g}\Big(\UC,\what\W(M,\VC\tens\mathbb C)\Big)$.
 We extend linearly the contraction operator $\iota$ on $\OC_{\g
g}(\UC,\what\W(M,\mathbb C))$ to $\OC_{\g g}\Big(\UC,\what\W(M,\VC\tens\mathbb
C)\Big)$. We put:
 \begin{equation}
\mathbb A_{\g g}=\mathbb A-\ii\iota.
\end{equation} 
 It is a differential operator on $\OC_{\g g}\Big(\UC,\what\W(M,\VC\tens\mathbb
C)\Big)$ which leaves $\what\W_G(\g g,M,\VC)$ invariant. We call it the {\em
equivariant superconnection} associated with $\mathbb A$.

 For $\a\in\what\W_G(\g g,M)$ and $\w\in\what\W_G(\g g,M,\VC)$ non zero and
homogeneous we have:
\begin{equation}
\mathbb A_{\g g}(\w\a)=(\mathbb A_{\g g} \w)\a+(-1)^{p(\w)}\w(d_{\g g}\a).
\end{equation}

 The superconnection $\mathbb A_{\g g}$ 
acts on 
$\what\W_G(\UC,M, \g{gl}(\VC))$ by means of  the supercommutation bracket of 
the algebra of endomorphisms of $\what\W_G(\g g,M,\VC)$. 

Let  $\eta$  be the application which sends a form $\th\in\what\W_G(\g
g,M,\g{gl}(\VC))$ on the multiplication on the left by $\th$ in $\what\W_G(\g
g,M,\g{gl}(\VC))$. The {\em equivariant curvature} $F_{\g g}$ of $\mathbb A_{\g
g}$ is the element of $\what\W_G(\g g,M,\g{gl}(\VC))$ that satisfies:
\begin{equation}
\eta(F_{\g g})=\mathbb A_{\g g}^2+\ii\LC^\VC.
\end{equation}
The Bianchi identity is proved as in the classical case(cf {\cite{BGV92}
proposition 7.4} p.210): 
\begin{equation}
\mathbb A_{\g g}F_{\g g}=0.
\end{equation}
We define the {\em equivariant momentum} (for the $G$-invariant superconnection
$\mathbb A$) as:
\begin{equation}
\mu_{\mathbb A}=\LC^\VC-[\mathbb A,\iota].
\end{equation}
It is an element of $\what\W_G(\g g,M,\g{gl}(\VC))$ linear on $\g g$ which
satisfies:
\begin{equation}
F_{\g g}=\mathbb A^2+\ii\mu_{\mathbb A}=F+\ii\mu_{\mathbb A}.
\end{equation}

\section{Equivariant Thom form}\label{Thom.sec}

Let $G$ be a supergroup and $\pi:\VC\rightarrow M$ be a $G$-equivariant
supervector bundle. When $\VC$ is not purely even, there is no smooth
equivariant closed pseudodifferential form $\th$, 
integrable along the fibres and such that $\pi_*\th=1_M$ (cf.
\cite{Lav98}). Nevertheless, if the action of 
$G$  is sufficiently non-trivial, we shall show that there is a closed
equivariant form with generalized coefficients satisfying this propriety. 

Here is an example where such a form does not exists. Let $M$ be a point and
$V=\mathbb R^{(0,2)}$. We take $G=\{e\}$. A smooth equivariant form is a
function:
\begin{multline}
\w(\xi^1,\xi^2,d\xi^1,d\xi^2)=\w_{(0,0)}(d\xi^1,d\xi^2)\\
+\xi^1\w_{(1,0)}(d\xi^1,d\xi^2) +\xi^2\w_{(0,1)}(d\xi^1,d\xi^2)
+\xi^1\xi^2\w_{(1,1)}(d\xi^1,d\xi^2).
\end{multline} The condition $\w$ integrable means that $\w_{(0,0)},\w_{(1,0)},\w_{(0,1)}$ and
$\w_{{(1,1)}}$ are rapidly decreasing in $d\xi^1$ and $d\xi^2$. and the condition
$d\w=0$ implies that $\w_{(1,1)}=0$. Thus $\pi_{\*}\w=\int_V\w=0.$

\subsection{Preliminaries}
\subsubsection{Definitions}
\begin{defi} Let $G=(G_{\0},\g g)$ be a supergroup. Let $M$ be a
$G$-supermanifold. Let $\pi:\VC\rightarrow M$ be an equivariant supervector
bundle. We assume that $\VC$ and $M$ are globally oriented.

 An {\em equivariant Thom form} on $\VC$ is an equivariantly closed form
$\th\in\what\W_{G,\pi_*}^{-\oo}(\g g,\VC)$ which is integrable along the fibres
and such that $\pi_*\th=1.$
\end{defi}

This last equality means that if $t$ is a smooth compactly supported
distribution $\g g$, then
\begin{equation} 
(\pi_*\th)(t) =\int_{\g g}t(X)
\in \what\W(M,\mathbb C).
\end{equation}

Following Matha\"\i-Quillen (cf. \cite{MQ86}, see also \cite{BGV92} and
\cite{DV88}) we will construct such a Thom form for an Euclidean equivariant
supervector bundle $\VC\rightarrow M$ with a sufficiently non trivial action of
$G$.

\medskip

Let $Q$ be  a $G$-invariant Euclidean structure on $\VC$.

 We recall that $\OddId Q$ is symplectic form on $\Pi\VC$ (cf. equation
(\ref{eq:PiB})).

We recall that we denote by $\mu$ the moment application from $\g{spo}(\Pi\VC)$
to $S^2(\Pi\VC^*)\subset\OC(\Pi\VC)$ (cf. section \ref{sec:moment}).

We suppose that $\VC$ has a $G$-invariant superconnection $\mathbb A$ which
leaves  the Euclidean structure invariant. Therefore, the curvature $\mathbb
A^2$ is a pseudodifferential form on $M$ with values in $\g{osp}(\VC)$: $\mathbb
A^2\in\what\W(M,\g{osp}(\VC))$ and $\mathbb A^{\OddId
2}\in\what\W(M,\g{spo}(\Pi\VC))$ (cf. section \ref{sec:Superconn-Def}).

\subsubsection{First condition} We describe the first condition we need to
construct a Thom form.

We put $\mathbb A^{\OddId 2}=(\mathbb A^{\OddId })^2$.

The application $\mu$ (cf. formula (\ref{eq:moment})) sends $\g{spo}(\Pi\VC)$ on
$S^2(\Pi\VC^*)$. Therefore, $\mu(\mathbb A^{\OddId 2})$ is a function on $\what
M$ with values in $S^2(\Pi\VC^*)$. In other words, $\mu(\mathbb A^{\OddId 2})$
is      a function  on the bundle $\what
M\mathop\times\limits_M\Pi\VC\rightarrow \what M$ which is polynomial and
homogeneous of degree $2$ in the fibres. So 
$\mu(\mathbb A^{\OddId 2})_{\mathbb R}$ is a function on $(\what
M\mult_M\Pi\VC)_\0=(\what M)_\0\mathop\times\limits_{M_\0}(\Pi\VC)_\0$. The
condition is:
\begin{equation*} (*)\qquad\mu(\mathbb A^{\OddId 2})_{\mathbb R}\le 0.
\end{equation*}

\medskip

Let $v\in\C_{\VC}(M_{\0})_{\1}$ be an odd section of $\VC$. Then $\OddId v$ is
an even section of $\Pi\VC$. We have:  
\begin{equation}
\mu(\mathbb A^{\OddId 2})(\OddId v)=-\frac 12 \OddId Q(\OddId v,\mathbb
A^{\OddId 2}\OddId v)=\frac 12 Q(v,\mathbb A^2 v)\in\what\W(M)=\OC_{\what
M}((\what M)_{\0}).
\end{equation}
 Condition $(*)$ means that for any $v\in\C_{\VC}(M_{\0})_{\1}$, the restriction
of $\frac12Q(v,\mathbb A^2 v)$ to $(\what M)_0$, is non positive.

\medskip

This condition is trivially satisfied when $\mathbb A^2=0$. Moreover, this
condition depends of the choice of the particular superconnection we choose.

\subsubsection{More notations} We denote by $\mathbb A_{\g g}$ the equivariant
superconnection associated with $\mathbb A$. We denote by $\pi^*S(\Pi \VC^*)$ the
bundle with $\VC$ as base space which is  the pullback of $S(\Pi \VC^*)$ by
$\pi$:  $\pi^*S(\Pi \VC^*)=S(\Pi \VC^*)\mathop\times\limits_{M}\VC$.  

Since $\Pi\VC^*=S^1(\Pi\VC^*)$, there is a natural injection of $\pi^*\Pi\VC^*$
in $\pi^*S(\Pi\VC^*)$. 

To help the reader we fix (abuse of) notations in the following commutative
diagram:

\begin{equation}
\begin{CD} 
p^*\pi^*(\Pi\VC)=\Pi\VC\mult_{M}\what\VC @>>>\what \VC \cr
@VVV  @VV{p}V \cr
\pi^*(\Pi\VC)=\Pi\VC\mult_{M}\VC@>>>  \VC \cr
@VV{\pi}V  @VV{\pi}V \cr
\Pi\VC @>>> M
\end{CD}
\end{equation}

\medskip

The left vertical projections correspond to the  right vertical ones with the
same letter and thus ``act only on the second factor''.
 In the preceding diagram we can change $\Pi\VC$ into $\Pi\VC^*$ without any
other changes.

\medskip

We denote by $\bv$ be the tautological section with change of parity of
$\Pi\VC\mathop\times\limits_{M}\VC\to\VC$. Let $\UC\subset M_{\0}$ be a
trivialization  subset of $\VC$. Let $(e_{i})_{i\in I}$ be an homogeneous basis
of sections of $\VC$. Let $(x^{i})$ be its dual basis. Then,  $(\OddId e_{i})$ 
the basis of sections of $\Pi\VC$. By abuse of notations we also denote by  
$(\OddId e_{i})$  the basis of sections  $(\pi^*\OddId e_{i})$ of
$\Pi\VC\mult_{M}\VC$. We have:
\begin{equation}
\bv=\som_{i}(\OddId e_{i})
x^{i}\in\C_{\Pi\VC}(\UC)\tens_{\OC_{M}(\UC)}\C_{\VC^*}(\UC)\subset
\C_{\pi^*\Pi\VC}(\pi_{\0}^{-1}(\UC)).
\end{equation}
Since the above formula does not depend on the choice of $(e_{i})_{i\in I}$,
this defines a global section $\bv\in\C_{\pi^*\Pi \VC}(\VC_{\0})$.

 We recall (cf. section \ref{sec:bilin} for the definition of $\OddId Q^*$)
that:
\begin{equation}
\OddId Q^*(\bv)\in\C_{\pi^*\Pi\VC^*}(\VC_{\0}).
\end{equation}

Let $\PC$ be a near superalgebra. We  explicit $\OddId Q^*(\bv)$ for
$\PC$-points. For $w\in(\pi^*\Pi\VC)_{\PC}$, we recall that $\pi(w)$ is its
projection on $(\Pi\VC)_{\PC}$. For any  $w\in(\pi^*\Pi\VC)_{\PC}$ and any
 $v\in \VC_{\PC}$ we have:

\begin{equation}
\OddId Q^*(\bv(v))(w)=\OddId Q(\OddId v,\pi(w)).
\end{equation}
In particular $\OddId Q^*(\bv)(\bv)=\OddId Q(\bv,\bv)$.

\bigskip

We recall (cf. section \ref{IndConn}) that $\mathbb A$ determines a
superconnection $\mathbb A^{\OddId}$ on $\Pi\VC$ and $\mathbb A^{\OddId*}$ on $S(\Pi\VC^*)$. This determines  equivariant
superconnections $\mathbb A_{\g g}^{\OddId}$ on $\what\W_G(\g g_{\0},M,\Pi\VC)$
and $\mathbb A_{\g g}^{\OddId*}$ on the algebra $\what\W_G(\g
g_{\0},M,S(\Pi\VC^*))$.

\medskip

The differential operator $\pi^*\mathbb A_{\g g}^{\OddId}$  on $\what\W_G(\g
g,\VC, \pi^*\Pi\VC)$ is an equivariant superconnection defined for $\w\in\what\W_G(\g
g,M,\Pi\VC)$ by:

\begin{equation}
\pi^*\mathbb A_{\g g}^{\OddId}(\pi^*\w)=
\pi^*(\mathbb A_{\g g}^{\OddId}\w).
\end{equation}

 We have:
\begin{equation}
\what\W_G(\g g_{\0},\VC,\pi^*\Pi\VC)=\what\W_G(\g
g_{\0},M,\Pi\VC)\tens_{\what\W_{G}(M)} \what\W_{G}(\VC).
\end{equation}
With respect to this decomposition we have:
\begin{equation}\label{eq:ConnTens}
\pi^*\mathbb A_{\g g}^{\Pi}=\mathbb A_{\g g}^{\Pi}\tens 1+1\tens d_{\g g},
\end{equation}
where $(1\tens d_{\g g})(\a\tens\be)=(-1)^{p(\a)}\a\tens d_{\g g}\be$ and
$( \mathbb A^{\OddId}_{\g g}\tens 1)(\a\tens\be)=\mathbb
A^{\OddId}_{\g g}\a\tens  \be$.

\medskip

We define similarly the equivariant connection $\pi^*\mathbb A_{\g g}^{\OddId*}$
 on $\what\W_G(\g g,\VC, \pi^*S(\Pi\VC^*))$.

\medskip

Now, 
\begin{displaymath}
\what\W_G(\g g_{\0},\VC,\pi^*S(\Pi\VC^*))\subset \OC_{\g g}\Big(\g g_{\0}
,\OC\big(\Pi\VC\mult_{M}\what\VC\big)\Big)^G
\end{displaymath} and the inclusion is dense. The operator  $\pi^*\mathbb A_{\g
g}^{\OddId*}$  can be extended by continuity to this last space. 

We stress  that, if $\w\in\what\W_G(\g g_{\0},\VC,\pi^*S(\Pi\VC^*))_{\0}$ and
$f$ is an entire function of the complex variable, the function $f(\w)$ is well
defined in $\OC_{\g g}\Big(\g g_{\0}
,\OC\big(\Pi\VC\mult_{M}\what\VC\big)\Big)^G $.
\medskip

Finaly, let $v\in\C_{\pi^*\Pi\VC}(\VC_{\0})$ be non zero and homogeneous. We
denote by $\de_{v}$ the derivation of $\C_{\pi^*S(\Pi\VC^*)}(\VC_{\0})$ such
that for $\phi\in\C_{\pi^*\Pi\VC^*}(\VC_{\0})$ non zero and homogeneous:
\begin{equation}
\de_{v}\phi=(-1)^{p(v)p(\phi)}\phi(v).
\end{equation}

We extend it $\what\W_{G}(\g g_{\0},\VC)$-linearly to a derivation of
$\what\W_{G}(\g g_{\0},M,\pi^*S(\Pi\VC^*))$ and then by continuity to a
derivation of $\OC_{\g g}\Big(\g g_{\0}
,\OC\big(\Pi\VC\mult_{M}\what\VC\big)\Big)^G $.

\bigskip

We denote by $v^{\OddId}$ the generic point of the $\OC( \VC)$-module
 $\C_{\Pi\VC\mult_{M}\VC}(\VC_{\0})$. Let $((\OddId e_{i})^*)$ be the
dual basis of $(\OddId e_{i})$. We have:
\begin{equation}
\label{eq:v} v^{\OddId}=\som_{i}\OddId e_{i}(\OddId e_{i})^*
\end{equation}
We have:
\begin{equation}
\de_{\bv}v^{\OddId}=\bv.
\end{equation}

\bigskip

\subsubsection{The form $\w_\VC$} 
The operator
\begin{equation}
\mathbb B=\pi^*\mathbb A^{\OddId *}_{\g g}+\ii\de_{\bv}
\end{equation}
 defines an equivariant superconnection on $\pi^*S(\Pi\VC^*)$.
 
We denote by $F_{\g g}^{\OddId}$ the equivariant curvature of $\mathbb A_{\g
g}^{\Pi}$.

 Let $\mu_{\mathbb A}$ be the equivariant moment of $\mathbb A_{\g g}$:
\begin{displaymath}
\mu_{\mathbb A}\in\what\W_G(\g g,M,\g{osp}(\VC)).
\end{displaymath} Let $\mu_{\mathbb A^{\OddId}}$ be the equivariant moment of
$\mathbb A^{\OddId}$. We recall that we identified $\what\W_G(\g
g,M,\g{osp}(\VC))$ and $\what\W_G(\g g,M,\g{spo}(\Pi\VC))$. Under this
identification we have $\mu_{\mathbb A^{\OddId}}=\mu_{\mathbb A}$.

\medskip

We put:
\begin{equation}
\begin{split}
\w_{\mathcal V}&=\frac12 \OddId Q(\mathbb B v^{\OddId},\mathbb B v^{\OddId})
+\ii\mu(\mu_{\mathbb A}(X))\\
&=-\frac12\OddId Q(\bv,\bv)+\ii \pi^*\mathbb A_{\g g}^{\OddId*}\OddId
Q^*(\bv)+\mu(F_{\g g}^{\OddId})\in\what\W_G(\g g,\VC,\pi^*S(\Pi\VC^*)).
\end{split}
\end{equation}

\bigskip

\begin{prop} 
\begin{equation}
\label{Bian}
\big(\pi^*\mathbb A_{\g g}^{\OddId*}+\ii\de_{\bv}\big)\w_{\mathcal V}=0,
\end{equation} 

 If $f$ is an entire function of one complex variable, we have:
\begin{equation}
  (\pi^*\mathbb  A^{\OddId*}_{\g g}+\ii\de_{\bv})f(\w_{\mathcal V})=0.
\end{equation} 

\end{prop}

\begin{proof}[Proof] Since $\OddId Q$ is antisymmetric we have:
\begin{equation}
\mathbb B\OddId Q(\mathbb B v^{\OddId},\mathbb B v^{\OddId})= \OddId Q(\mathbb
B^2 v^{\OddId},\mathbb B v^{\OddId})+ \OddId Q(\mathbb B v^{\OddId},\mathbb B^2
v^{\OddId})=0.
\end{equation}
On the other hand, since $Q$ is $\mathbb A$-invariant:
\begin{equation}
\pi^*\mathbb A^{\OddId*} \mu(\mu_{\mathbb A}(X))= -\frac12\pi^*\mathbb
A^{\OddId*} \OddId Q(v^{\OddId},\mu_{\mathbb A}(X) v^{\OddId})=0.
\end{equation}
Finally, since $\OddId Q$ is antisymmetric and $\mu_{\mathbb
A}(X)\in\C_{\g{spo}(\Pi\VC)}(M_{\0})$:
\begin{equation}
\de_{\bv} \mu(\mu_{\mathbb A}(X))= -\frac12\OddId Q(\bv,\mu_{\mathbb A}(X)
v^{\OddId}) -\frac12\OddId Q(v^{\OddId},\mu_{\mathbb A}(X) \bv)=0.
\end{equation}
Hence $\mathbb B\mu(\mu_{\mathbb A}(X))=0$ and equality (\ref{Bian}) follows. 

 The last equality  follows immediately, since $\pi^*\mathbb A_{\g
g}^{\OddId*}+\ii\de_{\bv}$ is a derivation of the superalgebra $\OC_{\g g}(\g
g_{\0},\OC(\Pi\VC\mult_{M}\what\VC))^G.$
\end{proof}

\medskip

\subsubsection{An equivariantly closed form} Since $Q$ is an euclidean structure
on $\VC$, $\OddId Q$ is a symplectic structure on $\Pi\VC$. Moreover, since
$\VC$ is globally oriented, $\Pi\VC$ is an oriented symplectic supervector
bundle in sense of section \ref{sec:SympOrBun}. Let $\DC_{\Pi\VC}$ be the
Liouville
  volume form along the fibres of $\Pi\VC\mult_{M}\what\VC\to \what\VC$. We
denote by 
\begin{displaymath}
\phi\mapsto T(\phi)=\int_{\Pi\VC\mult_{M}\what\VC/\what\VC}\DC_{\Pi\VC}\phi
\end{displaymath}
 the corresponding  distribution (cf. section \ref{IntFib}).
 
 Recall that $\w_{\VC}\in \what\W_{G}(\g g,\VC,\pi^*S(\Pi\VC^*))\subset\OC(\g g,\Pi\VC\mult_{M}\what\VC)$.
 Thus, if all the functions involved are integrable along the fibres of 
 $\OddId\VC\times_{M}\what\VC$, we obtain:
\begin{equation}
(-1)^k T\Bigl((\pi^*\mathbb  A^{\OddId*}_{\g g}+\ii\de_{\bv})f(\w_{\mathcal
V})\Bigr)=d_{\g g}T\bigl((f(\s^*(\w_{\mathcal V}))\bigr),
\end{equation}
where $(k,l)$ is the dimension of the fibres of $\VC$. This follows from formula
(\ref{eq:IntFibEq}) and  from the fact that since $\bv$ is a section of $\pi^*\Pi
V$, $\de_{\bv}$ is a derivation along the fibres of
$\Pi\VC\mult_{M}\what\VC\to\what\VC$ with constant coefficients in the direction
of the fibres     and thus $T\big(\de_{\bv}f(\w_{\mathcal V})\big)=0$.

Since $(\pi^*\mathbb A^{\OddId*}_{\g g}+\ii\de_{\bv})f(\w_\VC)=0$, it follows
that for any  entire function $f$ of one complex variable such that
$f(\w_{\VC})$ is rapidly decreasing along the fibres of
$\Pi\VC\mult_{M}\what\VC\to\what\VC$:
\begin{equation} 
d_{\g g} T\Bigl(f\bigl(\w_{\mathcal V}\bigr)\Bigr)=0.
\end{equation}

Therefore the form $T\bigl(f(\w_{\mathcal V})\bigr)$ is  an equivariantly closed
form on $\VC$. To find a Thom form  we just have to find a ``good'' function
$f$. The choice $f(z)=\exp(z)$ gives, up to a multiplicative constant, a Thom
form on $\VC$.

Moreover, the equivariant cohomology class of this form  in $\what
H_{G\int}^{\oo}(\VC)$ does not depends on the choice of the $G$-invariant
superconnection $\mathbb A$. Let $\w'\in\what\W(M,\g{osp}(\VC))_\1^G$. We put
$\mathbb A(t)=\mathbb A+t\w'$  and denote  its equivariant curvature by
$F(t)_{\g g}$. We put
\begin{equation}
\w_\VC(t)=-\frac12\OddId Q(\bv,\bv)+\ii\pi^*\mathbb A(t)^{\OddId*}_{\g g}
\OddId Q^*(\bv)+\mu(F(t)_\g g).
\end{equation} 
We have:
\begin{equation}
\begin{split}
\frac d{dt}\w_{\VC}(t)&=\ii\,(\pi^*{\w'}^*)\OddId Q^*(\bv)+\t\big([\mathbb
A(t)_{\g g},\w']\big)\\
&=\big(\pi^*\mathbb A(t)^{\OddId*}_{\g g}+\ii\de_{\bv}\big)\mu(\w').
\end{split}
\end{equation}
and since $(\pi^*\mathbb A(t)^{\OddId *}+\de_{\bv})\w_{\VC}(t)=0$:
\begin{equation}
\frac d{dt}\exp(\w_\VC(t))=
\big(\pi^*\mathbb A(t)^{\OddId*}_{\g g}+\ii\de_{\bv}\big)\Big(\mu(\w')
\exp(\w_\VC(t))\Big)
\end{equation}
Thus, since $\w_\VC(0)$ and $\w_\VC(1)$ corresponds to the superconnections
$\mathbb A$ and $\mathbb A+\w'$:

\begin{equation}
\exp(\w_{\VC}(1))-\exp(\w_{\VC}(0)) =\int_0^1 \big(\pi^*\mathbb
A(t)^{\OddId*}_{\g g}+\ii\de_{\bv}\big)\Big(\mu(\w')
\exp(\w_\VC(t))\Big)dt,
\end{equation}
and:
\begin{multline} T\Big(\exp(\w_\VC(1))-\exp(\w'_\VC(0))\Big)\\
=T\Big(\int_0^1 \big(\pi^*\mathbb A(t)^{\OddId*}_{\g
g}+\ii\de_{\bv}\big)\big(\mu(\w')
\exp(\w_\VC(t))\big)dt\Big)\\
=\int_0^1 T\Big(\big(\pi^*\mathbb A(t)^{\OddId*}_{\g
g}+\ii\de_{\bv}\big)\Big(\mu(\w')
\exp(\w_\VC(t))\Big)\Big)dt\\
=d_{\g g}\int_0^1 T\big(\mu(\w')
\exp(\w_\VC(t))
\big)dt.
\end{multline}

\subsubsection{Two additional notations}

Let $\UC\subset M_{\0}$ be an open subset. 
We put

\medskip

\begin{equation}
\label{U0m} U_{+}^{\VC_{\1}}(\UC)=\Big\{X\in\g g_\0\,\Big/\,\big(\forall
v\in\C_{\VC}(\UC)_\1\,/ \,
\forall x\in \UC,  v_{\mathbb R}(x)\not=0\big),\ 
 Q(v,\mu_{\mathbb A}(X)v)_{\mathbb R}>0\Big\}.
\end{equation} 

\medskip

 For $X\in\g g_\0$,  and $v\in\C_{\VC}(\UC)$, we have $Q(v,\mu_{\mathbb
A}(X)v)\in\what\W_{M}(\UC)=\OC_{\what M}((\what M(\UC))_{\0})$. Thus
 $Q(v,\mu_{\mathbb A}(X)v)_{\mathbb R}$ is a function on $(\what M(\UC))_{\0}$.
Thus $Q(v,\mu_{\mathbb A}(X)v)>0$ means
\begin{equation}
\forall {\bm}\in (\what M(\UC))_{\0},\, Q(v,\mu_{\mathbb A}(X)v)_{\mathbb
R}({\bm}) =Q_{x_{\bm}}(v_{\mathbb R}(x_{\bm})v,\mu_{\mathbb A}(X)(\bm)v_{\mathbb
R}(x_{\bm}))>0.
\end{equation}
where $x_{\bm}$ is the projection of $\bm\in (\what M(\UC))_{\0}$ on $\UC$.

 Since $v$ is supposed to be odd, $v_{\mathbb R}(x_{\bm})\in
(\VC_{x_{\bm}})_{\1}$. On the other hand, $Q_{x_{\bm}}|_{(\VC_{x_{\bm}})_{\1}}$
is symplectic, thus $w\mapsto Q_{x_{\bm}}(w,\mu_{\mathbb A}(X)(\bm)w)$ is a
quadratic form on $(\VC_{\bm})_{\1}$ and the condition says that this quadratic
form is positive definite for any ${\bm}\in(\what M(\UC))_{\0}$.

\medskip

On the other hand we put:

\medskip

\begin{equation}
\label{U} U^{\VC_{\1}}(\UC)=\text{Interior of }\Big\{X\in\g g_0\,\Big/
\forall \bm\in (\what M(\UC))_{\0},
\mu_{\mathbb A}(X)(\bm)\big|_{(\VC_{x_{\bm}})_{\1} }\text{ is invertible}\Big\}.
\end{equation}

\bigskip

We describe  these subsets in the particular case when $M$ is a point and
$\mathbb A$ is the trivial superconnection $d$.

In this case $\VC$ is a Euclidean supervector space with  a representation
$\r=(\r_{\0},\r)$ of $G=(G_{\0},\g g)$. In particular $\r:\g
g\to\g{osp}(\VC,Q)$. Then for $X\in\g g_{\0}$, $\mu_{\mathbb A}(X)=\r(X)$. We
have (since $M$ is a point we omit $\UC$):
\begin{equation}
U^{\VC_{\1}}=\big\{X\in\g g_{\0}\big/\r(X)|_{\VC_{\1}}\text{ is invertible.}\big\}
\end{equation}
It is an open subset of $\g g_{\0}$.

 Since $(\VC_{\1},Q|_{\VC_{\1}})$ is, as an ungraded vector space, a symplectic
vector space, $w\mapsto B(w,\r(X)w)$ is a quadratic form on (the ungraded vector
space) $\VC_{\1}$. Thus $U_{+}^{\VC_{\1}}$ is the open subset of $X\in\g g_{\0}$ such
that this quadratic form is positive definite; and $U^{\VC_{\1}}$ is the open subset
of $X\in\g g_{\0}$ such that this quadratic form is non degenerate.

\subsection{Construction of an equivariant Thom form}

We put:
\begin{equation}
\th=T\Bigl(\exp(\w_{\mathcal V})\Bigr).
\end{equation} 

\begin{theo}\label{Thom} Let $G=(G_{\0},\g g)$ be a supergroup. Let $M$ be a
$G$-supermanifold. Let $\pi:\VC\rightarrow M$ be an equivariant supervector
bundle on $M$ with rank $(k,l)$. We assume that $\VC$ and $M$ are globally
oriented, and that $\VC$ is endowed with an Euclidean structure denoted by
$Q(.,.)$, all these structures being $G$-invariants.

Finally we assume that there is a $G$-invariant superconnection which leaves the
Euclidean structure invariant. We assume:
\begin{itemize}
\item[$(*)$] $\mu(\mathbb A^{\OddId2})_{\mathbb R}\le 0$ on $(\what M\mult_M\Pi\VC)_\0$.
\item[$(**)$] There is a covering of $M_{\0}$ by open subsets $\UC$  such that
$U_{+}^{\VC_{\1}}(\UC)$ contains a non empty open subset (cf. formula (\ref{U0m})).
\end{itemize} Then, the equivariant form $\th$ defined above is a Thom form. It
does not depends on the choice of the superconnection $\mathbb A$.

Moreover, let $\UC$ be an open subset  of $M_{\0}$. Then the restriction of
$\th$ to $U^{\VC_{\1}}(\UC)\mult\VC(\UC)$ is a smooth equivariant form.

\end{theo}

\begin{proof}[Proof] Up to the multiplicative constant we  have to check that
the form $T\bigl(\exp(\w_{\mathcal V})\bigr)$ has the required proprieties. We
already know that  it is an equivariantly closed form, so  we have to check that
the form $T\bigl(\exp(\w_{\mathcal V})\bigr)$ is a well defined equivariant form
with generalized coefficients, integrable along the fibres and to evaluate its
integral.

\subsubsection{A form with generalized coefficients}\label{sec:ThomGene}

First we have to show that $T\bigl(\exp(\w_{\mathcal V})\bigr)$ is a well
defined generalized function on $\g g$.

It means that for any smooth compactly supported distribution $t$ on $\g g$ the
integral
\begin{equation}
\label{I}
\int_{\g g}t(X)\exp(\w_{\VC}(X))\in\OC(\Pi\VC\mult_{M}\what\VC).
\end{equation}
 is a rapidly decreasing function along the fibres of
$\Pi\VC\mult_{M}\what\VC\to \what\VC$.
 
In this case, by definition, we put:
 \begin{equation}
<\th,t> =T\Big(\int_{\g g}t(X)\exp(\w_{\VC}(X))\Big).
\end{equation}

\medskip

Since $F_{\g g}^{\OddId }(X)=\mathbb A^{\OddId 2}+\ii\mu_{\mathbb A}(X)$ and
$\pi^*\mathbb A_{\g g}^{\OddId*}(X)\OddId Q^*(\bv)=\pi^*\mathbb
A^{\OddId*}\OddId Q^*(\bv)$, we have:
\begin{equation}
\w_{\VC}(X)=-\frac12\OddId Q(\bv,\bv)+\ii\pi^*\mathbb A^{\OddId*}\OddId
Q^*(\bv)+\mu(\mathbb A^{\OddId 2})+\ii\mu(\mu_{\mathbb A}(X)).
\end{equation}
Thus the preceding integral (\ref{I}) is equal to:
\begin{equation}
\exp\big(-\frac12\OddId Q(\bv,\bv)+\ii\pi^*\mathbb A^{\OddId *}\OddId
Q^*(\bv)+\mu(\mathbb A^{\OddId 2})\big)
\int_{\g g}t(X)\exp(\ii\mu(\mu_{\mathbb A}(X))).
\end{equation}

\medskip

Since the problem is entirely  local on $M$, we restrict us to a trivialization
subset $\UC$ of $\VC$. Moreover hypotheses $(**)$ allows us to assume that $
U_{+}^{\VC_{\1}}(\UC)$ contains a non empty open subset.  To avoid boring notations
we assume that $\VC$ is trivial and $\UC=M_{\0}$.

 We have to show that $\int_{\g g}t(X)\exp(\ii\mu(\mu_{\mathbb A}(X)))$ is
rapidly decreasing along the fibres of $\Pi\VC\mult_{M}\what\VC\to\what\VC$. We
stress that this problem  is completely similar to the one of sections 3.1 and
3.7 of \cite{Lav03}. We reproduce here the
argument.

\medskip

Let $(e_{i},f_{j})$ be a standard basis of sections of
$\VC$ on $M$. Then $(\OddId f_{j},\OddId e_{i})$ is a standard basis of sections
of $\Pi\VC$. By abuse of notations we also denote by   $(\OddId f_{j},\OddId
e_{i})$  the basis of sections  $(\pi^*\OddId f_{j},\pi^*\OddId e_{i})$ of
$\Pi\VC\mult_{M}\what\VC$.  Let  $((\OddId f_{j})^*,(\OddId e_{i})^*)$ be its
dual basis. We have
\begin{equation}
\label{eq:v2} v^{\OddId }=\som_{j}\OddId f_{j}(\OddId f_{j})^*+\som_{i}\OddId
e_{i}(\OddId e_{i})^*.
\end{equation}
 We put $v_{\0}=\som_{j}\OddId f_{j}(\OddId f_{j})^*$ and $v_{\1}=\som_{i}\OddId
e_{i}(\OddId e_{i})^*$.

We have
\begin{multline}
\mu(\mu_{\mathbb A}(X))=-\frac12\OddId Q(v^{\OddId},\mu_{\mathbb
A}(X)v^{\OddId})\\
=-\frac12
\Big(\OddId Q(v_{\0},\mu_{\mathbb A}(X)v_{\0})+\OddId Q(v_{\1},\mu_{\mathbb
A}(X)v_{\1})+2\OddId Q(v_{\0},\mu_{\mathbb A}(X)v_{\1})\Big).
\end{multline} Thus, 
\begin{multline}
\exp(\mu(\mu_{\mathbb A}(X)))=\exp\Big(-\frac12
\Big(\OddId Q(v_{\1},\mu_{\mathbb A}(X)v_{\1})+2\OddId Q(v_{\0},\mu_{\mathbb
A}(X)v_{\1})\Big)\Big)\\
\exp\Big(-\frac12
\big(\OddId Q(v_{\0},\mu_{\mathbb A}(X)v_{\0})\big)\Big).
\end{multline}

\medskip

Since $v_{\1}$ is a linear combinaison of nilpotent elements, 
\begin{displaymath}
\exp\Big( -\ii\frac12\big(\OddId Q(v_{\1},\mu_{\mathbb A}(X)v_{\1})+2\OddId
Q(v_{\0},\mu_{\mathbb A}(X)v_{\1})\big)\Big)
\end{displaymath}
 is a polynomial function on $\g g$ with values in $\C_{S(\Pi\VC^*)\mult_{M}
\what\VC}(\what\VC_{\0})$.  

Now we put 
\begin{displaymath}
\r(X)=t(X)\exp\Big( -\ii\frac12\big(\OddId Q(v_{\1},\mu_{\mathbb
A}(X)v_{\1})+2\OddId Q(v_{\0},\mu_{\mathbb A}(X)v_{\1})\big)\Big);
\end{displaymath} It is a a smooth compactly supported distribution on $\g g$
with values in $\C_{S(\Pi\VC^*)\mult_{M} \what\VC}(\what\VC_{\0})$.

Thus:

\begin{equation}
\int_{\g g}t(X)\exp(\ii\mu(\mu_{\mathbb A}(X)))=
\int_{\g g}\r(X)\exp(-\frac\ii2\OddId Q(v_{\0},\mu_{\mathbb A}(X)v_{\0}))
\end{equation}

\medskip

Since $ U_{+}^{\VC_{\1}}(M_{\0})$ contains a non empty subset, there exists a basis
$(G_{i}) $  of $\g g_{\0}$ such that $G_{i}\in U_{+}^{\VC_{\1}}(M_{\0})$. Let
$(g^{i})$ be its dual basis. Let $(H_{j})$ be a basis of $\g g_{1}$ and
$(h^{j})$ be its dual basis. We put $X_{\0}=\som_{i}G_{i}g^{i}$ and
$X_{\1}=\som_{j}H_{j}h^{j}$. We recall that a distribution can be integrated on
$\g g_{\1}$ and then on $\g g_{\0}$. We put:
\begin{equation}
\s(X_{\0})=\int_{\g g_{\1}}\r(X_{\0}+X_{\1})\exp(-\frac\ii2\OddId
Q(v_{\0},\mu_{\mathbb A}(X_{\1})v_{\0})).
\end{equation}

It is  a smooth compactly supported distribution on $\g g_{\0}$ with values in
$\C_{S(\Pi\VC^*)\mult_{M} \what\VC}(\what\VC_{\0})$. Now, we have:
\begin{equation}
\begin{split}
\int_{\g g}t(X)\exp\big(\ii\mu(\mu_{\mathbb A}(X))\big)&=
\int_{\g g}\r(X)\exp(-\frac\ii2\OddId Q(v_{\0},\mu_{\mathbb A}(X)v_{\0}))\\
&=\int_{\g g_{\0}}\s(X_{\0})\exp(-\frac\ii2\OddId Q(v_{\0},\mu_{\mathbb
A}(X_{\0})v_{\0}))\\
&=\what\s\Big(\frac12 \som_{i}\OddId Q(v_{\0},\mu_{\mathbb
A}(G_{i})v_{\0})g^{i}\Big)
\in\OC(\Pi\VC\mult_{M}\what\VC).
\end{split}
\end{equation}

Since $\s$ is smooth and compactly supported, $\what\s$ is rapidly decreasing on
$\g g_{\0}^*$. Since $G_{i}\in U_{+}^{\VC_{\1}}(M_{\0})$, for any
$\bm\in\what\VC_{\0}$, 
\begin{displaymath}
\OddId Q(v_{\0},\mu_{\mathbb
A}(G_{i})v_{\0})(\bm)=-\Big(\som_{r,s}Q(f_{r},\mu_{\mathbb
A}(G_{i})f_{s})(\OddId f_{s})^*(\OddId f_{r})^*\Big)(\bm)
\end{displaymath}
 is a negative definite form on
$(\VC_{x_{\bm}})_{\1}=((\Pi\VC\mult_{M}\what\VC)_{\bm})_{\0}$. Thus, thanks to
the Taylor formula (\ref{FonctLisse}), the above function is rapidly decreasing
along the fibres of $\Pi\VC\mult_{M}\what\VC$.

\medskip

Now, $-\frac12 \OddId Q(\bv,\bv)$ is constant along  the fibres of 
 $\Pi\VC\mult_{M}\what\VC\to\what\VC$, and $\exp(\ii\pi^*\mathbb A^{\OddId*}\OddId
Q^*(\bv))$ has (also by Taylor formula (\ref{FonctLisse})) almost polynomial
growth. Finaly,  hypothesis $(*)$ ensures that $\exp(\mu(\mathbb A^{\OddId 2}))$
has also almost polynomial growth. 

It follows that 
\begin{multline}
\int_{\g g}t(X)\exp(\w_{\VC}(X)) =\\
\exp\big(-\frac12\OddId Q(\bv,\bv)+\ii\pi^*\mathbb A^{\OddId*}\OddId
Q^*(\bv)+\mu(\mathbb A^{\OddId 2})\big)
\int_{\g g}t(X)\exp(\ii\mu(\mu_{\mathbb A}(X))).
\end{multline} is a rapidly decreasing function along the fibres of
$\Pi\VC\mult_{M}\what\VC\to \what\VC$.

\subsubsection{$\th$ is integrable along the fibres of $\VC$} This means that
for any smooth compaclty supported distribution $t$ on $\g g$, $<\th,t> $ is an
integrable form along the fibres of $\VC$. The problem being local on $M$ we
still assume that $\VC$ is trivial.

\medskip

Since $-\frac12 \OddId Q(\bv,\bv)$ is constant along  the fibres of
$\Pi\VC\mult_{M}\what\VC\to\what\VC$,
\begin{equation}
<\th,t>=\exp(-\frac12 \OddId Q(\bv,\bv))T\Big(\exp\big(\ii\pi^*\mathbb
A^{\OddId*}\OddId Q^*(\bv)+\mu(\mathbb A^{\OddId 2})\big)
\int_{\g g}t(X)\exp(\ii\mu(\mu_{\mathbb A}(X)))\Big).
\end{equation}
We recall that $v^{\OddId }$ is the generic point of $\C_{\Pi\VC\mult_{M}\what
M}({\what M}_{\0})$ (cf. formula (\ref{eq:v2})) and that $\DC_{\Pi\VC}$ is the
Liouville volume form along the fibres of $\Pi\VC\mult_{M}\what\VC$. We put:
\begin{equation}
\begin{split}
\psi(v^{\OddId})&=\DC_{\Pi\VC}(v^{\OddId })\exp\big(\mu(\mathbb A^{\OddId
2})(v^{\OddId})\big)
\int_{\g g}t(X)\exp(\ii\mu(\mu_{\mathbb A}(X))(v^{\OddId }))\\
&=\DC_{\Pi\VC}(v^{\OddId})\int_{\g g}t(X)\exp(\mu(F_{\g g}(X))(v^{\OddId }))
\in\C_{\vol(\Pi\VC)\mult_{M}\what \VC}(\what\VC_{\0}).
\end{split}
\end{equation}

The preceding section shows that $\psi$ is a rapidly decreasing volume form on the fibres of $\Pi\VC\mult_{M}\what\VC$
\medskip

Before going further, we have to consider the objects involved under another
point of view.

First, locally, we consider
$\bv\in\C_{\Pi\VC}(\UC)\tens_{\OC_{M}(\UC)}\C_{\VC^*}(\UC)\simeq
\Hom_{\OC_{M}(\UC)}(\C_{\VC}(\UC),\C_{\Pi\VC}(\UC))$. Thus $\bv$ is identified
to odd canonical isomorphism of sheafs of $\OC_{M}$-modules 
$\OddId:\C_{\VC}\to\C_{\Pi\VC}$.

 Now, since  $\mathbb A=d+\w$ with $\w\in\what\W_{M}(M_{\0},\g{spo}(\VC))_{\1}$,
$\pi^*\mathbb A^{\OddId}\bv$ realizes an even isomorphism of sheafs of 
 $\OC_{\VC\mult_{M}\what{M}}$-modules between $\C_{\what\VC}$ and
$\C_{\Pi\VC\mult_{M}(\VC\mult_{M}\what M)}$.

Finaly $\pi^*\mathbb A^{\OddId*}\OddId Q^*(\bv)=\pi^*\OddId Q^*(\mathbb A^{\OddId}\bv)$
realizes an even linear isomorphism of sheafs of 
 $\OC_{\VC\mult_{M}\what{M}}$-modules between $\C_{\what\VC}$ and
$\C_{\Pi\VC^*\mult_{M}(\VC\mult_{M}\what M)}$.

It follows that:

\begin{equation}
\begin{split} T\Big(\exp\big(\ii\pi^*\mathbb A^{\OddId*}\OddId
Q^*(\bv)+\mu(\mathbb A^{\OddId2})\big)
\int_{\g g}t(X)\exp(\ii\mu(\mu_{\mathbb A}(X)))\Big)=\what\psi (\pi^*\mathbb
A^{\OddId*}\OddId Q^*(\bv))
\end{split}
\end{equation}
is rapidly decreasing along the fibres of $\what \VC\to \VC\mult_{M}\what M$.

\medskip

Since $\exp\big(-\frac12 \OddId Q(\bv,\bv)\big)$ is constant along these fibres
and rapidly decreasing along the fibres of $\VC\mult_{M}\what M\to\what M$, it
follows that 
\begin{displaymath} <\th,t>=\exp\big(-\frac12 \OddId Q(\bv,\bv)\big)\what\psi
(\pi^*\mathbb A^*\OddId Q^*(\bv))
\end{displaymath}
 is a pseudodifferential  form on $\VC$ that is integrable along the fibres.

\subsubsection{Evaluation of $\pi_*\th$} Now, we have  to check that $\pi_*\th$
is a constant generalized function on $\g g$.

We stress that in $\Pi\VC\mult_{M}\what\VC$ there are two copies of $\Pi\VC$, we
denote by $\what\Pi\VC$ the copy of $\Pi\VC$ in $\what\VC$:
$\what\VC=\VC\mult_{M}\what\Pi\VC\mult_{M}\what M$.

We denote by $v'$ the generic point of $\C_{\VC\mult_{M}\what M}(\what M_{\0})$
and by $v^{\what\OddId}$ the generic point of $\what\Pi\VC\mult_{ M}\what M$.

We recall that
$\bv\in\C_{\Pi \VC\mult_{M}\VC}(\VC_{\0})\subset\what\W(\VC,\Pi\VC\mult_{M}\VC)$.
It follows from formula (\ref{eq:ConnTens}) with $\mathbb A$ in place of $\mathbb{A}_{\g g}$
that:
\begin{equation}
\pi^*\mathbb A^{\OddId*}\OddId Q^*(\bv)=\OddId Q^*(\pi^*\mathbb
A^{\OddId}\bv)=\OddId Q^*((1\tens d)\bv)+\OddId Q^*((\mathbb A^{\OddId}\tens
1)\bv).
\end{equation}
We put $\bd=1\tens d$. Since $\OddId Q^*((\mathbb A^{\OddId}\tens 1)\bv)$ is
constant along the fibres of $\what\VC\to\VC\mult_{M}\what M$ and
$d_{(d\xi,dx)}$ is invariant by translation, we obtain:

\begin{equation}
\begin{split} <\pi_{*}\th,t>&=\pi_{*}\big(<\th,t>\big)\in\what\W_{M}(M_{\0})\\
&=\frac{1}{(-1)^{\frac{k+l}{2}}}
\int_{\what\VC/\what M}d_{(x,d\xi,dx,\xi)}(v',v^{\what\OddId})\\
&\hskip 4cm
\exp\big(-\frac12 \OddId Q(\bv,\bv)\big)(v')
\what\psi (\pi^*\mathbb A^*(\OddId Q^*(\bv)))(v',v^{\what\OddId})
\\
&=\frac{1}{(2\pi)^{\frac{k+l}{2}}}\int_{\VC\mult_{M}\what M/\what M}d_{(x,\xi)}(v')\exp\big(-\frac12 \OddId
Q(\bv,\bv)\big)(v')\\
&\hskip 4cm
\int_{\what\Pi\VC\mult_{M}\what M/\what M}
d_{(d\xi^{j},dx^{i})}(v^{\what\OddId})\what\psi (\OddId
Q^*(\bd\bv))(v^{\what\OddId}).
\end{split}
\end{equation}

Now, as above for $\pi^*\mathbb A^{\OddId*}\OddId Q^*(\bv)$,
$\bd\bv:\what\VC=\what\Pi \VC\mult_{M}\VC\mult_{M}\what
M\to\Pi\VC\mult_{M}\VC\mult_{M}\what M$ is an even isomorphism of supervector
bundles. It induces an isomorphism $\ber(\bd\bv):\ber(\what\Pi\VC)\to\ber(\Pi\VC)$.

We recall that $(e_{i},f_{j})$ is a standard basis of $\C_{\VC}(\UC)$ and
$(x^{i},\xi^{j})$ is its dual basis. Then $(\OddId f_{j},\OddId e_{i})$ is a
standard basis of $\C_{\Pi\VC}(\UC)$. Its dual basis is $(-\OddId \xi^{j},\OddId
x^{i})$. Thus we have:
\begin{equation}
\ber(\bd\bv)\big(d_{(d\xi^{j},dx^{i})}\big)=(-1)^{l}(d_{(\OddId \xi^{j},\OddId x^{i})}).
\end{equation}
Thus:
\begin{equation}
\int_{\what\Pi\VC\mult_{M}\what M/\what M}
d_{(d\xi^{j},dx^{i})}(v^{\what\OddId})\what\psi (\OddId
Q^*(\bd\bv))(v^{\what\OddId})
=
(-1)^l \int_{\Pi\VC\mult_{M}\what M/\what M}
d_{(\OddId\xi^{j},\OddId x^{i})}(v^{\OddId})\what\psi (\OddId
Q^*(v^{\OddId})).
\end{equation}

\medskip

 On the other hand, $d_{(x^{i},\xi^{j})}\in\ber(\VC)$. We recall that we
canonically identified $\ber(\VC)$ and $\ber(\Pi\VC^*)$.   We assume that
$(-\OddId f_{j},\OddId e_{i})$ is an oriented symplectic basis of  $\Pi\VC$.
Under this identification:
\begin{equation}
d_{(x^{i},\xi^{j})}=d_{(-\OddId f_{j},\OddId e_{i})}
=\frac{(2\pi)^l}{(-1)^{\frac{k(k-1)}2}\ii^{k}}(d_{(\OddId\xi^{j},\OddId
x^{i})})^* ;
\end{equation}
and
\begin{equation}
\DC_{\Pi\VC}=\frac{1}{(2\pi){\frac{l}{2}}}d_{(-\OddId f_{j},\OddId e_{i})}
=\frac{1}{(2\pi)^{\frac{l}{2}}}d_{(x^{i},\xi^{j})}
\end{equation}

It follows from Fourier inversion formula (since $l$ is even $(-1)^{l}=1$):

\begin{equation}
\begin{split}
 d_{(x,\xi)}\int_{\what\Pi\VC\mult_{M}\what M/\what M}
d_{(d\xi^{j},dx^{i})}(v^{\what\OddId})\what\psi (\OddId
Q^*(\bd\bv))(v^{\what\OddId})
&=
\frac{(2\pi)^{l}}{(-1)^{\frac{k(k-1)}2}\ii^{k}}\DC_{\Pi\VC}
\int_{\g g}t(X)\\
&=\frac{(2\pi)^{\frac{l}{2}}}{(-1)^{\frac{k(k-1)}2}\ii^{k}}d_{(x,\xi)}
\int_{\g g}t(X).
\end{split}
\end{equation}

Since $k$ is even, $(-1)^{\frac{k(k-1)}{2}}\ii^k=1$. Moreover, $\OddId Q(\bv,\bv)\big)(v')=
 Q(v',v')$. Thus:

\begin{equation}
\begin{split} <\pi_{*}\th,t>&=\frac{1}{(2\pi)^{\frac{k+l}{2}}}
\frac{(2\pi)^{\frac{l}{2}}}{(-1)^{\frac{k(k-1)}2}\ii^{k}}
\int_{\VC\mult_{M}\what M/\what M}d_{x,\xi}(v')
\exp\big(-\frac12  Q(v',v')\int_{\g g}t(X)\\
&=\frac{1}{(2\pi)^{\frac{k+l}{2}}}
\frac{(2\pi)^{\frac{l}{2}}}{(-1)^{\frac{k(k-1)}2}\ii^{k}}(2\pi)^{\frac k2}\int_{\g
g}t(X)\\
&=\int_{\g g}t(X)
\end{split}
\end{equation}

\medskip

 Consequently $\th$
  is a Thom form on $\g g$. Moreover its restriction to $\pi^{-1}(W)$ is 
$\CC^\oo$ on $U^{\VC_{\1}}(W).$

\end{proof}

{\bf Remark:} If $V=V_1$, the hypotheses imply that the moment map $\mu_{\mathbb
A}$ is proper. Let $\phi$ be smooth compactly supported function
 on $\g g_0$. The fact that $\int_{\g g_0}dX
\exp\ii\mu\big(\mu_{\mathbb A}(X)\big)\phi(X)$ defines  an integrable function 
on the vector space $V_1$ is a particular case of Lemma 12 in \cite{Ver97}.

\vskip 2mm

{\it Example:} Let us consider the case where $M$ is a point, $\VC=\mathbb
R^{(0,2)}$ ($\C_\VC(M)=\mathbb R^{(0,2)}$) and $\g g=\mathbb R X$. The action of $X$ is
defined by the field of vectors:
\begin{equation}
X_{\VC}=\xi{\frac{\de}{\de \eta}}-\eta{\frac{\de}{\de\xi}}.
\end{equation}
In the canonical basis $(f_1,f_2)$ of $\mathbb R^{(0,2)},$ 
$\LC^\VC(X)$ is represented by
\begin{equation}
X=
\begin{pmatrix}0&1\\
-1&0
\end{pmatrix}.
\end{equation}
The Euclidean structure $Q$ is given by $Q(f_1,f_2)=1$ We denote by $z$ the
element of $\g g^*$ such that $z(X)=1$. We denote by $(\xi,\eta)$  the dual
basis of $(f_{1},f_{2})$. Thus $(\OddId\xi,\OddId\eta)$ is a basis of
$\Pi(\mathbb R^{(0,2)})^*=(\mathbb R^2)^*$. We have:
\begin{equation}
\w_\VC=\xi\eta-\ii(d\xi \OddId\eta -d\eta
\OddId\xi)-{\frac{\ii}2}z\big((\OddId\xi)^2+(\OddId\eta )^2\big).
\end{equation}

On the open set $U=(\mathbb R\setminus\{0\}) X$, we have
\begin{equation}
\begin{split}
\th&=\int_{\Pi\mathbb R^{(0,2)}}|d(\OddId\xi) d(\OddId\eta)| \,\exp(\w_\VC)\\ 
&
\begin{aligned}
=(1+\xi\eta)\int_{\mathbb R^2}|d(\OddId\xi) d(\OddId\eta)| 
\exp\Big(-{\frac{\ii}2}z\Big((\OddId\xi -&{\frac{d\eta} z})^2 +(\OddId\eta
+{\frac{d\xi} z})^2\Big)\\
&+{\frac{\ii}2}\big({\frac{d\xi^2} z}+{\frac {d\eta^2} z}\big)\Big)
\end{aligned}
\\  
&={\frac{2\ii\pi}{ z}}\exp\Big(\xi\eta-{\frac{\ii}{
2z}}\big(d\xi^2+d\eta^2\big)\Big).
\end{split}
\end{equation}

It appears that $\th$ is rapidly decreasing on $\what{\mathbb R^{(0,2)}} $ but
only as a generalized function on $\g g.$ More precisely, for any compactly
supported function $\phi\in\CC^\oo(\mathbb R)$ which vanishes at $0$, we have:
\begin{equation}
\int_{\mathbb R}|dz|\,\th(z)\phi(z) =-2\ii\pi
\exp(\xi\eta)\what{\Big({\frac{du}{u}}{\phi\big({\frac{1}{u}}\big) }\Big)}\big({-\frac1
2}(d\xi^2+d\eta^2)\big),
\end{equation}
where $\Big(\what{{\frac{du}{u} }\phi({\frac1 u})}\Big)$ denotes the Fourier
transform of the smooth distribution ${\frac{du}{u}}\phi({\frac{1}{u}})$ on $\mathbb
R\setminus
\{0\}$. Since $\phi$ is compactly supported, this integral is rapidly decreasing
in 
$d\xi$ and $d\eta$.

\subsection{A relation between cohomology classes}\label{RelCoho}

Let $G=(G_{\0},\g g)$ be a supergroup. Let $\pi:\VC\mapsto M$ be an oriented
$G$-equivariant supervector bundle. We assume that $\VC$ and $M$ are globally
oriented.

\begin{lemme} Let $\UC$ be a $G$-invariant open subset of $\g g_{\0}$. Let $\th$
be a Thom form on $\VC$. Let $\a\in\what\W^{\oo}_{G,\int}(\UC,M)$, then
\begin{equation}
\th(\pi^*\a)\in\what\W_{G,\int}^{-\oo}(\UC,\VC).
\end{equation}
\end{lemme} 

{\bf Remark:} This product is well defined because $\a$ has smooth coefficients.

\begin{proof}[Proof] We have to show that this product defines a generalized
function in the integrable forms. It is enough to check locally that we can
integrate along the fibres of $\VC\rightarrow M$, and then to check that the
resulting form is integrable on $M$. But since $\pi^*(\a)$ is constant along the
fibres, the integrability along the fibres comes from that of $\th$. Thus we
have $\pi_*( \th\ \pi^*(\a))=\pi_*(\th)\a=\a$. Hence, the resulting form is
$\a$, which is integrable by hypothesis.
\end{proof}

\begin{theo} Let $\UC$ be a $G_{\0}$-invariant open subset of $\g g_{\0}$. Let
$\th$ be a Thom form on $\VC$. Let $\a$ be an equivariantly closed form with
smooth coefficients ($\a\in\what\W_{G,\int}^\oo(\UC,\VC)$). We have the
following equality between cohomology classes in  $\what
H_{G,\int}^{-\oo}(\UC,\VC)$:
\begin{equation}
\a\equiv \th\pi^*(\pi_*\a).
\end{equation}
\end{theo}

\begin{proof}[Proof] The proof in the even situation (cf.
\cite{KV93}) can be repeated here. We consider the supervector bundle $\VC\mult_{M}\VC$
over $M$. For $t\in\mathbb R$, we denote  by $\s_t$ the linear transformation in
the fibres of $\VC\mult_{M}\VC$ which is defined for any section $x$ and $y$ of
$\VC$ by:

\begin{equation}
\s_t(x,y)=\Bigl((\cos t)x+(\sin t)y,-(\sin t)x+(\cos t)y\Bigr).
\end{equation}
We have $\s_0=Id$, and 
$\s_{\frac\pi 2}=\s$ where, for any sections $x$ and $y$ of $\VC$,
$\s(x,y)=(y,-x)$. We put:
\begin{equation}
 S={\frac{d} {dt}}\s_t\in\C_{T(\VC\mult_{M}\VC)}((\VC\mult_{M}\VC)_{\0}).
\end{equation} 
Since for all $t$, $\s_t$ commutes with the action of $G$, $\LC(S)$ and
$\iota(S)$ leave $\what\W^{-\oo}_{G,\int}(\g g,\VC\mult_{M}\VC)$ invariant. Thus we have the following relation (between derivations of $\what\W^{-\oo}_{G,\int}(\g g,\VC\mult_M\VC)$):
\begin{equation}
\LC(S)=d_{\g g}\iota(S)+\iota(S)d_{\g g},
\end{equation}
($S$ is an even field of vectors and thus  for $X\in\g g$
$\iota(X)\,\iota(S)+\iota(S)\,\iota(X)=[\iota(X),\iota(S)]=0$).

We put for $\nu\in\what\W_{G,\int}^{-\oo}(\g g,\VC\mult_M\VC)$:
\begin{equation}
 H\nu(X)=\int_0^{\frac{\pi}2}\,\s_t^*(\iota(S)\nu(X))\, dt .
\end{equation}
It defines an application $H:\what\W^{-\oo}_{G,\int}(\g g,\mathcal
V\mult_{M}\mathcal V)\rightarrow\what\W^{-\oo}_{G,\int}(\g g,\mathcal
V\mult_{M}\mathcal V)$.

For any $\nu\in\what\W^{-\oo}_{G,\int}(\g g,\VC\mult_M\VC)$ we have:
\begin{equation}
\s^*\nu-\nu=(d_{\g g}H-Hd_{\g g})\nu.
\end{equation}

We denote by $p_i:\VC\mult_M\VC\rightarrow\mathcal V, \ i=1$ (resp. 2), the
projection on the first (resp. the second) factor. 

\medskip

Let $\a$ be an integrable
equivariantly closed form on $\VC$. Then (for similar reasons as those we used
in the preceding lemma) $(p_1^*\th )(p_2^*\a)$ is a well defined in 
$\what\W^{-\oo}_{G,\int}(\g g,\mathcal V\mult_{M}\mathcal V)$. 

It is a closed form because $\a$ and $\th$ are closed. The previous relation
shows that $(p_1^*\th )(p_2^*\a)$ is in the same cohomology class as
$\s^*((p^*_1\th) (p^*_2\a))$. 

For any form $\be$ on $\VC$, we denote by $\overline \be$ the form $\t^*\be$
where $\t:\VC\rightarrow\VC,$ $v\mapsto -v.$ So we have
$\s^*((p^*_1\th)(p^*_2\a))=(p^*_2\th)(\ p^*_1\overline\a)$.
 We integrate  these two cohomologically equivalent forms along the fibres of
$p_1$. We obtain ($\pi_{*}$ and $\th$ have same parity $p(\pi_{*})=p(\th)=\frac{n+m}{2}(mod 2)$):
\begin{align}
 (p_1)_*(p_1^*\th\ p_2^*\a)&=\pi^*(\pi_*\th)\a,\\
\begin{split} (p_1)_*(p_2^*\th
p_1^*\overline\a)&=(-1)^{(n+m)p(\th)}\th\pi^*(\pi_*\overline\a),\\ 
&=(-1)^{(n+m)(p(\th)+1)}\th\pi^*(\pi_*\a),\\ 
&=\th\pi^*(\pi_*\a).
\end{split}
\end{align} Since $\pi_*\th=1,$ we obtain the requested equality in cohomology:
\begin{equation}
\th\pi^*(\pi_*\a)\equiv\pi^*(\pi_*\th)\a=\a.
\end{equation}
\end{proof}

\section{Equivariant Euler form}\label{Eul.sec}

\begin{defi} Let $G=(G_{\0},\g g)$ be a supergroup. Let $M$ be a
$G$-supermanifold. Let $\VC\rightarrow M$ be an 
equivariant supervector bundle. We assume that $\VC$ and $M$ are globally
oriented. Let $j$ be the injection of $M$ into $\VC$ by means of the zero
section. We assume that there is an equivariant Thom form $\th$ on $\VC$. We
put:
\begin{equation}
\EC_{\g g}=j^*\th\in\what\W_{G}(M).
\end{equation}
We say that $\EC_{\g g}$ is an equivariant Euler form on $\VC$.

\end{defi}

Now, we assume  that there is a $G$-invariant Euclidean structure denoted by $Q$
on the fibres and an equivariant superconnection denoted by $\mathbb A_{\g g}$
on $\VC$. We denote by $F_{\g g}=\mathbb A_{\g g}^2+\ii\LC^{\mathcal V}$ its
equivariant curvature. It is an $\g{osp}(\VC,Q)$-valued equivariant form.

We recall that we identified $\g{osp}(\VC,Q)$ and $\g{spo}(\Pi\VC,\OddId Q)$.
This identifies $F_{\g g}$ and $F_{\g g}^{\OddId}$. Thus, we can consider 
$\mu(F_{\g g})$ instead of $\mu(F_{\g g}^{\OddId})$.

We assume that the hypothesis in Theorem \ref{Thom} are satisfied (in particular
$(*)$ and $(**)$).

Let $\DC_{\Pi\VC}$ be
 Liouville volume form on the oriented symplectic bundle $(\Pi\VC\mult_{M}\what M,\OddId Q)$.

We put (compare with \cite{Lav03}):

\begin{equation}
\Spf(-\ii F_{\g g}(X))=\ii^{\frac{k-l}2} \int_{\Pi\VC\mult_{M}\what M/\what M}\DC_{\Pi\VC}(v^{\OddId})\exp\big(\mu(F_{\g
g}(X)(v^{\OddId}))\big).
\end{equation}
The hypotheses $(*)$ and $(**)$ ensure that this definition makes sense as a
generalized function on $\g g$ with values in the pseudodifferential forms on
$M.$ It means that for any smooth compactly supported  distribution $t$ on $\g
g$, the integral 
\begin{equation}
\int_{\g g}t(X)\,\exp(\mu(F_{\g g}(X)))\,\in\OC(\Pi\VC\mult_{M}\what M)
\end{equation}
is rapidly decreasing along the fibres of $\Pi\VC\mult_{M}\what
M\rightarrow\what M$ (cf. subsection \ref{sec:ThomGene} in demonstration of
Theorem
\ref{Thom}). Then by definition:
\begin{equation}
<\Spf(-\ii F_{\g g}(X)),t(X)>=\ii^{\frac{k-l}2} \int_{\Pi\VC\mult_{M}\what M/\what M}\DC_{\Pi\VC}(v^{\OddId})
\int_{\g g}t(X)\,\exp(\mu(F_{\g g}(X)(v^{\OddId}))).
\end{equation}

\begin{prop} We have the following equality between cohomology classes:

\begin{equation}
\EC_{\g g}\equiv \frac 1 {\ii^{\frac{k-l}2}} \Spf({-\ii F_{\g g}}),
\end{equation}
where $(k,l)$ is the dimension of the fibres of $\VC$.
\end{prop} {\bf Remark:} In \cite{BGV92} the equivariant Euler form is defined
by the preceding formula (in the purely even case).
\begin{proof}[Proof] It is enough to show that this equality is true for the
particular equivariant Euler form obtained with the Thom form constructed in
Theorem \ref{Thom}.

 When $k$ is odd, $Spf(-\ii F_{\g g})$ is zero. 
 
 When $k$ is even, we have (the existence of an Euclidean structure ensures that $l$ is even):
\begin{equation}
\begin{split} j^*\th&=\int_{\Pi\VC\mult_{\what M}\what M/\what M}\DC_{\Pi\VC}(v^{\OddId})\exp(\mu(F_{\g g})(v^{\OddId}))\\ 
&= \frac 1 {\ii^{\frac{k-l}2}} \Spf(-\ii F_{\g g}).
\end{split}
\end{equation}
\end{proof}

\section{First localization formulas}

\subsection{The linear situation} Let $G=(G_{\0},\g g)$ be  a supergroup and $V$
be a supervector space of dimension $(m,n)$. We assume that we have a
representation of $G$ in $V$. We  denote by $\r$ the associated representation
of $\g g$. We assume  that $V$ is  $G$-invariantly globally oriented and has a
$G$-invariant Euclidean structure $Q$. In particular $m$ is even.

Let $U$ be the open subset of $\g g_{\0}$ definede by:
 \begin{equation}
U=\big\{X\in\g g_\0\,/\,\r(X) \text{ is invertible}\big\}.
\end{equation}

Let $U_{+}$ be the open subset of $\g g_{\0}$ such that:
\begin{equation}
U_{+}=\big\{X\in \g g_\0\,\big/\,v\mapsto Q(v,\r(X)v)\text{ is a scalar product
on $V_{\1}$.}\big\}.
\end{equation}

We have $U_{+}\subset U$. In order to satisfy condition $(**)$, we assume that
$U_{+}$ is not empty. In particular, it implies that $m$ is even.

\begin{prop}\label{FormLin} Let $\a\in\what\W^{\oo}_{G,\int}(U,V)$ be an
equivariantly closed form on $V$. We denote by $j$ the canonical injection of
$\{0\}$ in $V$. Then we have in $\CC^{\oo}(U)$:
\begin{equation}
\int_V \a(X)=\ii^{\frac{m-n}2}(2\pi)^{{\frac{n+m}2}} \frac{(j^*\a)(X)}{
Spf(\r(X))}.
\end{equation}

\end{prop}

\begin{proof}[Proof] Following section \ref{RelCoho}, we know that the
cohomology class of $\a$ in  $\what H^{\oo}_{G,\int}(U,V)$ is equal to that of
$\th(\pi_{*}\a)=\frac 1{(2\pi)^{\frac {m+n}2}}\th(\int_V\a)$, where $\th$ is an
equivariant Thom form on $V$ ($\th$ is smooth on $U$). 

We use the notation of the preceding section, just replacing $\VC$ by $V$ and
paying attention to the fact that, in this case, we can take $\mathbb A=0$ ($M$
is a point). The equivariant moment $\mu_{\mathbb A}$ is the representation $\r$
of $\g g$ in $V$ ($\r:\g g\rightarrow \g{osp}(V)$) and $F_{\g g}=\ii\r$.

We just have to evaluate $j^*\th$. 
In this case condition $(*)$ is trivial and condition $(**)$ follows from  the assumption
that $U_{+}$ is not empty.

We obtain:
\begin{equation}
 j^*\th(X)=\frac1{\ii^{\frac{m-n}2}}Spf(\r(X)).
\end{equation} 
 
 Since on $U$, the function $X\mapsto Spf(\r(X))$ is smooth and invertible, the formula follows.
\end{proof}

\subsection{The fibered situation} We can generalize the preceding formula to
the case of a $G$-equivariant oriented Euclidean supervector bundle. 

We put:
\begin{equation}
\label{eq:Utot} U^{\VC}(\UC)=\text{Interior of }\Big\{X\in\g g_0\,\Big/
\forall \bm\in (\what M(\UC))_{\0},
\mu_{\mathbb A}(X)(\bm)\text{ is invertible}\Big\}.
\end{equation}

We have $U^{\VC}(\UC)\subset U^{\VC_{\1}}(\UC)$.

\begin{prop}\label{FormFib} Let $G=(G_{\0},\g g)$ be a supergroup and $M$ be a
$G$-supermanifold. Let $\VC\rightarrow M$ be a $G$-equivariant oriented
Euclidean supervector bundle. 

We  assume that there is a superconnection $\mathbb A$ satisfying 
condition $(*)$.  We denote by $\mu_{\mathbb A}$ the equivariant moment of
$\mathbb A$.

We assume that  $U_{+}^{\VC_{\1}}(M)_0\subset \g g_0$ (cf formula (\ref{U0m})) 
contains a non empty open subset, which implies 
condition $(**)$ (cf. Theorem
\ref{Thom} for conditions $(*)$ and $(**)$). 

Let $\a\in\what\W^{\oo}_{G,\int}(U^\VC(M),\VC)$ be an equivariant integrable
form with smooth coefficients ($U^\VC(M)$ defined by formula (\ref{eq:Utot})). 

We denote by $j$ the injection of $M$ into $\VC$ by means of  the zero section
and by $\EC_{\g g}$ an equivariant Euler form on $M$ associated to $\VC$. We
have the following equality of functions in $\CC^\oo(U^\VC(M))$:
\begin{equation}
\int_{\VC}\a(X)=(2\pi)^{\frac{m+n}2}\int_M\frac{j^*\a(X)}{\EC_{\g g}(X)}
\end{equation}
(On $U^\VC(M),$ $\EC_{\g g}$ is smooth and invertible.)
\end{prop}
\begin{proof}[Proof] Exactly the same as in the case  where $M$ is just a point.
\end{proof}

\section{Preliminaries to the localization formula}\label{preLoc}

\subsection{Introduction} Let $G=(G_{\0},\g g)$ be a connected supergroup. Let
$M$ be  a globally oriented $G$-supermanifold. We also assume that $M$ has a
$G$-invariant weak Euclidean structure denoted by $Q$. Let $\a$ be an
equivariantly closed form with smooth coefficients  defined on  $\g g_{\0}.$ Let
$X\in\g g_\0$ such that $\overline{\exp(\mathbb R X)}\subset G_\0$ is compact.
Let $M(X)$ be the supermanifold of zeroes of $X$ in $M$ (it will be defined in the
next subsection).  We express $\int_M\a$ on a neighborhood of $X$ in $\g
g(X)_{\0}$ in terms of an integral on $M(X)$. This method comes from
\cite{BV83a},\cite {BV83b} and \cite{Bis86}.

 To achieve this we have to assume that the normal bundle $T_NM(X)$ of $M(X)$ in
$M$ (cf. section \ref{sec:local} for a definition) has a $G(X)$-invariant
Euclidean structure, and has a $G(X)$-invariant superconnection, $\mathbb A$,
which preserves  the Euclidean structure of $T_NM(X)$.

Moreover, we assume:
\begin{enumerate}
\item [$(*)$] ${\mu(\mathbb A^{\OddId 2})}_{\mathbb {R}}\le 0$ on $\big(\Pi
T_NM(X)\big)_0$;
\item[(**)] there is a covering of $M_{\0}$ by open subsets $W$ such that
$U_{+}^{T_NM(X)_{\1}}(W)$ contains a non empty open subset (cf. formula (\ref{U0m})
for $U_{+}^{T_NM(X)_{\1}}(W)$).
\end{enumerate}

\vskip 2mm

First, we define $M(X)$, the manifold of zeroes of $X$ in $M$.

Then, we construct:
\begin{enumerate}
\item an open neighborhood $U$ of $X$ in $\g g_{\0}$,
\item  an equivariant form $\be\in\what\W_G^{\oo}(U,M)$,
\item an open neighborhood $V$ of $M(X)_{\0}$ in $M$,
\end{enumerate}

 such that $d_{\g g}\be$ is invertible as equivariant form on $M$ defined on the
complement of $\overline V$ in $M_\0$ (where $\overline V$ is the closure of
$V$).

Then, we prove the localization formula in the case of isolated zeroes and in
the general case.

\subsection{Supermanifold of zeroes}

Let $X\in\g g_\0$ be such that $\overline{\exp(\mathbb R X)}\subset G_\0$ is
compact. Let $\IC$ be the sheaf of ideals of $\OC_M$ such that for any open
subset $\UC$ of $M$, $\IC(\UC)$ is the ideal generated by
\begin{displaymath}
\{ X_Mf,\  f\in\OC(\UC)\} .
\end{displaymath}

We recall that $(X_{M})_{\mathbb R}$ is the canonical projection of $X_{M}$ in
$\C_{TM}(M_{\0})\tens_{\OC_{M}(M_{\0})}\CC^{\oo}(M_{\0})$. Since
$\overline{\exp(\mathbb R X)}$ is compact, the set of
 zeroes of  $(X_M)_{\mathbb R}$ in $M_\0$ is a manifold denoted by $M(X)_\0$ (cf. \cite{BGV92}). 

\begin{prop} The sheaf of ideals $\IC$ defines a subsupermanifold of $M$ denoted
by $M(X)$ which  underlying manifold is $M(X)_\0$ and such that for any open
subset $\UC\subset M_{\0}$:
\begin{equation}
\OC_{M(X)}(\UC\cap M(X)_{\0})=\OC_{M}(\UC)/\IC(\UC).
\end{equation}
\end{prop}

Then we put:

\begin{defi} We call supermanifold of  zeroes of $X\in\g g_\0$ in $M$ the subsupermanifold
$M(X)$ of $M$ defined in the preceding proposition.
\end{defi}

{\em Example:} Let $M$ be a supervector space with a representation
$(\r_{\0},\r)$ of $G=(G_{\0},\g g)$. Let $X\in\g g_{\0}$. Then
$M(X)=\ker(\r(X))$.

\medskip

Let $G_{\0}(X)$ (resp. $\g g(X)$) be the centralizer of $X$ in $G_{\0}$ (resp.
$\g g$). We put $G(X)=(G_{\0}(X),\g g(X))$. It is a subsupergroup of $G$ called
the centralizer of $X$ in $G$.

The supermanifold $M(X)$ is stable under the action of the centralizer $G(X)$ of
$X$ in $G.$

\subsection{An equivariant form with invertible  equivariant differential}

\subsubsection{Construction of the form} Let $X\in\g g_\0$ be a central element in
$\g g$ ($\g g(X)=\g g$). Let $M(X)$ be the set of  zeroes of $X$.

Let $\a$ be an equivariantly closed form on $M$ with smooth coefficients 
defined on $\g g_\0$. We denote by $Supp(\a)$ its support. It is the smallest
closed subset $F$ of $M_\0$ such that $\a$ vanishes on $M_\0\setminus F$. We
assume that $Supp(\a)$ is compact.

We will prove that if  $Supp(\a)\cap M(X)_\0=\O,$ there is a $G_{\0}$-invariant
open neighborhood $U'$ of $X$ such that the restriction of $\a$ to $U'$ is
$d_{\g g}$-exact.

 For any $Y$ in $\g g$, we put:  
\begin{equation}\label{eq:beta}
\be(Y)=\OddId Q^*(\OddId Y_M)\in \C_{\Pi
T^*M}(M_{\0})\subset\what\W_{M}(M_{\0}).
\end{equation}
This defines $\be\in\g g^*\tens\what\W(M)\subset\OC_{\g g}(\g
g_{\0},\what\W(M))$. Since $Q$ is $G$-invariant, this form is equivariant.

We have $d_{\g g}\be(Y)=-\ii Q(Y_M,Y_M)+d\be(Y)$. To ensure that $d_{\g g}\be(Y)$
is invertible it is sufficient to show that $Q(Y_M,Y_M)_{\mathbb R}$ is not
zero. First we remark that the function
\begin{displaymath} (m,Y)\mapsto Q_m(Y_M(m),Y_M(m))
\end{displaymath}
 is continuous and positive on $M_0\times \g g_0$. 

\subsubsection{Construction of $U_a(\a,X,V)$} We use the following lemma:

\begin{lemme} Let $W$ and $V$ be two topological spaces where $V$ is compact,
and $\phi$ be a continuous function on $W\times V$ into $\mathbb R$. Then the
set
\begin{equation}
\{m\in W,\forall v\in V\ \phi(m,v)>0\}
\end{equation}
is open in $W.$
\end{lemme}

Let $V$ be a closed neighborhood of $M(X)_{\0}$. Let $a>0$. Then by the
preceding lemma:
\begin{equation}
\label{Ua} U_a(\a,X,V)=\big\{Y\in\g g_0\,/
\forall m\in Supp(\a)\cap\overline{M_0\setminus V},
\,Q_m(Y_M(m),Y_M(m))>a\big\}
\end{equation}
is open (we recall that $\overline{M_{\0}\setminus V}$ designs the closure of
$M_{\0}\setminus V$ in $M_{\0}$).

Moreover  $Q(X_M,X_M)_{\mathbb R}>0$ on $M_\0\setminus M(X)_\0$. As $Supp(\a)$
is compact, $Q(X_M,X_M)_{\mathbb R}$ has a strictly positive lower bound $A$ on
$Supp(\a)\cap
\overline{M_\0\setminus V}$.  Then $U_a(\a,X,V)$  is not empty if $a\le A$.

\subsubsection{Exactitude of $\a$}

\begin{lemme} Let $X\in\g g_0$, central in $\g g$.

Let $\a\in\what\W_G^{\oo}(\g g_{\0},M)$ such that $d_{\g g}\a=0$, $Supp(\a)$ is
compact and  $Supp(\a)\cap M(X)_0=\O$. 

Then, there exists an open neighborhood $U$ of $X$ in $\g g_{\0}$ such that
 $\a$ is $d_{\g g}$-exact in $\what\W_G^{\oo}(U,M).$ 

Moreover, if $\a$ is integrable, $\a$ is $d_{\g g}$-exact in
$\what\W_{G,\int}^{\oo}(U,M).$
\end{lemme}

\begin{proof}[Proof] Since $Supp(\a)\cap M(X)_\0=\O$ and $X$ is central in $\g
g$ there exists a $G$-invariant closed neighborhood $V$ of $M(X)_\0$ such that
$V\cap Supp(\a)=\O.$ Then, it is enough to take $U=U_a (\a,X,V)$ for a small
enough $a$  and to check that $\a=d_{\g g}\Bigl({\frac{\be\a}{d_{\g
g}\be}}\Bigr)$. 

If $\a$ is integrable, ${\frac{\be\a}{d_{\g g}\be}}$ is integrable because
$\frac{\be }{d_{\g g}\be}$ is a smooth bounded function on
$U_a(\a,X,V)\times\what{M}(M_{\0}\setminus V).$
\end{proof}

{\bf Remark:} Since $\be$ is smooth, the lemma is still true for an equivariant
form $\a$ with generalized coefficients.

\section{Localization formula}

\subsection{The localization}\label{sec:local}
 As in the preceding paragraph, we consider a connected supergroup $G=(G_{\0},\g
g)$ and a globally oriented $G$-supermanifold $M$ with a weak $G$-invariant
superstructure.

Let $\a$ be an integrable equivariant form with smooth coefficients. In
particular this implies that $Supp(\a)$ is compact.

\subsubsection{$X$ and $T_NM(X)$} We fix $X\in\g g_0$ such that
$\overline{\exp(\mathbb R X)}\subset G_0$ is compact. We assume that $X$ is
central in $\g g$. Otherwise, we replace $\g g$ by $\g g(X)$, the centralizer of
$X$ in $\g g$, and $G$ by the connected component of  identity of the
centralizer $G(X)$ of $X$ in $G$.

Let $T_NM(X)=(TM|_{M(X)})/TM(X)\to M(X)$ be the normal bundle of $M(X)$ in $M$.
Its sheaf of sections is given for any open subset  $\UC\subset M_{\0}$ by:
\begin{equation}
\C_{T_NM(X)}(\UC\cap M(X)_{\0})=\Big(\C_{TM}(\UC)\tens_{\OC_{M(\UC)}}
\OC_{M(X)}(\UC\cap M(X)_{\0})\Big)/\C_{TM(X)}(\UC\cap M(X)_{\0}).
\end{equation}

\medskip

 Then $X$ determines a supervector bundle orientation of $T_NM(X)$  by the
following (cf. for example  \cite{BGV92} in the purely even case). 

Since $Q$ is $G$-invariant, $Q$ determines by restriction a weak Euclidean structure
on $TM(X)$ and then on $T_{N}M(X)$. We still denote the quotient weak Euclidean
structure by $Q$. Then $\LC^{\VC}(X)\in\C_{\g{osp}(T_{N}M(X))}(M(X)_{\0})$ is
invertible. Thus $(v,w)\mapsto Q(v,\LC^{\VC}(X)w)$ defines a symplectic
structure on $T_{N}M(X)$. This gives as usual an orientation of $T_{N}M(X)\to
M(X)$.

We recall that since $M$ is oriented, $TM\to M$ is oriented. We give to $M(X)$
the orientation such that the resulting  the quotient supervector bundle
orientation on $T_NM(X)=(TM|_{M(X)})/TM(X)\to M(X)$ is the one determined by
$X$.  

\bigskip

In order to ensure the existence of a Thom form on $T_{N}M(X)$, we assume that $T_NM(X)$ has an Euclidean structure (not just a weak Euclidean structure). For simplicity we assume and a $G$-invariant
superconnection $\mathbb A$ which preserves the Euclidean structure of $
T_NM(X)$. Moreover, we assume that $\mu(\mathbb A^{\OddId 2})_{\mathbb R}\le 0$
on $\Big(\what{M(X)}\mult_{M(X)}\Pi T_NM(X)\Big)_\0$ (condition $(*)$ with $\mu$
defined by  formula (\ref{eq:moment})  with $\AC=\OC(M(X))$, $V=\C_{\Pi
T_N(M(X)}(M(X)_{\0})$ and $B=\OddId Q$).

\subsubsection{A neighborhood of $M(X)$} Let $V$ be a $G_{\0}$-invariant open
neighborhood of $M(X)_{\0}$ in $M_{\0}$  such that $M(V)$ is isomorphic to $T_N
M(X)(V')$ where $V'$ is an open subset of $ T_N M(X)_{\0}$. Since $M_\0$ has a
$G_\0$-invariant weak Euclidean structure, such an open neighborhood exists. The
action of $G(X)$ on $V$ can be transferred to $V'$ by means of the isomorphism
between $M(V)$ and $M(X)(V').$ 

Moreover, we choose a closed neighborhood $V''$ of $M(X)_{\0}$ in $M_{\0}$ such
that $V''\subset V$. We denote by $U_a(\a,X,V'')$ the $G$-invariant open
neighborhood of $X$ in $\g g$ defined in (\ref{Ua}) for 
$a>0$ small enough. Since $\a$ is integrable,
$\a\in\what\W^\oo_{G,\int}(U_a(\a,X,V''),M)$.

\subsubsection{A limit formula} We put for $t\in \mathbb R$, $\phi_t(x)=\exp(\ii
tx)$. Let $\be$ be the form constructed in the preceding section (cf. formula (\ref{eq:beta})). Then $d_{\g
g}\be$ is even and we can define $\phi_t(d_{\g g}\be)$ by means of  the Taylor
formula (\ref{FonctLisse}).
 
 Let  $\psi_t$ be the analytic function such that
$\psi_t(x)={\frac{(1-\phi_t(x))}{x}}$ for $x\not= 0$ and $\psi_t(0)=\ii
t\exp'(0)=\ii t$. We put $\gamma_t=\be\psi_t(d_{\g g}\be)$. We have
$(1-\phi_t(d_{\g g}\be))\a=d_{\g g} (\gamma_t\a).$ Thus we have as functions in
$\CC^{\oo}(\g g)$:
\begin{equation}
\forall t\in\mathbb R,\quad \int_M\a(Z)=\int_M\phi_t(d_{\g g}\be(Z))\a(Z).
\end{equation}
(This is equivalent to say that this equality holds for any $Z\in\g g_{\PC}$
where $\PC$ is any near superalgebra.) In particular:
 \begin{equation}
\int_M\a(Z)=\lim\limits_{t\rightarrow +\oo}\int_M\phi_t(d_{\g g}\be(Z))\a(Z).
\end{equation}
Since $\phi_t(d_{\g g}\be)$ is smooth and bounded and $\a$ is integrable, these
integrals make sense. 
 
\subsubsection{A partition of unity} Let $\chi_1+\chi_2=1$ be a partition of
unity on $M$ such that $\chi_1$ is equal to $1$ on $V''$ and vanishes on the
complementary of $V$. Then we have:
\begin{equation}
\int_M\phi_t(d_{\g g}\be)\a=\int_{M(M_{\0}\setminus \overline
V'')}\chi_2\phi_t(d_{\g g}\be)\a+\int_{M(V)}\chi_1\phi_t(d_{\g g}\be)\a.
\end{equation} 
By construction of $U_a(\a,V'',X),$ the function $\Im m ( (\iota\be)_{\mathbb R})$ ($\Im m$ denoted the imaginary part of the complex)
has a strictly positive lower bound on
$U_a(\a,V'',X)\times Supp(\a)\cap\overline{M_{\0}\setminus V''}.$ Therefore
\begin{equation}
\lim_{t\rightarrow \oo}\int_{M(M_{\0}\setminus \overline V'')}\chi_2\phi_t(d_{\g
g}\be(Z))\a(Z)=0
\end{equation}

\subsection{The case of  isolated zeroes and smooth coefficients} Now we assume
that the zeroes of $X$ are isolated.

Let $p$ be such a zero. We can assume that $M(X)=\{p\}$ and that $V$ is a
neighborhood of $p$. Let $V'\subset (T_{p}M)_{\0}$ open such that $M(V)$ is
isomorphic to  $\VC'=T_{p}M(V')=V'\times (T_{p}(M))_{\1}$. We denote by $\t_{p}:\VC'\to M(V)$
this isomorphism. We transform  the action of $G$ on $M(V)$ to an action on
$\VC'$ by means of $\t_{p}$. 

We have to evaluate the limit of the integral of $\t_{p}^*(\chi_1\phi_t(d_{\g
g}\be(Z))\a(Z))$ on $\VC'.$ We denote by $\d_t$ the contraction of $T_pM$ by a 
factor of $\frac1{\sqrt t}$. We transform the action of $G$ on $\VC'$ into an
action on $\d_t(\VC')$ by means of $\d_t$. We have
$\mathop{\lim}\limits_{t\rightarrow\oo}\d_t(\VC')=T_pM$ and the action of $G$ which is obtained
in the limit is equal to the tangent action of $G$ on $T_pM$. 

Then the limit we are looking for is:
\begin{equation}
\lim_{t\rightarrow\oo}\int_{\d_t(V')}\d_t^*\t_{p}^*(\chi_1\phi_t(d_{\g
g}\be(Z))\a(Z)).
\end{equation}
We denote by $j_p$ the injection of the origin into $T_p(M).$ We have:
\begin{equation}
\begin{matrix}
\lim\limits_{t\rightarrow\oo}&\d_t^*\t_{p}^*(\chi_1)&=&1,\\ 
\lim\limits_{t\rightarrow\oo}&\d_t^*\t_{p}^*(\a)(Z)&=&j_p^*(\a)(Z).\\ 
\end{matrix}
\end{equation}
In the last equality, $j_p^*(\a)(Z)$ is seen as a constant function on  $T_pM.$
Finally $\lim\limits_{t\rightarrow\oo}t\d_t^*\t_{p}^*(Q)$ is the weak Euclidean
structure $Q_p$ on $T_pM$. It is $G$-invariant. 

For  $Z$ in $\g g$, we recall that $Z_{T_pM}$ denotes the vector field generated
by the linear action of $Z$ in $T_pM.$ Let $\be_{p}$ be the equivariant form on 
$T_pM$ defined by
 for $Z$ in $\g g$ by $\be_{p}(Z)=(\OddId Q_{p})^*(\OddId Z_{T_pM})$. We have:
\begin{equation}
\begin{matrix}
\lim\limits_{t\rightarrow\oo}&\d_t^*\t_{p}^*(\phi_t(d_{\g
g}\be(Z)))&=&\phi_1(d_{\g g}\be_{p}(Z)).\\ 
\end{matrix}
\end{equation}
Therefore,
\begin{equation}
\lim\limits_{t\rightarrow\oo}\int_{\d_t(V')}\d_t^*\t_{p}^*(\chi_1\phi_t(d_{\g
g}\be(Z))\a(Z))=j_p^*\a(Z)\int_{T_pM}\exp(-id_{\g g}\be_{p}(Z)).
\end{equation}

But the form under the integral on the right hand side is equivariantly closed
and integrable; so the evaluation of this integral has already been done in
Proposition
\ref{FormLin}. We obtained:

\begin{theo} Let $G=(G_{\0},\g g)$ be a supergroup 
 and $M$ be a globally oriented real $G$-supermanifold which has a $G$-invariant
Euclidean structure denoted by $Q$.

Let $\a\in\what\W_{G,\int}^{\oo}(\g g_{\0},M)$ such that $d_{\g g}\a=0$.
 
Let $X\in\g g_0$ such that $\overline{\exp(\mathbb R X)}\subset G_0$ is compact.
We assume that the zeroes of $X$ in $M$ are isolated.

Let  $p\in M(X)$. We denote by $\t_p$ the representation of $\g g(X)$ in $T_pM.$

We assume that $T_{p}M $ has a $G$-invariant weak Euclidean structure denoted by $Q'_{p}$ 

We assume that $\forall p\in M(X)$:
\begin{equation}
\big\{Z\in \g g_\0(X)\,\big/\,v\mapsto Q'_{p}(v,\t_p(Y)v)\text{ is positive
definite on $(T_pM)_\1$.}\big\}\not=\O
\end{equation}

Then, there exists an open subset $O$ of $\g g_{\0}(X)$ such that $X\in O$ and
as functions in $\CC^{\oo}_{\g g}(O)$:
\begin{equation}
\int_M\a(Z)=
\ii^{\frac{m-n}2}(2\pi)^{\frac{n+m}2}\som_{p\in M(X)}\frac{j_p^*(\a)(Z)}
{Spf(\t_{p}(Z))}.
\end{equation}

\end{theo}
\begin{proof}[Proof] As the fixed points are isolated, $\t_p(Z)$ is invertible
in $\g{gl}(T_{p}M)$ on a neighborhood of $X$. We denote by $\OC(p)$ the open
subset of $\g g_{\0}(X)$ defined by:
\begin{equation}
\OC(p)=\big\{Z\in\g g_0(X)\,\big/\,\t_p(Z)\text{ is invertible  in
}\g{gl}(T_{p}M)\big\}.
\end{equation}

Let us consider $O_p=U_a(\a,X,V'')\cap\OC(p)$.  On $O_p,$ $Spf\circ\t_{p}$ is
smooth and invertible. As $j_p^*(\exp(-id_{\g g}\be_{p}(Z)))=1$ and as the
hypotheses imply that $m$ is even, we obtain the formula considering
\begin{equation}
O=\bigcap\limits_{p\in M(X)_{\0}\cap Supp(\a)}O(p)
\end{equation}
($M(X)_{\0}\cap Supp(\a)$ is finite).
\end{proof}

\subsection{Localization formula for non-isolated zeroes} In
 section \ref{preLoc} we restricted  the problem to an integration in a
$G$-invariant neighborhood of $M(X)$ in $M$. But such a neighborhood is
isomorphic to an open neighborhood of $M(X)$ in $T_N M(X)$ on which we can
transfer the action of $G$. As in the case of  isolated zeroes we can restrict
the problem to an integration on $T_NM(X)$. Here, instead of applying 
Proposition
\ref{FormLin} we apply Proposition \ref{FormFib}.

Let $\a$ be an integrable equivariant form with smooth coefficients. Let $C$ be
a relatively compact $G_{\0}$-invariant open neighborhood of the support of $\a$
in $M_{\0}$. As the support of $\a$ is $G_{\0}$-invariant and as $M_{\0}$ has a
$G_{\0}$-invariant Rimanian structure, such a neighborhood exists. 

In order to apply Proposition
\ref{FormFib} we have to check that $X$ is in the open subsupermanifold
$U^{T_N{M(X)}}\big(M(X)_{\0}\cap C\big)$ (cf. formula (\ref{U}) for a
definition). Indeed, since $M(X)$ is the manifold of zeroes of $X$, the action
of $X$ on the fibres of $T_NM(X)$ is invertible.

Finally, we have:

\begin{theo}\label{Formula} Let $G=(G_{\0},\g g)$ be a supergroup and $M$ be a
globally oriented real $G$-supermanifold with  a weak Euclidean structure.

Let $\a\in\what\W_{G,\int}^{\oo}(\g g_{\0},M)$ such that $d_{\g g}\a=0$.

Let $X\in\g g_\0$ such that $\overline{\exp(\mathbb R X)}\subset G_\0$ is
compact.

We denote by $T_N M(X)$ the  normal bundle to the manifold of zeroes of $X$ in
$M$. We assume $T_N M(X)$, has a $G(X)$-invariant Euclidean structure $Q'$ and  a
$G(X)$-invariant superconnection, $\mathbb A$ which preserves the Euclidean
structure of $T_N M(X)$.

Moreover we assume:
\begin{enumerate}
\item[$(*)$] $\mu(\mathbb A^{\OddId 2})_{\mathbb R}\le 0$ on $\Big(\Pi
T_NM(X)\mult_{M(X)}\what{M}(X)\Big)_\0$ (cf. formula (\ref{eq:moment}) for the
definition of $\mu$ with $\AC=\OC(\what M(X))$, $V=\C_{\Pi
T_NM(X)\mult_{M(X)}\what{M}(X)}(\what M(X)_{\0})$ and $B=\OddId Q'$);
\item[$(**)$] there is a covering of $M$ by open subsupermanifolds $W$   such
that $U_{+}^{T_NM(X)_{\1}}(W)$ contains a non empty open subset (cf. formula
(\ref{U0m}) for the definition of $U_{+}^{T_NM(X)_{\1}}(W)$ with $\UC=W$ and
$\VC=T_NM(X)$).
\end{enumerate}

We denote by $j$ the canonical injection of $M(X)$ into $M$ and by $\EC_{\g g}$
a representative  of the $G(X)$-equivariant Euler class of the normal bundle
$T_N M(X)$.

Then, there exists an open subset $O$ of $\g g(X)_{\0}$ such that $X\in O$ and
as a function on $\g g(X)(O)$ we have:
\begin{equation}
\int_M\a(Z)=(2\pi)^{\frac{n+m}2}\int_{M(X)}\frac{j^*\a(Z)}{
\EC_{\g g}(Z)}.
\end{equation}
(This is equivalent to say that this equality holds for any $Z\in\g
g(X)_{\PC}(O)$ where $\PC$ is any near superalgebra.)
\end{theo}

\begin{proof}[Proof] As in the case of isolated zeroes the proof is similar to
that of the non-super case (cf. for example
\cite{BGV92}).
\end{proof}

 We can note that the integral on the right hand side is a sum of integrals on
the connected components of $M(X)$. As the support of $\a$ is compact the number
of such integrals that are not equal to zero is finite.

\end{document}